\documentclass[aop, preprint, noinfoline]{imsart}

\setattribute{journal}{name}{}

\usepackage[OT1]{fontenc}	
\usepackage{graphicx}
\usepackage{amsmath, amsthm, amssymb}
\usepackage{mathrsfs}
\usepackage{stmaryrd}			
		\SetSymbolFont{stmry}{bold}{U}{stmry}{m}{n}

\usepackage{pdfsync}
\usepackage[a4paper,top=3cm,bottom=3cm,left=2.5cm,right=2.5cm, bindingoffset=0mm]{geometry}

\usepackage[usenames,dvipsnames]{xcolor}
\usepackage{enumerate}								
\usepackage[sort,numbers]{natbib}								
\usepackage[hyphens]{url}									
\usepackage{verbatim}								

\usepackage{dsfont}									

\usepackage[pdftex,bookmarks,
						colorlinks,				
						breaklinks,unicode,
						citecolor=OliveGreen,
						linkcolor=Maroon
						]{hyperref} 

\usepackage{mathtools}

\usepackage{tikz}
\usetikzlibrary{calc}

\arxiv{}

\startlocaldefs

\newcommand{\mfe}{{\mathfrak e}}
\newcommand{\mff}{{\mathfrak f}}
\newcommand{\mfg}{{\mathfrak g}}
\newcommand{\mfh}{{\mathfrak h}}
\newcommand{\mfi}{{\mathfrak i}}

\newcommand{\mfl}{{\mathfrak l}}

\newcommand{\mfs}{{\mathfrak s}}

\newcommand{\mfz}{{\mathfrak z}}

\newcommand{\mfB}{{\mathfrak B}}
\newcommand{\mfC}{{\mathfrak C}}

\newcommand{\mfP}{{\mathfrak P}}

\newcommand{\mfS}{{\mathfrak S}}

\newcommand{\mcC}{{\mathcal C}}
\newcommand{\mcD}{{\mathcal D}}

\newcommand{\mcF}{{\mathcal F}}
\newcommand{\mcG}{{\mathcal G}}

\newcommand{\mcK}{{\mathcal K}}
\newcommand{\mcL}{{\mathcal L}}

\newcommand{\mcO}{{\mathcal O}}
\newcommand{\mcP}{{\mathcal P}}

\newcommand{\mcR}{{\mathcal R}}
\newcommand{\mcS}{{\mathcal S}}

\newcommand{\mcU}{{\mathcal U}}

\newcommand{\msB}{{\mathscr B}}
\newcommand{\msC}{{\mathscr C}}
\newcommand{\msD}{{\mathscr D}}

\newcommand{\msF}{{\mathscr F}}

\newcommand{\msH}{{\mathscr H}}			
\newcommand{\msK}{{\mathscr K}}
\newcommand{\msP}{{\mathscr P}}

\newcommand{\msR}{{\mathscr R}}

\newcommand{\mbfL}{{\mathbf L}}

\newcommand{\mbfW}{{\mathbf W}}
\newcommand{\mbfX}{{\mathbf X}}
\newcommand{\mbfY}{{\mathbf Y}}

\newcommand{\mbfb}{{\mathbf b}}
\newcommand{\mbfc}{{\mathbf c}}

\newcommand{\mbfe}{{\mathbf e}}
\newcommand{\mbff}{{\mathbf f}}

\newcommand{\mbfj}{{\mathbf j}}
\newcommand{\mbfk}{{\mathbf k}}

\newcommand{\mbfm}{{\mathbf m}}
\newcommand{\mbfn}{{\mathbf n}}

\newcommand{\mbfs}{{\mathbf s}}
\newcommand{\mbft}{{\mathbf t}}
\newcommand{\mbfu}{{\mathbf u}}

\newcommand{\mbfw}{{\mathbf w}}
\newcommand{\mbfx}{{\mathbf x}}
\newcommand{\mbfy}{{\mathbf y}}
\newcommand{\mbfz}{{\mathbf z}}

\newcommand{\boldalpha}{{\boldsymbol \alpha}}

\newcommand{\boldeta}{{\boldsymbol \eta}}

\newcommand{\boldlambda}{{\boldsymbol \lambda}}

\newcommand{\mbbE}{{\mathbb E}}

\newcommand{\eps}{\varepsilon}
\newcommand{\kai}{\varkappa}

\newcommand{\defeq}{\eqqcolon}

\newcommand{\mathsc}[1]{\text{\textsc{#1}}}

\newcommand{\emparg}{{\,\cdot\,}}

\DeclareMathOperator{\Span}{span}
\DeclareMathOperator{\Homeo}{Homeo}
\DeclareMathOperator{\Diffeo}{Diff}

\DeclareMathOperator{\Dir}{Dir}								
\DeclareMathOperator{\conv}{conv}

\DeclareMathOperator{\law}{law}

\DeclareMathOperator{\ev}{ev}

\DeclareMathOperator{\pr}{pr}

\DeclareMathOperator{\diag}{diag}

\DeclareMathOperator{\eqdef}{\coloneqq}

\let\epsilon\varepsilon
\let\subset\subseteq

\newcommand{\rar}{\rightarrow}

\newcommand{\nlim}{\lim_{n}}									

\newcommand{\hlim}{\lim_{h }}

\DeclareMathOperator*{\dirlim}{\underrightarrow{\lim}}							

\DeclareMathOperator{\length}{length}							
\DeclareMathOperator{\supp}{supp}								
\newcommand{\diff}{\mathop{}\!\mathrm{d}}						
	
\newcommand{\abs}[1]{\left\lvert#1\right\rvert}						
\newcommand{\norm}[1]{\left\lVert#1\right\rVert}					
\newcommand{\set}[1]{\left\{#1\right\}}							
\newcommand{\tonde}[1]{\left(#1\right)}							
\newcommand{\quadre}[1]{\left[#1\right]}							
\newcommand{\scalar}[2]{\left\langle #1 \,\middle |\, #2\right\rangle}		

\newcommand{\hodot}{\overline{\odot}}
\newcommand{\hoplus}{\overline{\oplus}}

\DeclareMathOperator{\cl}{cl}									
\DeclareMathOperator{\diam}{diam}								
\newcommand{\seq}[1]{\tonde{#1}}								

\DeclareMathOperator{\Meas}{\mathscr M}
\DeclareMathOperator{\Mp}{\mathscr M^+}
\newcommand{\Mb}{\mathscr M_b}
\newcommand{\Mbp}{\mathscr M_b^+}

\newcommand{\TV}{{\mathsc{tv}}}

\newcommand{\pfwd}{\sharp}

\newcommand{\quotient}[2]{\left.\raisebox{.1em}{$#1\!$}\middle/\raisebox{-.1em}{$#2$}\right.}

\DeclareMathOperator{\car}{\mathds 1}

\DeclareMathOperator{\emp}{\varnothing} 

\DeclareMathOperator{\N}{{\mathbb N}}
\DeclareMathOperator{\R}{{\mathbb R}}

\DeclareMathOperator{\Q}{{\mathbb Q}}
\renewcommand{\C}{{\mathbb C}} 	
\DeclareMathOperator{\Z}{{\mathbb Z}}

\newcommand{\K}{\Bbbk}

\newcommand{\restr}{\big\lvert}

\newcommand{\functionnn}[5]{\begin{align*}#1\colon#2&\longrightarrow#3\\ #4&\longmapsto#5\end{align*}}		

\usepackage{booktabs}

\allowdisplaybreaks

\usetikzlibrary{shapes.misc}
\tikzset{cross/.style={cross out, draw=black, minimum size=2*(#1-\pgflinewidth), inner sep=0pt, outer sep=0pt},
cross/.default={4pt}}

\usepackage{youngtab}

\newcommand\catalannumber[4]{
  (#1)
  \foreach \dir in {#4}{
    \ifnum\dir=0
    -- ++(1,0)
    \else
    -- ++(0,1)
    \fi
  } |- (#1);
  \draw[help lines] (#1) grid +(#2,#3);
  \coordinate (prev) at (#1);
  \foreach \dir in {#4}{
    \ifnum\dir=0
    \coordinate (dep) at (1,0);
    \else
    \coordinate (dep) at (0,1);
    \fi
    \draw[line width=1pt] (prev) -- ++(dep) coordinate (prev);
  };

}

\newcommand\independent{\protect\mathpalette{\protect\independenT}{\perp}}
\def\independenT#1#2{\mathrel{\rlap{$#1#2\,\,$}\mkern2mu{#1\,#2}}}

\newcommand{\comma}{\quad\textrm{,}\quad}
\newcommand{\semicolon}{\quad\textrm{;}\qquad}
\newcommand{\fstop}{\quad\textrm{.}}

\newcommand{\MultiZero}[1]{M_1\!\quadre{#1}}
\newcommand{\Mnom}[1]{M_2\!\quadre{#1}}
\newcommand{\Bella}[1]{M_3\!\quadre{#1}}

\DeclareMathOperator{\zero}{{\mathbf 0}}
\DeclareMathOperator{\uno}{{\mathbf 1}}

\newcommand{\av}[1]{\left\langle#1\right\rangle}

\DeclareMathOperator{\Gam}{Gam}
\DeclareMathOperator{\Poi}{Poi}

\DeclareMathOperator{\Isom}{Isom}

\DeclareMathOperator{\EGF}{EGF}

\newcommand{\ext}{{\mathrm{ext}}}

\newcommand{\n}[1]{{\overline{#1}}}

\newcommand{\compo}{\diamond}
\newcommand{\contra}{\star}

\newcommand{\Poch}[2]{\left\langle{#1}\right\rangle_{#2}}

\newcommand{\mtbinom}[2]{\big(\!\!\tbinom{#1}{#2}\!\!\big)}

\renewcommand{\length}[1]{{\abs{#1}}}
\newcommand{\bag}[1]{\llbracket{#1}\rrbracket}

\renewcommand{\P}{{\mathbb P}}

\let\temp\phi
\let\phi\varphi
\let\varphi\temp

\newcommand{\Beta}{\mathrm{B}}

\newcommand{\DF}{{\mcD}}
\newcommand{\GP}{{\mcG}}
\newcommand{\PP}{{\mcP}}

\numberwithin{equation}{section}
\theoremstyle{plain}
\newtheorem{thm}{Theorem}[section]
\newtheorem{prop}[thm]{Proposition}
\newtheorem{lem}[thm]{Lemma}

\theoremstyle{definition}
\newtheorem{conj}[thm]{Conjecture}

\theoremstyle{remark}
\newtheorem{rem}[thm]{Remark}
\newtheorem{ese}[thm]{Example}
\endlocaldefs

\begin{document}

\begin{frontmatter}
\title{Moments and cumulants in infinite dimensions\\ with applications to\\ Poisson, Gamma and Dirichlet--Ferguson\\ random measures\thanksref{T1}
}
\runtitle{Moments \& cumulants in infinite dimensions}

\begin{aug}
\author{\fnms{Lorenzo} \snm{Dello Schiavo}\thanksref{t1}\ead[label=e1]{delloschiavo@iam.uni-bonn.de}}

\thankstext{T1}{This work has been jointly supported by the CRC 1060 \emph{The Mathematics of Emergent Effects} at the University of Bonn, funded through the Deutsche Forschungsgemeinschaft, and by the Hausdorff Center for Mathematics at the University of Bonn, funded through the German Excellence Initiative.}

\thankstext{t1}{I~gratefully acknowledge fruitful conversations with P.~Ferrari; \mbox{K.-T.}~Sturm and M.~Huesmann;~E. Lytvynov, M. Gubinelli and L.~Borasi; S.~Albeverio, M.~Gordina and C.~Stroppel, on some aspects of this work respectively concerning combinatorics, probability theory, infinite-dimen\-sion\-al analysis and representation theory. I~am also indebted to A.~Nota for having introduced me to the concept of cumulants, to E.~Lytvynov and W.~Miller for pointing out to me the references~\cite{KonLyt00} and~\cite{Mil77} respectively, and finally to M.~Huesmann for partly proofreading a preliminary version of this work.
}
\runauthor{L. Dello Schiavo}

\affiliation{University of Bonn
}

\address{Institut f\"ur Angewandte Mathematik\\
Rheinische Friedrich-Wilhelms-Universit\"at Bonn\\
Endenicher Allee 60\\
DE 53115 Bonn\\
Germany\\
\printead{e1}
}

\end{aug}

\begin{abstract}
We show that the chaos representation of some Compound Poisson Type processes displays an underlying intrinsic combinatorial structure, partly independent of the chosen process.
From the computational viewpoint, we solve the arising combinatorial complexity by means of the moments/cumulants duality for the laws of the corresponding processes, themselves measures on  distributional spaces, and provide a combinatorial interpretation of the associated `extended' Fock spaces.
From the theoretical viewpoint, in the case of the Gamma measure, we trace back such complexity to its `simplicial part', i.e. the Dirichlet--Ferguson measure, hence to the Dirichlet distribution on the finite-dimensional simplex. We thoroughly explore the combinatorial and algebraic properties of the latter distribution, arising in connection with cycle index polynomials of symmetric groups and dynamical symmetry algebras of confluent Lauricella functions.
\end{abstract}

\vspace{.5cm}
\today
\vspace{.5cm}



\end{frontmatter}

\section{Introduction}
Given a probability distribution~$\nu$ on a linear space~$\Phi$, denote by 
\begin{align*}
\mcL[\nu](tX)\eqdef& \int_{\Phi} \exp[ t X(x)] \diff \nu(x)\quad, & \textrm{resp. by} && \mcK[\nu](tX)\eqdef& \ln \mcL[\nu](tX) \comma
\end{align*}
its {moment}, resp. {cumulant}, {generating function} computed at a real-valued random variable~$X$. The coefficients in the corresponding McLaurin expansions in $t$, respectively the (\emph{raw}) \emph{moments}~$\mu'_n=\mu^{\prime\, \nu}_n[X]$ and the \emph{cumulants}~$\kappa_n=\kappa_n^\nu[X]$, represent a cornerstone in elementary probability, with both broad-ranging applications and great interest per se (e.g. --~for the cumulants~-- in \emph{free probability}~\mbox{\cite[and ref.s therein]{AriHasLehVar15, Leh04, Leh06}} and \emph{kinetic derivation theory}~\cite{LukMarNot16,LukMar16}; also, see~\cite{Hal00} for historical remarks).
Among the inherent problems, that of obtaining moments from cumulants is seemingly a simple one, usually reduced to the observation that
\begin{align*}
\mu'_n=B_n[\kappa_1,\dotsc, \kappa_n] \comma
\end{align*}
where $B_n$ denotes the $n$-variate (complete) Bell polynomial. Similar formulae hold for obtaining cumulants from moments, raw moments from central moments and vice versa, thus suggesting for these correspondences the --~somewhat informal~-- designation of \emph{dualities}.

In this article we show the following:
\begin{itemize}
\item appropriate generalizations of such dualities~(Thm.~\ref{t:MomCum}) provide a better understanding (e.g.~Rem.s~\ref{r:MomPoi}\&\ref{r:ExtFock}) and simpler proofs (e.g.~Thm.~\ref{t:Surgailis}, Prop.~\ref{p:Recursive}) of \emph{multiple stochastic integration} (MSI) and related results for the (compensated) Poisson measure~$\PP$, the Gamma measure~$\GP$ and the laws of other \emph{compound Poisson type} (CPT) processes;

\item the intrinsic \emph{exponential} nature of CPT laws (\S\S\ref{ss:LK}\&\ref{ss:CPT}) results in an expression for their cumulants (eq.~\eqref{eq:PoiCumGen}) in terms of the moments of the associated intensity and L\'evy measures and more generally prompts to cumulants as preeminent over moments for computational purposes;

\item if the moments of the CPT laws' L\'evy measures have a peculiar factorial part -- we discuss the case of~$\GP$ as a prototypical example~-- the combinatorial complexity of MSI, pertinent to set partition lattices, is reduced to that of permutations (\S\ref{s:Bell}, Rem.~\ref{r:CQES}, Prop.~\ref{p:MomGP});

\item in the case of~$\GP$, such combinatorial complexity is traced back to the simplicial part of the law, namely the Dirichlet--Ferguson measure~$\DF$, whence to the \emph{Dirichlet distribution} on the finite-dimensional standard simplex. We show how a detailed study of the latter cannot prescind from combinatorics and algebra and yields interesting connections with enumeration theory (\S\ref{ss:Polya}) and Lie representation theory (\S\S\ref{ss:Lauricella}\&\ref{app:A}), which we thoroughly explore in the finite-dimensional framework.
\end{itemize}

Throughout the paper we focus on~$\PP$, $\GP$ and $\DF$ as three major examples. Most of the results concerned with~$\DF$ are proved as generalizations of their analogues for the Dirichlet distribution in finite dimensions. 
Where such generalizations are not yet at hand, suitable conjectural statements and heuristics are provided (\S\ref{ss:UBDSA}), on which we plan to base future work along the line sketched in~\S\ref{s:Conclusions}.
\smallskip

\paragraph{Motivations} Our study is partly motivated by that of chaos representation~(\S\ref{ss:Chaos}) for square-integrable functionals of~CPT random measures, which provide a large class of interesting examples. As noticed in Kondratiev~et al.~\cite[4.5]{KonDaSStrUs98} for the \emph{Gamma measure}~(\S\ref{ss:PGDFrms}), CPT processes do not posses the \emph{chaos representation property} (CRP). As a consequence, decompositions of the space of square-integrable functionals obtained via generalized Fourier transforms become in this case particularly involved (cf. e.g.~\cite{KonDaSStrUs98, Lyt03, BarOueRia11}).
A useful tool in the `simpler' framework of CRP for processes of Poisson and Gaussian type is {Engel--Rota--Wallstrom MSI theory}, based on the combinatorial properties of {set partitions lattices} (see the standard reference monograph {Peccati--Taqqu}~\cite{PecTaq11}). Despite the great generality achieved and some strikingly simple formulations of otherwise highly non-trivial expressions for moments and cumulants of multiple stochastic integrals (e.g. Rota--Wallstrom~\cite{RotWal97}), such theory heavily relies on set partitions equations, usually hardly tractable in explicit non-recursive forms.

Partly because of such complexity, the chaos representation for CPT processes has been variously addressed (cf. e.g.~\cite{DiNOksPro04, Lyt03}) by MSI with respect to power jump processes (i.e. Nualart--Schoutens chaos representation~\cite{NuaSch00}), via orthogonalization techniques for polynomials on distributions spaces~\cite{Lyt03}, through vector valued Gaussian white-noise (Tsilevich--Vershik~\cite{TsiVer03}), or by MSI with respect to Poisson processes, which is mostly the case of our interest. In that framework, unitary isomorphisms are realized, via Jacobi fields, to Fock-type spaces. Whereas the combinatorial properties of such Jacobi fields are often intractable and usually settled by means of recursive descriptions (e.g.~\cite{BerLytMie03, KonLyt00, Lyt03, BarOueRia11}), the said \emph{extended/non-standard Fock spaces}~\cite{KonLyt00, KonDaSStrUs98, Lyt03,BarOueRia11, BerMie01} associated to CPT processes display an intrinsic structure, addressed in full generality in Lytvynov~\cite{Lyt03}, and connected to (the enumeration of) set partitions, thus hinting to Engel--Rota--Wallstrom theory. A basic example of this interplay between MSI and set partitions (lattices) in the CPT case is the `minimality' of the embedding~\cite{KonDaSStrUs98} of the standard Fock space associated to~$\PP$ into the extended Fock space of~$\GP$ (see Rem.~\ref{r:ExtFock}).

\smallskip

Set partitions enumeration appearing as a tool in the aforementioned approach to CRP turns out to be an incarnation of moments/cumulants dualities for the laws of the considered CPT processes, in which setting we show that cumulants (rather than moments) are the interesting quantities. We base this observation, of rather physical flavor, on the study of Gibbs measures, a {Leitmotiv} in heuristic derivations in probability from Feynman's derivation of the Wiener measure to that of the entropic measure in Wasserstein diffusion (von Renesse--Sturm~\cite{vReStu09}). Indeed, the \emph{Helmholtz free energy} (see Rem.~\ref{r:MeccaStat}) associated to a Gibbs measure is (up to sign and normalization) the cumulant generating function of the energy of the system, and the \emph{additivity} of the latter translates into a classical feature of cumulants: \emph{linearity} on independent random variables.
This property becomes of particular importance when the {partition function} of the system, here playing as a moment generating function, may be expressed --~in a `natural' way~-- in exponential form. Physical heuristics aside, in the CPT framework the required exponential form is clear, as it is provided by the well-known \emph{L\'evy--Khintchine formula}~\eqref{eq:LK}; hence, for example, the cumulant generating function of a~CPT random measure with no drift nor Brownian part is \emph{linear} in the L\'{e}vy measure. The computation of cumulants becomes then essentially trivial and moments are recovered via duality.

\smallskip

\paragraph{Results} Among other CPT processes, we compare the law~$\PP$ of the Poisson process with that,~$\GP$, of the Gamma process, showing how set partitions combinatorics entailed in the moments/\allowbreak cum\-ulants duality may be significantly simplified, thus making the set partitions lattices machinery unnecessary. Indeed, while the moments of~$\PP$ are given by Bell polynomials (computed at monomial powers, see Thm.~\ref{t:ConvLogPP} \& Rem.~\ref{r:MomPoi}) --~hence, naturally, by averages over \emph{set} partitions~--, those of~$\GP$ are given by augmented cycle index polynomials of symmetric groups (Prop.~\ref{p:MomGP}) --~related in the same way to \emph{integer} partitions.

On the one hand, provided this understanding of moments, we show how convoluted (recursive) expressions for the scalar product of the extended Fock space of~$\GP$ (and of other non-standard Fock spaces) follow immediately (cf. e.g. Prop.~\ref{p:Recursive}) from classical identities 
for Bell polynomials.
On the other hand, the correspondence between integer partitions and the cyclic structure of permutations (\S\ref{s:Bell}) suggests that the arising structure is purely a property of the `simplicial part' (cf.~\cite{TsiVerYor00}) of~$\GP$, i.e. of the \emph{Dirichlet--Ferguson measure}~$\DF$~(see~\cite{Fer73},~\S\ref{ss:PGDFrms}), supported on a space of probability measures.

\smallskip

Alongside with its finite-dimensional analogue, i.e. the \emph{Dirichlet distribution}~(see~\cite{NgTiaTan11},~\S\ref{ss:PoiGamDir})
, the Dirichlet--Ferguson measure~$\DF$ has been widely studied both per se (e.g. Lijoi--Regazzini~\cite{LijReg04}) and in connection with the Gamma process (Tsilevich--Vershik--Yor~\cite{TsiVerYor00,TsiVerYor01}), mostly owing to its numerous applications, ranging from population genetics to coalescent theory, number theory and Bayesian nonparametrics (see e.g. the monograph~Feng~\cite{Fen10} and the surveys Lijoi--Pr\"unster~\cite{LijPru10} and Berestycki~\cite{Ber09}).
In this work, we focus on the combinatorial and algebraic properties of the Dirichlet distribution, which we subsequently extend to the case of~$\DF$.

\smallskip

On the combinatorial side, deepening results by Kerov--Tsilevich~\cite{KerTsi01}, we compute the moments of the Dirichlet distribution (Thm.~\ref{l:lemma}) as cycle index polynomials of symmetric groups and interpret this result in light of \emph{P\'{o}lya enumeration theory}~(Prop.~\ref{p:AggrPolya}). As a byproduct, we obtain  a new representation (Prop.~\ref{p:LaurEGF}) and asymptotic formulae (Prop.~\ref{l:Asympt}) for the moment generating function of the distribution, i.e. the \emph{second multivariate Humbert function}~${}_k\Phi_2$, a confluent form of the Lauricella function~${}_kF_D$; we subsequently generalize this limiting behavior to the Dirichlet--Ferguson measure, matching on arbitrary compact separable spaces a result in~\cite{vReStu09} for the entropic measure on the 1-sphere (Prop.~\ref{p:AsymptDF}). A key r\^ole in this combinatorial perspective is assumed by the actions of symmetric groups on supports of Dirichlet distributions (i.e. the standard simplices), for which we provide different complementary interpretations (Rem.~\ref{r:Models}, Prop.~\ref{p:chased}).

On the algebraic side, a metaphor of urns and beads (\S\ref{ss:Polya}) points to the \emph{dynamical symmetry algebra} of the second Humbert function. By means of the general theory in Miller~\cite[and ref.s therein]{Mil73b} we show (Thm.~\ref{t:dsa}) that the said algebra is the special linear Lie algebra of square matrices with vanishing trace. The result is then --~partly heuristically~-- extended by means of Lie theory to the infinite-dimensional case of~$\DF$ (Thm.~\ref{t:Borel0}, Conj.~\ref{t:Borel1}); the construction is reminiscent of Vershik's construction of the \emph{infinite-dimensional Lebesgue measure}~\cite{Ver07}.
We point out that this connection with Lie theory arising from the study of Dirichlet(--Ferguson) measures is not entirely surprising: For instance, several results concerned with the Gamma measure (e.g.~Kondratiev--Lytvynov--Vershik~\cite[and ref.s therein]{KonLytVer15}) were obtained as part of a much wider program~(cf. e.g.~\cite[\S1.4]{TsiVer03}) to study infinite-dimensional representations of (measurable) $\mfs\mfl_2(\R)$-current groups.

\smallskip

We postpone some final remarks to~\S\ref{s:Conclusions}.

\paragraph{Summary} 

Firstly, basic facts on partitions and permutations, exponential generating functions and multi-sets are reviewed~(\S\ref{s:Prelim}). 
Secondly, some notions are recalled on infinitely divisible distributions and the L\'{e}vy--Khintchine formula for measures on finite-dimensional linear spaces~(\S\ref{ss:LK}) and on cylindrical measures on nuclear spaces~(\S\ref{ss:Radon}), turning then to dualities between moments and cumulants for measures on such spaces (\S~\ref{ss:MomCum}) and their applications to CPT processes.
The main results are contained in~\S\S\ref{s:Dir}\&\ref{s:DirFer}, where we comparatively explore some properties of the Poisson, Gamma and Dirichlet distributions, respectively of the laws~$\PP$,~$\GP$ and~$\DF$ of Poisson, Gamma and Dirichlet--Ferguson processes, mostly focusing on the Dirichlet case.
The combinatorial properties of (extended) Fock spaces are addressed in~\S\ref{ss:Fock}. Some conclusive remarks are found in~\S\ref{s:Conclusions}. Finally, complementary results and proofs, mainly concerned with Lie theory, are collected in the Appendix~\S\ref{app:A}.

\section{Combinatorial preliminaries}\label{s:Prelim}
\paragraph{Notation} Let $n$ be a positive integer (in the following usually implicit) and set\begin{align*}
\mbfx\eqdef& \seq{x_1,\dotsc, x_n} & \mbfe^i\eqdef&\seq{{}_10,\dotsc,0,{}_i1,0,\dotsc,{}_n0}\\
\uno\eqdef& \seq{{}_11,\dotsc,{}_n1} & {}_i \mbfx_{j}\eqdef& \seq{x_{i+1}, x_{i+2},\dotsc, x_{j-1}, x_j}\\
\length{\mbfx}\eqdef& x_1+\cdots+x_n
 &\mbfx_{\hat \imath}\eqdef& \seq{x_1,\dotsc, x_{i-1}, x_{i+1}, \dotsc, x_n}\\
\vec{\mbfn}\eqdef& \seq{1,\dotsc,n} & \mbfx \oplus \mbfy\eqdef&\seq{x_1,\dotsc,x_n,y_1,\dotsc, y_n}\\
[n]\eqdef& \set{1,\dotsc, n} & \pi\in\mfS_n\eqdef& \set{\textrm{permutations of }[n]}\\
\mbfx_\pi\eqdef& \seq{x_{\pi(1)}, \dotsc, x_{\pi(n)}} & \mbfx\compo\mbfy\eqdef& (x_1y_1,\dotsc, x_ny_n)\\
\mbfx^{\compo n}\eqdef& \underbrace{\mbfx\compo\dots\compo \mbfx}_{n \textrm{ times}} & \mbfx\cdot\mbfy\eqdef & x_1y_1+\cdots + x_ny_n \comma
\end{align*}
where we stress the position of an element in a vector with a \emph{left} subscript. We stress that~$\length{\mbfx}$ \textbf{is always a signed quantity} and the symbol~$\length{\mbfx}$ \textbf{never denotes a norm of the vector~$\mbfx$}.

For $f\colon \C\rar \C$ or $\mbff\eqdef (f_1,\dotsc, f_n)$ with $f_i\colon \C\rar\C$, write
\begin{align*}
f^{\compo}[\mbfx]\eqdef& \seq{f[x_1],\dotsc,f[x_n]} & \mbff^{\compo}[\mbfx]\eqdef& \seq{f_1[x_1],\dotsc,f_n[x_n]} \\
f[\mbfx]\eqdef& (f^{\compo}[\mbfx])^{\uno} & \mbff[\mbfx]\eqdef& (\mbff^{\compo}[\mbfx])^{\uno} \comma
\end{align*}
e.g. $\mbfx^{\compo\mbfy}=\seq{x_1^{y_1},\dotsc, x_n^{y_n}}$ while $\mbfx^\mbfy= x_1^{y_1}\dots x_n^{y_n}$.
Finally, denote by $\Gamma$ the \emph{Euler Gamma function}, by $\Poch{\alpha}{k}\eqdef \Gamma[\alpha+k]/\Gamma[k]$ the \emph{Pochhammer symbol} of~$\alpha>0$, by
\begin{align*}
\Beta[x,y]\eqdef \frac{\Gamma[x]\Gamma[y]}{\Gamma[x+y]} && \textrm{resp.} && \Beta[\mbfx]\eqdef& \frac{\Gamma[\mbfx]}{\Gamma[\abs{\mbfx}]}
\end{align*}
the \emph{Euler Beta function}, resp. its multivariate analogue. We reserve the upright typeset `$\Beta$' for the Beta function, so that no confusion with Bell polynomials or Bell numbers (see below) may arise.

\subsection{Partitions and permutations combinatorics}
\emph{We briefly review some basic facts about partitions, permutations and multi-sets. Concerning set partitions, resp. Bell and Touchard polynomials, a more exhaustive exposition may be found in~\cite[\S2.3, resp.~\S2.4]{PecTaq11}; however, we rather put the emphasis on the comparison among enumeration of set partitions, integer partitions and permutations. Further related results in enumerative combinatorics may be found in~\cite{Sta01}}.

\paragraph{Set and integer partitions}\label{s:Bell}  A \emph{set partition} of $[n]$ is a tuple~$\mbfL\eqdef\seq{L_1,\dotsc, L_r}$ of disjoint subsets~$\emp\subsetneq L_i\subset [n]$, termed \emph{clusters} or blocks, and such that $\sqcup_i L_i=[n]$. For any such partition write~$\mbfL\vdash [n]$ and $\mbfL \vdash_r [n]$ if $\length{\mbfL}=r$, i.e. if $\mbfL$ has $r$ clusters.

A (\emph{integer}) \emph{partition~$\boldlambda$ of~$n$ into~$r$ parts} (write: $\boldlambda\vdash_r n$) is a vector of \emph{non-negative integer} solutions of the system $\abs{\vec\mbfn\compo\boldlambda}=n$, $\abs{\boldlambda}=r$. Term $\boldlambda$ a (\emph{integer}) \emph{partition of $n$} (write:~$\boldlambda\vdash n$) if the second requirement is dropped. We always regard a partition in its \emph{frequency representation}, i.e. as the vector of its ordered frequencies (see e.g.~\cite[\S1.1]{And76}).

To a set partition $\mbfL\vdash_r [n]$ one can associate in a unique way the partition $\boldlambda(\mbfL)\vdash_r n$ by setting $\lambda_i(\mbfL)\eqdef \#\set{h\in [r]\mid \#L_h=i}$.
The number of set partitions $\mbfL\vdash [n]$ with subsets of given cardinalities $\boldlambda\eqdef\boldlambda(\mbfL)\vdash n$ is counted by the \emph{Fa\`a di Bruno's coefficient} (or \emph{multinomial number of the third kind})
\begin{align*}
\Bella{\boldlambda}\eqdef n!\tonde{\prod_i^n (i!)^{\lambda_i}\lambda_i!}^{-1}=\frac{n!}{\boldlambda! (\vec\mbfn !^\compo)^\boldlambda} \fstop
\end{align*}

An interpretation of~$\Bella{\emparg}$ in terms of set partitions lattices may be found in~\cite[(2.3.8)]{PecTaq11}.

A permutations $\pi\in \mfS_n$ is said to have \emph{cyclic structure~$\boldlambda$} if the lengths of its cycles coincide with $\boldlambda$ (a permutation is thus always understood in its \emph{one-line notation}). Denote by~\mbox{$\mfS_{n}(\boldlambda)\subset \mfS_n$} the set of permutations with cyclic structure~$\boldlambda$ and recall that two permutations have the same cyclic structure if and only if they are conjugate to each other as group elements of~$\mfS_n$. Also recall (cf.~\cite[I.1.3.2]{Sta01}) that
\begin{align*}
\Mnom{\boldlambda}\eqdef\#\mfS_{n}(\boldlambda)=n!\tonde{\prod_i^n i^{\lambda_i}\lambda_i!}^{-1}=\frac{n!}{\boldlambda!\, \vec{\mbfn}^{\boldlambda}} \comma
\end{align*}
termed \emph{multinomial number of the second kind}. Given $\pi\in \mfS_n$, denote by $\boldlambda(\pi)$ the unique integer partition of $n$ encoding the cycle structure of $\pi$, so that $\pi\in \mfS_{n}(\boldlambda(\pi))$.

Finally, a \emph{multi-index} $\mbfm$ of \emph{size} $k$ and \emph{length} $n$ is any vector in $\N_0^k$ such that $\length{\mbfm}=n$. Such multi-indices encode the functions $f\colon [n]\rar [k]$ by regarding $m_i$ as $\#f^{-1}[i]$. Their count is given by the multinomial coefficient~$\tbinom{n}{\mbfm}$ of the multi-index.

In a similar fashion, a {positive multi-index} of size $k$ and length $n$ is any vector in $\N_1^k$ such that $\length{\mbfm}=n$. Such multi-indices encode \emph{surjective} functions in $[k]^{[n]}$, since (strict) positivity of every entry in $\mbfm$ entails that every element of $[k]$ is targeted by $f$. Up to permutation of the elements of $[k]$, surjective functions in $[k]^{[n]}$ are bijective to set partitions of $[n]$ into $k$ blocks, hence the latter ones' total count is given by $\tfrac{1}{k!}\tbinom{n}{\mbfm}$. 

\begin{rem}
For the sake of completeness, let us point out that --~in this context~-- the multinomial coefficient $\MultiZero{\mbfm}\eqdef \tbinom{\length{\mbfm}}{\mbfm}$ ought to be thought of as the first element in a sequence~$\MultiZero{\emparg}$,~$\Mnom{\emparg}$,~$\Bella{\emparg}$ and deserves the name of \emph{multinomial number of the first kind}. We prefer however to keep the usual notation and terminology, but we shall say `multinomial numbers' when collectively referring to the three of them.
\end{rem}

As well as being encoded by a partition $\boldlambda\vdash_r n$, the cyclic structure of a permutation~$\pi\in \mfS_n$ with~$r$ cycles is also encoded by a multi-index $\mbfm\in \N_0^k$ of length~$\length{\mbfm}=n$. Visually, the set of such multi-indices is bijective to that of monotone excursions on the lattice $\Z^2$ starting at $(0,1)$ and ending at $(n,k)$. In this case, $r$ counts the number of non-zero elements in $\mbfm$, which in turn count the lengths of cycles. The general statement is readily deduced from the following example via a \emph{Catalan-type diagram}
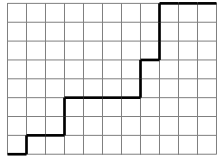
\begin{align*}
\begin{tikzpicture}[scale=0.25]
  \catalannumber{0,-9}{11}{8}{0,1,0,0,1,1,0,0,0,0,1,1,0,1,1,1,0,0,0};
\end{tikzpicture}
\end{align*}

The horizontal lines describe the cycles of the following permutation, resp. partition, multi-index 
\begin{align*}
\pi=&(\emparg)(\emparg\emparg)(\emparg\emparg\emparg\emparg)(\emparg)(\emparg\emparg\emparg)\in \mfS_{11}\comma\\
\boldlambda=&\seq{{}_1 2,{}_2 1, {}_3 1, {}_4 1, {}_5 0, \dotsc, {}_{11} 0}\vdash_5 11\\
\mbfm=&\seq{1,2,0,4,0,1,0,0,3}\in \N_0^9\fstop
\end{align*}

The latter is obtained by counting the length of horizontal lines on the lower side of each square, separated by as many $0$'s as the number of vertical lines on the right side of each square, reduced by one. The number of parts in which $n$ is subdivided, given by $\length{\boldlambda}$, accounts for the number of non-zero entries in $\mbfm$, hence for the number of different groups of contiguous horizontal lines in the diagram.

\begin{rem}\label{r:Duality} When regarded as encoding the cyclic structure of a permutation, multi-indices $\mbfm$ and integer partitions $\boldlambda$ ought to be understood as `dual' notions, the `duality' being given by flipping (along the main diagonal) the Catalan type diagram for~$\mbfm$ constructed above. Loosely speaking, if a property involving a multi-index~$\mbfm$ in~$\N_0^k$ and of length~$n$ holds, then one would expect some `dual' property to hold for the partition $\boldlambda$ corresponding to~$\mbfm$. As the correspondence between integer partitions and multi-indices is not bijective --~hence the quotation marks~--, we shall refrain from addressing it further. It will however be of guidance in discussing instances of the \emph{aggregation property} of the Dirichlet distribution in~\S\ref{s:Dir} below.
\end{rem}

\paragraph{Bell polynomials and the cycle index polynomial of $\mfS_n$} 
The cluster structure of a partition~\mbox{$\mbfL\vdash_r [n]$}, resp.~$\mbfL \vdash [n]$, is encoded by the \emph{partial}, resp. \emph{complete}, \emph{Bell polynomial}
\begin{align}\label{eq:Bell}
B_{n\,r}[\mbfx]\eqdef&\sum_{\boldlambda\vdash_r n} \Bella{\boldlambda} {\mbfx^\boldlambda}\comma &  \textrm{resp. } && B_n[\mbfx]\eqdef& \sum_{\boldlambda\vdash n} \Bella{\boldlambda} {\mbfx^\boldlambda} \comma
\end{align}
where $\boldlambda\eqdef\boldlambda(\mbfL)$ and the index of each variable in the monomial $\mbfx^\boldlambda=x_1^{\lambda_1}\cdots x_n^{\lambda_n}$ indicates the size of the cluster, i.e. there are~$\lambda_1$ clusters of size $1$, up to $\lambda_n$ clusters of size $n$. 
Bell polynomials count the \emph{Bell number} $B_n\eqdef B_n[\uno]$ of partitions of~$[n]$ and satisfy the identities
\begin{align}\label{eq:MultiBell}
B_{n\,r}[b(a\uno)^{\compo \vec \mbfn}\compo \mbfx]=&a^nb^r B_{n\,r}[\mbfx] & B_n[(a\uno)^{\compo\vec\mbfn}\compo\mbfx]=&a^nB_n[\mbfx] \comma
\end{align}
the recursive identity
\begin{equation}\label{eq:RecBell}
\begin{aligned}
B_{n+1}[\mbfx]=&\sum_{k=0}^n \tbinom{n}{k} B_{n-k}[\mbfx_{n-k}] x_{k+1} & \qquad B_0[\emp]\eqdef&1 \comma & \mbfx_\ell\eqdef& (x_1,\dotsc, x_\ell) \comma
\end{aligned}
\end{equation}
and the \emph{binomial type}~(cf.~\cite{RomRot78}) identity
\begin{align}\label{eq:BinomType}
B_n[\mbfx+\mbfy]=&\sum_{k=0}^n \tbinom{n}{k} B_k[\mbfx_k] B_{n-k}[\mbfy_{n-k}] \fstop
\end{align}

In the same way as the Bell polynomial $B_n$ encodes (the cardinalities of) set partitions of~$[n]$, the \emph{cycle index polynomial of~$\mfS_n$} defined by
\begin{align}\label{eq:CI}
Z_n[\mbfx]\eqdef& \frac{B_n[\Gamma^\compo[\vec\mbfn]\compo \mbfx]}{n!}= \frac{1}{n!}\sum_{\boldlambda\vdash n} \Mnom{\boldlambda} \mbfx^\boldlambda
\end{align}
encodes integer partitions of $n$ and satisfies the recurrence relation
\begin{align}\label{eq:RecBellZ}
(n+1)! Z_n[\mbfx]=\sum_{k=0}^n n! Z_{n-k}[\mbfx_{n-k}]x_{k+1} &\qquad Z_0[\emp]\eqdef 1  \comma & \mbfx_\ell\eqdef& (x_1,\dotsc, x_\ell) \fstop 
\end{align}

\begin{rem}[Cycle index polynomials of permutation groups]\label{r:CIP}
Cycle index poly\-no\-mials may be defined for any permutation group~$G\subset \mfS_n$ by the formula
\begin{align*}
Z^G[\mbfx]\eqdef \frac{1}{\# G} \sum_{\pi\in G} \mbfx^{\boldlambda(\pi)}	\fstop
\end{align*}

One of their numerous applications is \emph{P\'olya Enumeration Theory}, whose main result is recalled in~\S\ref{ss:Polya} below.
\end{rem}


A combinatorial proof of \eqref{eq:RecBell} and \eqref{eq:RecBellZ} is best given in terms of \emph{Young tableaux}. Indeed, let e.g. $\boldlambda=\seq{1,1,0,2,0,\dots,{}_{11}0}\vdash_4 11$, which induces by \emph{augmentation} the following partitions of~$12$
\begin{align*}
\boldlambda={\tiny\Yvcentermath1 \yng(4,4,2,1)} \qquad \longrightarrow \qquad {\tiny\Yvcentermath1 \yng(5,4,2,1) \qquad  \yng(4,4,3,1) \qquad \yng(4,4,2,2) \qquad \yng(4,4,2,1,1)} \comma
\end{align*}
obtained by adding to each subset of $\boldlambda$ one additional element in such a way that the associated tableau remains of Young type. Labelling each box in a tableau allows to consider set partitions rather than integer partitions.

For arbitrary $n$, the set partitions of $[n+1]$ are thus obtained by fixing~$S\subset [n]$ a subset of~$k$ elements and adding to it the additional element~$n+1$, so that the size $\#S\cup\set{n+1}=k+1$ of the resulting cluster is encoded in the variable $x_{k+1}$, while the partitions of the remaining~$n-k$ elements are encoded in the Bell polynomial $B_{n-k}$. Since there are~$\binom{n}{k}$ possible choices for the subset~$S$, equation~\eqref{eq:RecBell} follows. Taking this choice to be irrelevant, for each of the so chosen clusters $S$ has cardinality~$k+1$, we get \eqref{eq:RecBellZ}.

\subsection{Exponential generating functions}\label{ss:EGFs}
Given a sequence $\seq{a_n}_n$ of real numbers, we denote by $\EGF[a_n](x)\eqdef \sum_{n\geq0} a_n x^n/n!$ its \emph{exponential generating function}.

\paragraph{EGF's and partitions} It is a well-known fact in enumerative combinatorics that exponential generating functions provide yet another way to account for multinomial numbers, hence to count e.g. set partitions.

\begin{prop}
Let $f[x]\eqdef \EGF[a_n](x)$ and $f_\ell[x]\eqdef \EGF[a_n^{(\ell)}](x)$ for $\ell\in [k]$. It holds
\begin{subequations}
\begin{align}
\label{eq:EGF1} \prod_\ell^k f_\ell[x] =& \EGF[a_n](x) & a_n\eqdef& \sum_{\substack{\mbfm\in \N_0^k\\ \length{\mbfm}=n}} \tbinom{n}{\mbfm} a_{m_1}^{(1)}\cdots a_{m_k}^{(k)} \\
\label{eq:EGF2} \frac{f[x]^k}{k!}=&\EGF[b_n](x) & b_n\eqdef & \sum_{\substack{\mbfm\in \N_1^k\\ \length{\mbfm}=n}} \frac{1}{k!} \tbinom{n}{\mbfm} a_{m_1}\cdots a_{m_k} \\
\label{eq:EGF3} e^{f[x]}=&\EGF[c_n](x) & c_n\eqdef & \sum_{\boldlambda\vdash n} \Bella{\boldlambda} a_1^{\lambda_1} \cdots a_n^{\lambda_n}=B_n[a_1,\dotsc, a_n]  \fstop
\end{align}
\end{subequations}
\begin{proof}
A proof of~\eqref{eq:EGF1} (whence of~\eqref{eq:EGF2}) follows by induction from the classical Cauchy product of power series. Proofs of~\eqref{eq:EGF1} and of the formula for the composition of EGF's (whence of~\eqref{eq:EGF3}) are also found in~\cite[5.1.3, 5.1.4]{Sta01}, together with the respective combinatorial interpretations. We rather discuss the specific combinatorial interpretation of~\eqref{eq:EGF3}, which is as follows.
Firstly, we have from~\eqref{eq:EGF2}
\begin{align*}
e^{f[x]}=\sum_{k=0}^\infty \frac{f[x]^k}{k!}=\sum_{k=0}^\infty \sum_{\substack{\mbfm\in \N_1^k\\ \length{\mbfm}=n}} \frac{1}{k!} \tbinom{n}{\mbfm} a_{m_1}\cdots a_{m_k} \fstop
\end{align*}

As shown in the last paragraph, the coefficient~$\tfrac{1}{k!}\tbinom{n}{\mbfm}$ in~\eqref{eq:EGF2} accounts for the number of surjective functions in~$[k]^{[n]}$, which we can regard as set partitions of~$[n]$. Since~$k$ is itself varying in~$\N$, we are interested in set partitions~$\mbfL\vdash [n]$ into clusters of arbitrary cardinality (hence in set partitions with arbitrary number of clusters). Letting the subscript~$m_i$ of~$a_{m_i}$ denote the cardinality of the set~$L_i$ in $\mbfL$, the number of clusters in~$\mbfL$ with cardinality~$m_i$ is then given by coefficient~$\lambda_i$ in the integer partition $\boldlambda\eqdef\boldlambda(\mbfL)$, whence~\eqref{eq:EGF3} follows by the above discussion on Bell polynomials.
\end{proof}
\end{prop}

\begin{ese} The exponential generating function for the number of cyclic permutations of order $n$, given by $\#\mfS_n(\mbfe^n)=\Mnom{\mbfe^n}=\Gamma[n]$, is
\begin{align}\label{eq:EGFcycle}
\EGF[\Mnom{\mbfe^n}](x)=\sum_{n} \frac{\Gamma[n] \,x^n}{n!}= -\log(1-x) \fstop
\end{align}

As a simple check, since every permutation $\pi\in \mfS_n$ is uniquely \emph{partitioned} into cyclic permutations of maximal length, the exponential generating function for the total number of permutations is obtained by exponentiating the one above, viz. 
\begin{align*}
\EGF[\#\mfS_n](x)=\frac{1}{1-x}=\sum_{n=0}^\infty n! \frac{x^n}{n!} \fstop
\end{align*}
\end{ese}

\subsection{Multi-sets}\label{ss:Multi-sets}
We conclude this section by recalling some basic properties of multi-sets.

\paragraph{Notation} Given a set $S$, a \emph{finite $S$-multi-set} is any function $s\colon S\rar \N$ such that $\#s\eqdef\scalar{\#}{s}$ is finite, where~$\scalar{\#}{\emparg}$ denotes integration on $S$ with respect to the counting measure. We are mainly interested in \emph{real} multi-sets, which we denoted by
\begin{align*}
\{\underbrace{s_1,\dotsc, s_1}_{\alpha_1},\underbrace{s_2,\dotsc, s_2}_{\alpha_2},\dotsc\}\defeq \bag{\mbfs_\boldalpha} \comma
\end{align*}
where~$\boldalpha$ is the vector of positive values attained by the function~$s$ and~$\mbfs$ is a real valued vector with mutually different entries; since $\#\bag{\mbfs_\boldalpha}\eqdef \#s$ is finite,~$s$ is a simple function and thus $\boldalpha$ has finitely many entries; furthermore~$\#\bag{\mbfs_\boldalpha}=\length{\boldalpha}$.

We term the set~$\set{s_1,\dotsc, s_k}$ the \emph{underlying set} to~$\bag{\mbfs_\boldalpha}$, the number~$\alpha_i$ the multiplicity of~$s_i$,~$k$ the number of (different) \emph{types} in~$\bag{\mbfs_\boldalpha}$ and~$\length{\boldalpha}$ the cardinality of the multi-set.
Finally, recall that the total number of multi-sets of cardinality~$n$ with~$k$ types is given by the number of solutions in~$\N^k$ of the equation $\length{\boldalpha}=n$ and counted (see~\cite[\S{I.1.2}]{Sta01}) by the \emph{multi-set} coefficient 
\begin{align*}
\mtbinom{k}{n}\eqdef \Poch{k}{n}/n! \fstop
\end{align*}

For further purposes notice that the expression~$\mtbinom{\alpha}{n}$ is meaningful for any real~\mbox{$\alpha>0$}.


\section{Moments and cumulants on nuclear spaces}\label{a:MomCum}
We give here a short review of moments and cumulants in fairly broad generality, namely for (probability) measures on \emph{nuclear dual spaces} (i.e. nuclear spaces admitting a strongly continuous pre-dual, not to be confused with {dual nuclear} or {co-nuclear spaces}, i.e. spaces with nuclear dual) and more generally for cylindrical measures on topological linear spaces.

\subsection{Infinitely divisible distributions and the L\'{e}vy--Khintchine formula}\label{ss:LK}
\emph{For the sake of simplicity, we review here L\'{e}vy processes on~$\R^d$, mainly following~\cite{Sat13, App04, Str11}. Generalizations to measure-valued L\'evy processes on Riemannian manifolds, L\'evy processes on (separable) Banach and nuclear spaces may be found in~\cite{Lyt03},~\cite{App07} and in~\cite{Ust84} respectively. For convolution and infinite divisibility of measures on finite-dimensional spaces see~\cite[\S7]{Sat13} or~\cite[\S3.2.1]{Str11}; for the convolution calculus on spaces of distributions with applications to stochastic differential equations see e.g.~\cite{BCh02}. For those L\'{e}vy processes on Polish spaces relevant in the applications to random measures below specific references are provided in~\S\ref{s:DirFer}.}

\paragraph{Infinite divisibility} A~probability distribution $\nu$ on $\R^d$ is termed \emph{infinitely divisible} if for every positive integer~$k$ there exists a probability distribution~$\nu_k$ such that~$\nu=\nu_k^{\ast k}$, where~$\lambda^{\ast k}$ denotes the $k$-fold convolution of the measure~$\lambda$ with itself and we set $\lambda^{\ast 0}\eqdef \delta_{\zero_d}$. In words, a distribution~$\nu$ is infinitely divisible if it has convolution roots of arbitrary integer order in the convolution algebra of measures on~$\R^d$, in which case each convolution root is unique (see~\cite[7.5-6]{Sat13}) and it is possible to define arbitrary (non-integer) convolution powers (see~\cite[\S7]{Sat13}) of~$\nu$. Term further a (probability) distribution~$\nu$ on~$\R^d$ to be of Poisson type if there exists~$c>0$ and a probability distribution~$\lambda$, defined on the same space as $\nu$, such that
\begin{align}\label{eq:ConvExp}
\nu= \nu_{c\lambda } =\exp^*[c\lambda] \eqdef e^{-c}\sum_{n=0}^\infty \frac{c^k}{k!} \lambda^{\ast k} \comma
\end{align}
where~$e^{-c}$ is a normalization constant. We term~$\lambda$ the \emph{convolution logarithm} of~$\nu_{c\lambda}$ and denote it by~$\ln^*[\nu_{c\lambda}]$. Morally, a distribution is infinitely divisible if and only if it admits a well-defined convolution logarithm; it is indeed possible to show (see~\cite[\S3.2.1]{Str11}) that every infinitely divisible distribution is in the closure of Poisson type measures with respect to the narrow topology.

\paragraph{L\'{e}vy processes} Recall that a \emph{L\'{e}vy process} \emph{in law} on~$\R^d$ is any stochastically continuous process starting at~$\zero_d$ with stationary independent increments; any such process admits a modification with c\`adl\`ag paths (see~\cite[2.1.7]{App04}), usually termed a \emph{L\'{e}vy process}; since we are mainly interested in the laws of L\'{e}vy processes, the term `in law' is henceforth omitted and, whenever relevant, only modifications with c\`adl\`ag paths are considered.
A~distribution on $\R^d$ is the law of some L\'{e}vy process if and only if it is infinitely divisible (see~\cite[7.10]{Sat13}) and distributions as such are characterized by the following well-known representation theorem.

\begin{thm}[L\'{e}vy--Khintchine formula (see e.g.~{\cite[8.1]{Sat13}})]
A probability distribution~$\nu$ on~$\R^d$ is infinitely divisible if and only if there exist
\begin{itemize}
\item a constant~$\mbfc$ in~$\R^d$;
\item a non-negative definite quadratic form~$A$ on~$\R^d$;
\item a \emph{L\'{e}vy measure}~$\lambda$ on $\R^d$, i.e. satisfying 
\begin{align}\label{eq:LevyMeas}
\lambda\set{\zero_d}=&0 & \textrm{ and } && \int_{\R^d}\tonde{1\wedge \norm{\mbfx}_2^2}\diff\lambda(\mbfx)<&\infty \comma
\end{align}
\end{itemize}
such that
\begin{align}\label{eq:LK}
\mcF[\nu](\mbfz)=\exp\quadre{-\tfrac{1}{2}A[\mbfz,\mbfz]+i \, \mbfc\cdot \mbfz+ \int_{\R^d} \tonde{e^{i\, \mbfz\cdot \mbfx} -1 - i (\mbfz\cdot \mbfx) \car_{D}(\mbfx)} \diff\lambda(\mbfx) } \comma
\end{align}
where $\mbfz\in\R^d$ and $D$ denotes the closed ball in $\R^d$ of radius $1$, centered at the origin.
Furthermore, if this is the case, then $\mbfc$, $A$ and $\lambda$ are unique.
\end{thm}

It is well-known (see e.g.~\cite[2.5\&7.5]{Sat13}) that the Fourier transform of an infinitely divisible distribution on~$\R^d$ is a \emph{non-negative definite} (in the sense of~\eqref{eq:NonNegDef} below) functional on $\R^d$, continuous at~$\zero_d$ with value~$1$, and nowhere vanishing. Whereas these conditions are not sufficient to ensure infinite divisibility, we are interested in those infinite-dimensional linear spaces $\Phi$ (playing the r\^ole of~$\R^d$) such that any functional as above is the Fourier transform of some probability measure on~$\Phi$. A wide class of such spaces, seemingly sufficiently general to include most applications, is constituted by \emph{nuclear spaces}.

\subsection{Measures on nuclear spaces}\label{ss:Radon}
\emph{We follow~\cite[\S{IV.1.1-4}, IV.4.1-2]{GelVil64}. A modern treatment of analysis on nuclear spaces may be found in~\cite{BerKon95}.}
\emph{Topological linear spaces are always meant to be Hausdorff; the term is however omitted.}
\paragraph{Cylindrical measures}
Let $\Phi$ be a co-nuclear space and denote by $\Phi'$ its topological dual (a nuclear dual space) endowed with the weak*~topology induced by the canonical duality pair~$\scalar{\emparg}{\emparg}$.  Let further~$F$ be a finite-dimensional linear subspace of~$\Phi$ and denote by~$F^\perp$ its annihilator in~$\Phi'$. Given a subset~$A$ of the finite-dimensional vector space~$\quotient{\Phi'}{F^\perp}$ define the \emph{cylinder set} with base set $A$ and generating subspace~$\quotient{\Phi'}{F^\perp}$ as the pre-image of $A$ under the quotient map~$\pr\colon \Phi'\rar \quotient{\Phi'}{F^\perp}$. Since $\quotient{\Phi'}{F^\perp}$ is finite-dimensional, it may be endowed with a unique (locally convex) vector space topology and the induced Borel $\sigma$-algebra. By a \emph{cylindrical measure}~$\nu$ on~$\Phi'$ we mean a non-negative real-valued function on the cylinder sets of~$\Phi'$, countably additive on families of cylinder sets with same generating subspace and disjoint Borel measurable base sets, additionally satisfying the normalization~$\nu \Phi'=1$. A~cylindrical measure $\nu$ is termed \emph{continuous} if the functionals
\begin{align*}
\av{u\circ \mbfy}_\nu\eqdef \int_{\Phi'} u\quadre{\mbfy \compo (h\uno)} \diff \nu(h)
\end{align*}
 are sequentially continuous in the variables $\mbfy$ on $\Phi^{\times k}$ for every real-valued function~$u$ continuous and bounded on~$\R^k$ and arbitrary $k$ in $\N_1$. With common abuse of notation, we hereby implicitly denote by $y$ the linear functional $\scalar{\emparg}{y}$ for~$y$ in~$\Phi$, hence we set
\begin{align*}
\mbfy\compo (h\uno)\eqdef \seq{y_1 h, \dotsc, y_k h}=\seq{\scalar{h}{y_1}, \dotsc, \scalar{h}{y_k}}\comma
\end{align*}
and by~$\av{\emparg}_\nu$ the expectation with respect to~$\nu$.
Any continuous cylindrical measure on a nuclear dual space is countably additive (see~\cite[\S{IV.2.3}]{GelVil64}), thus in fact a probability measure.

\paragraph{Fourier \& Laplace transforms} Given a cyclindrical measure $\nu$ on $\Phi'$, define its \emph{Fourier}, resp. \emph{Laplace}, \emph{transform} in the variable~$y$ as
\begin{align*}
\mcF[\nu](y)\eqdef& \av{\exp i y}_\nu \comma &\textrm{resp. } && \mcL[\nu](y)\eqdef& \av{\exp y }_\nu \fstop%
\end{align*}

We borrow the term \emph{transform} from infinite-dimensional analysis; in the case~\mbox{$\Phi=\R^d$} those of \emph{characteristic}, resp. \emph{moment generating}, \emph{function} are of common use.
It is not difficult to show (see~\cite[\S{IV.4.1}]{GelVil64}) that the Fourier (or Laplace) transform of $\nu$ is actually characterized by the $\nu$-measure of half-spaces in $\Phi'$. Indeed, for any finite-dimensional subspace $F\subset \Phi$, every half-space of the form $y\leq c$ with $y$ in $F$ and $c$ a real constant consists of cosets induced by the quotient map $\pr\colon \Phi'\rar \quotient{\Phi'}{F^\perp}$ (see above). Thus, the $\nu$-measure of the half-space~$y\leq c$ in~$\Phi'$ coincides with the $\pr_\pfwd \nu$-measure of the half-space defined by the same inequality in the quotient space~$\quotient{\Phi'}{F^\perp}$ and it holds that
\begin{align*}
\mcF[\nu](y)=\mcF[\pr_\pfwd\nu](y) \fstop
\end{align*}

\paragraph{The Bochner--Minlos Theorem}
Among the reasons why we are interested in cylindrical measures on nuclear spaces is the following well-known realization theorem, which allows for the identification of those linear functionals on a space~$\Phi$ that are Fourier transforms of some (cylindrical) measure on the dual space~$\Phi'$.
Recall that a functional~$\mcF$ on a topological linear space~$\Phi$ is termed \emph{non-negative definite} if
\begin{align}\label{eq:NonNegDef}
\sum_{i,j}^k \mcF(y_i-y_j) \alpha_i \overline{\alpha}_j\geq 0
\end{align}
for every $\mbfy=\seq{y_1,\dotsc, y_k}$ in~$\Phi^{\times k}$,~$\boldalpha$ in~$\C^k$ and~$k$ in~$\N_1$.

\begin{thm}[Bochner--Minlos (see~{\cite[\S{IV.4.2}]{GelVil64}})]\label{t:BM}
Let $\mcF$ be a functional on a topological linear space $\Phi$. Then, $\mcF$ is the Fourier transform of a cylindrical measure~$\nu$ on~$\Phi'$ if and only if it is non-negative definite, sequentially continuous and such that~$\mcF[0_\Phi]=1$.

If, in addition, $\Phi$ is a nuclear space (hence $\Phi'$ is a nuclear dual space), the same statement holds for (countably additive, rather than only cylindrical) probability measures.
\end{thm}

\subsection{Moments, central moments and cumulants}\label{ss:MomCum}
\emph{We are solely concerned with \emph{linear} moments and cumulants, i.e. those of \emph{linear} functionals. However, when referring to moments and cumulants, the term `linear' is always omitted.}

\paragraph{Notation} Let $\nu$ be a cylindrical measure on a (topological) dual linear space $\Phi'$. By (\emph{raw}), resp. \emph{central}, (\emph{univariate}) \emph{moments} of $\nu$ we mean any integral of the form
\begin{align*}
\mu^{\prime \, \nu}_n[y]\eqdef& \av{y^n}_\nu\comma & & \textrm{resp. } & \mu_n^\nu[y]\eqdef &\av{(y-\av{y}_\nu)^n}_\nu \comma
\intertext{where $y^n\eqdef\scalar{\emparg}{y}^n$ for fixed $y$ in $\Phi$. We drop the superscript $\nu$ whenever the measure may be inferred from the context, thus writing only $\mu_n$ in place of $\mu_n^\nu$.
Let now $\mbfy\eqdef \seq{y_1,\dotsc,y_k}$ be a vector in $\Phi^{\times k}$ and denote the vector $\seq{\scalar{\emparg}{y_1}, \dotsc, \scalar{\emparg}{y_k}}$ in the same way. The multivariate Fourier, resp. Laplace, transform in the variables $\mbfy$ is defined as }
\mcF[\nu](\mbfy)\eqdef& \av{\exp[i \length{\mbfy}]}_\nu \comma & &\textrm{resp. } & \mcL[\nu](\mbfy)\eqdef& \av{\exp\length{\mbfy}}_\nu \fstop
\intertext{
\indent The corresponding multivariate moments \emph{of order~$n$} are defined analogously by setting}
\mu^{\prime \, \nu}_\mbfm[\mbfy]\eqdef& \av{\mbfy^\mbfm}_\nu\comma & & \textrm{resp. } & \mu_\mbfm^\nu[\mbfy]\eqdef&\av{(\mbfy-\av{\mbfy}_\nu^\compo)^\mbfm}_\nu \comma
\intertext{where $\mbfm$ is any multi-index in $\N_0^k$ of length $n$ and, consistently with the established notation and the aforementioned abuse thereof, we set}
\av{\mbfw}_\nu\eqdef& \av{w_1}_\nu\cdots \av{y_k}_\nu\comma & & \textrm{and } & \av{\mbfw}_\nu^\compo\eqdef& \seq{\av{w_1}_\nu,\dotsc, \av{w_k}_\nu} \fstop
\end{align*}

Whenever $\mbfm=\uno$ (whence $k=n$), we write $\mu_n[\mbfy]$ in place of $\mu_{\uno}[\mbfy]$ in order to keep track of the dimension of~$\mbfy$.
Finally, we denote collectively by $\mu^\circ$ both raw and central moments.
Given $\mbfy$ as above and $L\subset [n]$, denote by $\mbfy^L$ the product $\prod_{i\in L}\scalar{\emparg}{y_i}$. We define further the \emph{total multivariate cumulant} of order $n$ as~(cf.~\cite[3.2.19]{PecTaq11})
\begin{align}\label{eq:MultiCum1}
\kappa_n[\mbfy]=\kappa_{[n]}^\nu[\mbfy]\eqdef \sum_{r=1}^n (-1)^{r-1}\Gamma[r] \sum_{\mbfL\vdash_r [n]} \prod_{L\in \mbfL} \av{\mbfy^L}_\nu \comma
\end{align}
the corresponding univariate version being simply 
\begin{align*}
\kappa_n[y]\eqdef \kappa_n[y\uno] \fstop
\end{align*}

More general definitions of cumulants are possible, which can be straightforwardly adapted to the present setting. For instance (see e.g.~\cite[\S{A}]{LukMar16}), one can recursively define multivariate cumulants with arbitrary indices in terms of lower order moments by setting
\begin{align}\label{eq:MultiCum2}
\kappa_L[\mbfy]\eqdef& \sum_{\emp\subsetneq I\subset L} \av{\mbfy^{L\setminus I}}_\nu\kappa_I[\mbfy] \comma & \kappa_{\emp}[\mbfy]\eqdef&0 \comma
\end{align}
in which case $\kappa_n[\mbfy]=\kappa_{[n]}[\mbfy]$.

Finally, cumulants enjoy the following properties
\begin{enumerate}[$(i)$]
\item \emph{homogeneity}: $\kappa_L[\mbfc\compo\mbfy]=\mbfc^{\uno}\kappa_L[\mbfy]$ for~$\mbfc$ in $\R^k$;
\item \emph{additivity}: if $\mbfy\independent \mbfw$, then~$\kappa_L[\mbfy+\mbfw]=\kappa_L[\mbfy]+\kappa_L[\mbfw]$;
\item \emph{shift equivariance}: $\kappa_1[y+c]=\kappa_1[y]+c$ for $c\in \R$;
\item \emph{shift invariance}: $\kappa_n[y+c]=\kappa_n[y]$ for~$c\in \R$ and~$n\geq 2$;
\item \emph{independence}: $\kappa_{L}[\mbfy]=0$ as soon as $\mbfy_{L_1}\independent \mbfy_{L_2}$ for some non-trivial~\mbox{$\mbfL\eqdef \set{L_1,L_2}\vdash L$}.
\end{enumerate}

The same holds in fact for an arbitrary measure~$\nu$ and random variables~$\mbfY,\mbfW$ (see e.g.~\cite[\S3.1]{PecTaq11} for the easy proofs).

\paragraph{Dualities}
We stress that, \emph{a priori}, none of the above integrals is well-defined.

If the ($k$-vari\-ate) Laplace transform $\mcL[\nu]$ is well-defined and strictly positive on some domain $D\subset \Phi$ (resp. $D\subset \Phi^{\times k}$), one can define the \emph{cumulant generating function of $\nu$} by setting, on the same domain
\begin{align*}
\mcK[\nu]\eqdef \ln \mcL[\nu] \fstop
\end{align*}

In light of the discussion on half-spaces in~\S\ref{ss:Radon}, we shall focus on (cylindrical) measures $\nu$ on $\Phi'$ such that~$\mcL[\nu](ty)$ and~$\mcK[\nu](ty)$ (resp.~$\mcL[\nu](\mbft\cdot\mbfy)$ and~$\mcK[\nu](\mbft\cdot\mbfy)$) are analytic in the variable(s)~$t$ (resp.~$\mbft$) in a neighborhood of~\mbox{$t=0$} (resp. $\mbft=\zero$ in~$\R^k$) for every fixed $y$ in some affine open half-space~$H\subset \Phi$ containing~$0_\Phi$ (resp. $\mbfy$ in some affine hyper-octant of $\Phi^{\times k}$). We shall term these measures to be \emph{analytically of exponential type}, in which case one can check by straightforward computations that the (multivariate) raw moments, resp. (multivariate) cumulants, are precisely the coefficients in the McLaurin expansion of the respective generating functions; that is e.g.
\begin{subequations}
\begin{align}
\mu^{\prime}_\mbfm[\mbfy]=&\, (-i)^{\length{\mbfm}}\,\partial^{\mbfm}_\mbft\restr_{\mbft=\zero} \mcF[\nu](\mbft\cdot \mbfy) & &\mbft\in\R^k \comma\\
\label{eq:CumFourier}\kappa_\mbfm[\mbfy]=&\, (-i)^{\length{\mbfm}}\, \partial^{\mbfm}_\mbft\restr_{\mbft=\zero} \ln \mcF[\nu](\mbft\cdot \mbfy)& &\mbft\in\R^k \comma
\end{align}
\end{subequations}
where we denote by $\partial^{\mbfm}_\mbft$ the multivariate differential operator of indices $\mbfm$ in the variables $\mbft$. The identity~\eqref{eq:CumFourier} (usually taken as a definition) yields the \emph{cumulance} property, motivating the terminology. Namely, whenever~$\nu_1,\nu_2$ are finite Borel measures on a topological group with analytic Fourier transform  and such that their convolution is defined, then by properties of the latter
\begin{align*}
\mcK[\nu_1\ast \nu_2]=\mcK[\nu_1]+\mcK[\nu_2]
\end{align*}
so that~$\kappa_\mbfm^{\nu_1\ast \nu_2}=\kappa_\mbfm^{\nu_1}+\kappa_\mbfm^{\nu_2}$.

\begin{thm}[Moments/cumulants and central/raw moments dualities]\label{t:MomCum}
Let~$\nu$ be a (cylindrical) measure analytically of exponential type on a (topological) dual linear space~$\Phi'$.
Then, the moments $\mu_n'\eqdef\mu'_n[y]$ and the cumulants~$\kappa_n\eqdef \kappa_n[y]$ may be inferred from each other via the formulae
\begin{subequations}
\begin{align}
\label{eq:CumToMom} \mu'_n=&B_n[\kappa_1,\dotsc,\kappa_n] \comma\\
\label{eq:MomToCum} \kappa_n=&B_n^*[\mu'_1,\dotsc,\mu'_n]\eqdef\sum_{k=1}^n(-1)^{n-k}(n-k)!\,B_{n\,k}[\mu'_1,\dotsc,\mu'_k] \comma
\end{align}
\end{subequations}
where $B_n$ and $B_{n\,k}$ are the Bell polynomials \eqref{eq:Bell}.
Analogously, the raw moments $\mu'_n$ and the central moments $\mu_n\eqdef \mu_n[y]$ are related by the \emph{binomial transforms} (also \emph{M\"obius inversion formulae}~\cite{PecTaq11})
\begin{subequations}
\begin{align}
\label{eq:MomInvA} \mu_n=&\sum_{k=0}^n \tbinom{n}{k} (-1)^{n-k} \mu_1^{\prime \, n-k} \, \mu_k' & \mu'_0=&1 \comma \\
\label{eq:MomInvB} \mu'_n=&\sum_{k=0}^n\tbinom{n}{k} \mu^{\prime \, n-k}_1 \, \mu_k & \mu_0=&1, \mu_1=0 \fstop
%
\intertext{
\indent The corresponding (recursive) multivariate formulae for moments and cumulants are straightforwardly deduced from~\eqref{eq:MultiCum1}-\eqref{eq:MultiCum2}, while those for raw and central moments are as follows} 
\label{eq:MomInvMultiA} \mu_\mbfm[\mbfy]=&\sum_{\zero\leq \mbfj\leq \mbfm}\tbinom{\mbfm}{\mbfj} (-1)^{\length{\mbfm-\mbfj}} \av{\mbfy}_\nu^{\mbfm-\mbfj}\, \mu'_{\mbfj}[\mbfy] & \mu'_{\zero}=&1 \comma \\
\label{eq:MomInvMultiB} \mu_\mbfm'[\mbfy]=&\sum_{\zero\leq \mbfj\leq \mbfm}\tbinom{\mbfm}{\mbfj} \av{\mbfy}_\nu^{\mbfm-\mbfj}\, \mu_{\mbfj}[\mbfy] & \mu_{\zero}=&1, \mu_{\mbfe^i}=0 \fstop
\end{align}
\end{subequations}

\begin{proof} The assumptions on~$\nu$ grant the existence of moments and cumulants in~$y$ (resp.~$\mbfy$) for~$y$ (resp.~$y_i$) in an appropriate affine half-space of~$\Phi$. A proof is then straightforwardly deduced from the classical case~\mbox{$\Phi=\R^d$}. For example, \eqref{eq:CumToMom} is a consequence of the Fa\`a di Bruno-type formula
\begin{align}\label{eq:MomBell}
\diff_t \exp[f(t)]= B_n\quadre{\seq{f^{(1)}(t),\dotsc, f^{(n)}(t)}} \exp[f(t)]
\end{align}
restricted to $t=0$, obtained by differentiating~\eqref{eq:EGF3}. See e.g.~\cite[\S3.2]{PecTaq11}  and~\cite[\S2]{BalJohKot98} for the details in the classical case and~\cite{Smi95} for additional recursive formulae.
\end{proof}
\end{thm}

\section{The Dirichlet distribution}\label{s:Dir}
We take here a digression into some combinatorial aspects of the Dirichlet distribution on the \emph{finite-dimensional} standard simplex, instrumental to the discussion on Dirichlet--Ferguson and Gamma random measures in~\S\ref{s:DirFer} below. Our purpose is twofold: On one hand we show that the combinatorics of symmetric groups 
provides interpretation and, occasionally, new proofs for well-known properties of the Dirichlet distribution; on the other hand, we wish to draw attention on how these combinatorial aspects may be subsumed in the study of Lauricella hypergeometric functions, which we expect to constitute a link towards (Lie) representation theory, loosely in the spirit of~\cite{TsiVerYor00,TsiVerYor01,Ver07}.

\subsection{The Poisson, Gamma and Dirichlet distributions}\label{ss:PoiGamDir}
\emph{We minimally recall the main definitions and some facts on the Poisson, Gamma and Dirichlet distributions. Importantly, we describe the Dirichlet distribution on the standard simplex $\Delta^k\subset \R^{k+1}$ rather than on the corner simplex $\Delta^k_*\subset \R^k$ (in the Literature both choices are common)}.

\paragraph{Poisson distribution} The \emph{Poisson distribution} of \emph{expectation} $c>0$ is the \emph{discrete} distribution with density~$\diff\Poi[c](k)\eqdef e^{-c}c^k/k!$, so that
\begin{align}\label{eq:PoiDelta}
\Poi[c](k)\eqdef \exp^*[c\delta_1] \comma
\end{align}
which motivates the terminology of \emph{(C)PT distribution} in~\S\ref{ss:LK}.
The moments of a Poisson distribution of expectation $c$ are given by $T_n[c]$ where 
\begin{align*}
T_n[x]\eqdef B_n[x\uno]= e^{-x}\sum_{k=0}^\infty \frac{x^kk^n}{k!}
\end{align*}
denotes the $n^\textrm{th}$ \emph{Touchard polynomial}, with exponential generating function
\begin{equation}\label{eq:Touchard}
\EGF[T_n](x,t)\eqdef \sum_{n=0}^\infty \frac{T_n[x]t^n}{n!}=e^{x(e^{t}-1)} \fstop \qedhere
\end{equation}

The moments of the Poisson distribution with expectation $1$ are thus precisely the Bell numbers~$B_n$ (\emph{Dobi\'{n}ski formula}), while the cumulants of~$\Poi[c]$ are simply, cf.~\eqref{eq:CumToMom},
\begin{align*}
\kappa_n[s]=c=\mu_n^{\prime\, \ln^*\Poi[c]}[s] \fstop
\end{align*}

\paragraph{Gamma distribution} The \emph{Gamma distribution} of \emph{shape} $\theta>0$ and (\emph{inverse}) \emph{scale} $k>0$ is defined as
\begin{align*}
\Gam\quadre{\theta,k}(x)\eqdef \frac{k^\theta}{\Gamma[\theta]}x^{\theta-1}e^{-k x}\car_{[0,+\infty]}[x] \fstop
\end{align*}

The distribution is characterized by the following property (\emph{Lukacs'~characterization}, see~\cite{Luk55}): If~$Y$ and~$W$ are independent positive non-degenerate random variables, then $Y+W$ and $Y/(Y+W)$ are independent if and only if $Y$ and $W$ have Gamma distributions with the same scale parameter.
Furthermore, the distribution $\Gam[\theta,k]$ is moment determined, with moment generating function~$\mcL[\Gam[\theta,k]](t)=(1+k t)^{-\theta}$ and moments
\begin{align}\label{eq:MomGam}
\mu^{\prime \,\Gam[\theta,k]}_n=\Poch{\theta}{n} k^n \fstop
\end{align}

Given $\mbfY\eqdef \seq{Y_1,\dotsc, Y_k}$ i.i.d. random variables with $Y_i\sim \Gam[\alpha_i, \ell]$, their sum satisfies 
\begin{align*}
\length{\mbfY} \sim \Gam[\length{\boldalpha},\ell]\comma
\end{align*}
hence, inductively, $V$ is independent of the distribution of the $\Delta^{k-1}$-valued random vector $\mbfW\eqdef \length{\mbfY}^{-1}\mbfY$. The distribution 
of the latter is $\mbfW\sim \Dir[\boldalpha]$, where $\Dir$ is the following \emph{Dirichlet distribution} on the simplex $\Delta^{k-1}$.

\paragraph{Dirichlet distribution}
Let $\Delta^{k-1}\eqdef\conv\seq{\mbfe^i}_i^{k}\subset \R^{k}$ be the standard $(k-1)$-dimensional simplex. We denote by
\begin{align*}
\Dir\quadre{\boldalpha}(\mbfx)\eqdef \frac{\car_{\Delta^{k-1}}\quadre{\mbfx} }{\Beta\quadre{\boldalpha}} \prod_{i=1}^{k} x_i^{\alpha_i-1}= \car_{\Delta^{k-1}}\quadre{\mbfx} \frac{\mbfx^{\boldalpha-\uno}}{\Beta\quadre{\boldalpha}}
\end{align*}
the \emph{Dirichlet distribution with parameter~$\boldalpha>0$} (see e.g.~\cite{NgTiaTan11}).
The Dirichlet distribution is \emph{quasi-exchange\-able}, i.e.
\begin{align}\label{eq:Symmetry}
\Dir[\boldalpha](\mbfx)=&\Dir[\boldalpha_\pi](\mbfx_\pi)& \pi\in& \mfS_{k}
\end{align}
and satisfies the following \emph{aggregation property} (see e.g.~\cite[\S49.5]{JohKotBal00}) . Setting
\begin{align*}
\mbfy_{\contra+ i}\eqdef\mbfy_{\hat{\imath}}+y_{i}\mbfe^{i}=\seq{y_1, \dotsc, y_{i-1}, y_i+y_{i+1}, y_{i+2}, \dotsc, y_k} \comma
\end{align*}
a random vector $\mbfY$ satisfies
\begin{align}\label{eq:Aggregation}
\mbfY\sim \Dir[\boldalpha]\implies& \mbfY_{\contra+ i}\sim \Dir[\boldalpha_{\contra+ i}] \fstop
\end{align}

In the case $\boldalpha=\alpha\uno$ one has $\boldalpha_\pi=\boldalpha$ for every $\pi$, thus quasi-exchangeability coincides with exchangeability, hence the name.

\smallskip

By compactness of $\Delta^{k-1}$, the Dirichlet distribution is moment determined for any choice of the parameters~$\boldalpha$. As a consequence, any property holding for the Dirichlet distribution may be checked for a measure~$\nu$ with the same moments as~$\Dir[\boldalpha]$.
The moments are straightforwardly computed as
\begin{equation}\label{eq:MomDir}
\begin{aligned}
\mu'_n[\mbfs,\boldalpha]\eqdef& \int_{\Delta^{k-1}} (\mbfs\cdot \mbfx)^n \Dir[\boldalpha](\mbfx) \diff\mbfx \\
=&\sum_{\substack{\mbfm\in\N^k_0\\ \length{\mbfm}=n}} \tbinom{n}{\mbfm} \mbfs^\mbfm \Beta[\boldalpha+\mbfm] \Beta[\boldalpha]^{-1}=\frac{n!}{\Poch{\length{\boldalpha}}{n}} \sum_{\substack{\mbfm\in\N^k_0\\ \length{\mbfm}=n}} \frac{\mbfs^{\mbfm}}{\mbfm!} \Poch{\boldalpha}{\mbfm} \comma
\end{aligned}
\end{equation}
so that the \emph{Laplace transform} of the distribution satisfies
\begin{equation}\label{eq:LapDir}
\begin{aligned}
\mcL[\Dir[\boldalpha]](\mbfs)\eqdef&\int_{\Delta^{k-1}} \exp(\mbfs\cdot\mbfx)\Dir[\boldalpha](\mbfx) \diff\mbfx\\
=& \sum_{n=0}^\infty \frac{\mu'_n[\mbfs,\boldalpha]}{n!}= \sum_{\mbfm\in\N^k_0}\mtbinom{\boldalpha }{ \mbfm }\frac{ \mbfs^\mbfm}{\Poch{\length{\boldalpha}}{\length{\mbfm}}} \defeq {}_k\Phi_2[\boldalpha;\length{\boldalpha};\mbfs] \fstop
\end{aligned}
\end{equation}

We denoted here by ${}_k\Phi_2$ a confluent form of the $k$-variate Lauricella hypergeometric function of type $D$, a brief account of which is given in~\S\ref{ss:Lauricella} below. Let us point out that, since we denoted by~$\mbfs$ a vector in~$\R^k$ meant as a linear functional on the same space, our usual notation for univariate moments established in~\S\ref{ss:MomCum} would be conflicting with the one for multivariate moments and therefore be inconsistent. We thus write~$\mu_n'[\mbfs,\boldalpha]$ in place of~$\mu_n^{\prime\, \Dir[\boldalpha]}[\mbfs\cdot (\emparg)]$.
%

Since the vertices of $\Delta^{k-1}$ are in bijection with $[k]$ and the Dirichlet distribution is quasi-exchangeable, the question arises from Remark~\ref{r:Duality} wether the aggregation property has a natural generalization in terms of partitions of $k$. One notices that this is indeed the case by defining (for $\boldlambda\vdash k$) the \emph{$+$-contraction}
\begin{equation}\label{eq:AddContract}
\begin{aligned}
\mbfy_{\contra+\boldlambda} \eqdef& (\underbrace{y_1,\dotsc, y_{\lambda_1}}_{\lambda_1}, \underbrace{y_{\lambda_1+1}+y_{\lambda_1+2}, \dotsc, y_{\lambda_1+2\lambda_2-1}+y_{\lambda_1+2\lambda_2}}_{2\lambda_2},\dotsc, \\
&\quad \underbrace{y_{\length{(\vec\mbfk\compo\boldlambda)_{k-1}}+1}+\dots+y_{\length{(\vec\mbfk\compo\boldlambda)_{k-1}}+\lambda_k}, \dotsc, y_{\length{\vec\mbfk \compo\boldlambda}-\lambda_k+1}+\dots+y_{\length{\vec\mbfk\compo\boldlambda}}}_{k\lambda_k}) \comma
\end{aligned}
\end{equation}
whence inductively applying \eqref{eq:Aggregation} yields $\mbfY\sim \Dir[\boldalpha] \implies \mbfY_{\contra+\boldlambda}\sim \Dir[\boldalpha_{\contra+\boldlambda}]$ for $\boldlambda \vdash k$.
Combining the latter with the quasi-exchangeability~\eqref{eq:Symmetry}, the Dirichlet distribution satisfies then
\begin{align}\label{eq:AggrGen}
\mbfY\sim \Dir[\boldalpha] \implies& \mbfY_{\pi,\contra+\boldlambda}\sim \Dir[\boldalpha_{\pi,\contra+\boldlambda}] & \pi\in \mfS_k, \boldlambda \vdash k \fstop
\end{align}

The implication~\eqref{eq:Aggregation} is readily recovered by choosing $\boldlambda\eqdef (k-2,1,0,\dotsc)\vdash k$, thus $\mbfy_{\contra+\boldlambda}=\mbfy_{\contra+(k-1)}$, hence by quasi-exchangeability.

\begin{rem}[Models for $\Delta^{k-1}$ and $\mfS_{k}$]\label{r:Models}
Regarding~$\Delta^{k-1}$ endowed with the usual topology as (homeomorphic to) the space of probability measures over $[k]$ endowed with the narrow topology, the aggregation property may be given a natural measure-theoretical interpretation. Indeed, let~$\mfs^{i}\colon [k]\rar [k-1]$ denote the~$i^{\textrm{th}}$ degeneracy map of~$[k]$, i.e. the only weakly order preserving surjection `hitting'~$i$ twice, set $\boldeta\eqdef\sum_{\ell\in [k]} \eta_\ell \delta_\ell$ and interpret the Dirichlet distribution~$\Dir[\boldalpha]$ as induced by a (positive) reference measure~$\boldalpha$ on~$[k]$. Then~$\mfs^i_\pfwd \boldeta=\boldeta_{\contra+i}$ and
\begin{align*}
(\mfs^i_\pfwd)_\pfwd\mbfY_\pfwd \P=(\mfs^i_\pfwd \mbfY)_\pfwd \P=(\mbfY_{\contra+i})_\pfwd \P \comma
\end{align*}
hence, choosing $\mbfY\sim \Dir[\boldalpha]$, the aggregation property~\eqref{eq:Aggregation} translates into
\begin{align*}
(\mfs^i_\pfwd)_\pfwd \Dir[\boldalpha]=\Dir[\boldalpha_{\contra+i}]=\Dir[\mfs^i_\pfwd \boldalpha] \fstop
\end{align*}

In the same way, letting $g\colon [k]\rar [k]$, the aggregation property~\eqref{eq:AggrGen}, together with the quasi-exchangeability, translates into
\begin{align*}
(g_\pfwd)_\pfwd \Dir[\boldalpha]=\Dir[g_\pfwd \boldalpha] \fstop
\end{align*}

Next, we turn to the permutation group~$\mfS_k$ appearing in the definition of quasi-exchange\-abil\-ity above. The group --~or, rather, its action~-- may be given several different interpretations. Certainly it can be regarded as acting on the simplex $\Delta^{k-1}$ by permutation of its vertices, i.e. by permutation of co\"ordinates in the ambient space $\R^k$, or, equivalently, by permutation of the vectors $\set{\mbfe^1,\dotsc, \mbfe^k}$, supporting the Markov--Krein transform~\cite[4.8]{KerTsi01} of~$\Dir[\boldalpha]$. Another more involved and essentially algebraic interpretation is given in Proposition~\ref{p:chased} below.
For the purposes of the next sections, the right interpretation to keep in mind is however the following. Regarding again~$\Delta^{k-1}$ as the space of probability measures over $[k]$,~$\mfS_k$ naturally acts, by push-forward, as the space of measurable isomorphisms (or, equivalently, of homeomorphisms, diffeomorphisms) of~$[k]$. That is, any~$\pi$ in $\mfS_k$ acts on a measure~$\boldeta$ in~$\Delta^{k-1}$ by setting~$\pi.\boldeta\eqdef \pi_\pfwd \boldeta$.
\end{rem}

\subsection{Moments of the Dirichlet distribution}\label{ss:MomDir}
Let us now observe that the (raw) moments $\mu'_n[\mbfs,\boldalpha]$ (see~\eqref{eq:MomDir}) of the $k$-variate Dirichlet distribution are given by a summation indexed over the set of multi-indices in $\N_0^k$ with fixed length $n$. In view of~\S\ref{s:Bell}, it is possible to rearrange the sum indexing it over partitions, and to formulate the result in terms of the cycle index polynomials~\eqref{eq:CI}.

\begin{thm}[Moments of the Dirichlet distribution]\label{l:lemma} 
The following identity holds
\begin{align}\label{eq:l:MomDir}
\mu_n'[\mbfs,\boldalpha]=\frac{n!}{\Poch{\length{\boldalpha}}{n}} \sum_{\substack{\mbfm\in\N^k_0\\ \length{\mbfm}=n}} \frac{\mbfs^{\mbfm}}{\mbfm!} \Poch{\boldalpha}{\mbfm}=\frac{n!}{\Poch{\length{\boldalpha}}{n}} Z_n[\mbfs^{\compo 1}\cdot\boldalpha, \dotsc,\mbfs^{\compo n}\cdot \boldalpha] \defeq \zeta_n[\mbfs,\boldalpha] \fstop
\end{align}

\begin{proof}
We give here an elementary combinatorial proof, independent of any property of the Dirichlet distribution.
Let
\begin{align*}
\tilde\mu_n[\mbfs,\boldalpha]\eqdef& \frac{\Poch{\length{\boldalpha}}{n}}{n!} \mu_n'[\mbfs,\boldalpha]\comma & \tilde\zeta_n[\mbfs,\boldalpha]\eqdef& \frac{\Poch{\length{\boldalpha}}{n}}{n!} \zeta_n[\mbfs,\boldalpha] \fstop
\end{align*}

The statement is equivalent to $\tilde\mu_n=\tilde\zeta_n$, which we prove in two steps.

\medskip

\emph{Step 1}. The following identity holds
\begin{align}\label{eq:Dirichlet0.0}
\tilde\mu_{n-1}[\mbfs,\boldalpha+\mbfe^\ell]=\sum_{h=1}^n s_\ell^{h-1} \tilde\mu_{n-h}[\mbfs,\boldalpha] \fstop
\end{align}

By induction on $n$ with trivial (i.e. $1=1$) base step $n=1$.
\emph{Inductive step.} Assume for every $\boldalpha\in \R^k_+$ and $\mbfs\in\R^k$
\begin{align}\label{eq:Dirichlet0}
\tilde\mu_{n-2}[\mbfs,\boldalpha+\mbfe^\ell]=\sum_{h=1}^{n-1} s_\ell^{h-1} \tilde\mu_{n-1-h}[\mbfs,\boldalpha] \fstop
\end{align}

Notice that
\begin{equation}\label{eq:Dirichlet1}
\begin{aligned}
\partial_j \tilde\mu_n[\mbfs,\boldalpha] =&\sum_{\substack{\mbfm\in \N_0^k\\ \length{\mbfm}=n}} \frac{k_j \mbfs^{\mbfm-\mbfe^j}}{\mbfm !} \Poch{\boldalpha}{\mbfm}=
\sum_{\substack{\mbfm\in \N_0^k\\ \length{\mbfm}=n}} \frac{\mbfs^{\mbfm-\mbfe^j}}{(\mbfm-\mbfe^j) !} \alpha_j\Poch{\boldalpha+\mbfe^j}{\mbfm-\mbfe^j}\\
=& \alpha_j \sum_{\substack{\mbfm\in \N_0^k\\ \length{\mbfm}=n-1}} \frac{\mbfs^{\mbfm}}{\mbfm !}\Poch{\boldalpha+\mbfe^j}{\mbfm}= \alpha_j \tilde\mu_{n-1}[\mbfs,\boldalpha+\mbfe^j] \fstop
\end{aligned}
\end{equation}

If $k\geq 2$, we can choose $j\neq \ell$. Applying \eqref{eq:Dirichlet1} to both sides of \eqref{eq:Dirichlet0.0} yields
\begin{align*}
\partial_j \tilde\mu_{n-1}[\mbfs,\boldalpha+\mbfe^\ell]=& \alpha_j \tilde\mu_{n-2}[\mbfs,\boldalpha+\mbfe^j+\mbfe^\ell]\\
\partial_j \sum_{h=1}^n s_\ell^{h-1} \tilde\mu_{n-h}[\mbfs,\boldalpha]=& \sum_{h=1}^n s_\ell^{h-1} \alpha_j \tilde\mu_{n-h-1}[\mbfs,\boldalpha+\mbfe^j]\\
=& \alpha_j\sum_{h=1}^{n-1} s_\ell^{h-1} \tilde\mu_{n-h-1}[\mbfs,\boldalpha+\mbfe^j] \comma
\end{align*}
where the latter equality holds since $\tilde\mu_{-1}=0$. Letting now $\boldalpha'\eqdef \boldalpha+\mbfe^j$ and applying the inductive hypothesis~\eqref{eq:Dirichlet0} with $\boldalpha\leftarrow \boldalpha'$ yields
\begin{align*}
\partial_j \quadre{\tilde\mu_{n-1}[\mbfs,\boldalpha+\mbfe^\ell]- \sum_{h=1}^n s_\ell^{h-1}\tilde\mu_{n-h}[\mbfs,\boldalpha]}=0
\end{align*}
for every $j\neq \ell$. By arbitrariness of $j$, the bracketed quantity is a polynomial in the sole variables $s_\ell$ and $\boldalpha$ of degree at most $n-1$ (obviously, the same holds also in the case $k=1$). As a consequence of the above statement (or trivially if $k=1$), every summand in which $s_\ell$ does not appear freely cancels out by arbitrariness of $\mbfs$, yielding
\begin{align*}
\tilde\mu_{n-1}[\mbfs,\boldalpha+\mbfe^\ell]- \sum_{h=1}^n s_\ell^{h-1}\tilde\mu_{n-h}[\mbfs,\boldalpha]= \frac{s_\ell^{n-1}}{(n-1)!}\Poch{\alpha_\ell+1}{n-1}-\sum_{h=1}^n s_\ell^{h-1} \frac{s_\ell^{n-h}}{(n-h)!} \Poch{\alpha_\ell}{n-h} \fstop
\end{align*}

The latter quantity is proved to vanish as soon as
\begin{align*}
\frac{\Poch{\alpha+1}{n-1}}{(n-1)!}=\sum_{h=1}^n \frac{\Poch{\alpha}{n-h}}{(n-h)!} \comma  \textrm{or equivalently} \qquad \Poch{\alpha+1}{n-1}=\sum_{h=0}^{n-1} \frac{\Poch{\alpha}{h} (n-1)!}{h!}\comma
\end{align*}
in fact a particular case of the well-known \emph{Chu--Vandermonde identity} (e.g.~\cite[5.3.1]{RomRot78})
\begin{align*}
\Poch{\alpha+\beta}{n}=\sum_{k=0}^n \tbinom{n}{k}\Poch{\alpha}{k}\Poch{\beta}{n-k} \fstop
\end{align*}

\emph{Step 2.} It holds that $\tilde\mu_n=\tilde\zeta_n$. By strong induction on~$n$ with trivial (i.e. $1=1$) base step~$n=0$. \emph{Inductive step.}
Assume for every $\boldalpha\in \R^k_+$ and $\mbfs\in\R^k$ that $\tilde\mu_{n-1}[\mbfs,\boldalpha]=\tilde\zeta_{n-1}[\mbfs,\boldalpha]$. Then
\begin{align*}
\partial_j \tilde\mu_n[\mbfs,\boldalpha]=& \sum_{\boldlambda \vdash n} \frac{\Mnom{\boldlambda}}{n!} \sum_{h=1}^n \frac{\partial_j (\mbfs^{\compo h}\cdot \boldalpha)^{\lambda_h}}{(\mbfs^{\compo h}\cdot \boldalpha)^{\lambda_h} } \prod_{i=1}^n (\mbfs^{\compo i}\cdot \boldalpha)^{\lambda_i} \\
=&\sum_{\boldlambda \vdash n} \frac{\Mnom{\boldlambda}}{n!} \sum_{h=1}^n \frac{h\lambda_h s_j^{h-1}\alpha_j}{\mbfs^{\compo h}\cdot \boldalpha} \prod_{i=1}^n (\mbfs^{\compo i}\cdot \boldalpha)^{\lambda_i}\\
=&\alpha_j \sum_{h=1}^n s_j^{h-1} \sum_{\boldlambda \vdash n} \frac{h\lambda_h}{1^{\lambda_1}\lambda_1! \dotsc h^{\lambda_h}\lambda_h! \dotsc n^{\lambda_n}\lambda_n!} \frac{1}{\mbfs^{\compo h}\cdot \boldalpha} \prod_{i=1}^n (\mbfs^{\compo i}\cdot \boldalpha)^{\lambda_i}
\\
=&\alpha_j \sum_{h=1}^n s_j^{h-1} \sum_{\boldlambda \vdash n-h} \frac{\Mnom{\boldlambda}}{(n-h)!} \prod_{i=1}^{n-h} (\mbfs^{\compo i}\cdot\boldalpha)^{\lambda_i} \\
=&\alpha_j \sum_{h=1}^n s_j^{h-1} \tilde\zeta_{n-h}[\mbfs,\boldalpha] \fstop
\end{align*}

The inductive hypothesis, \eqref{eq:Dirichlet0.0} and \eqref{eq:Dirichlet1} yield
\begin{align*}
\partial_j \tilde\zeta_n[\mbfs,\boldalpha]=\alpha_j \sum_{h=1}^n s_j^{h-1} \tilde\mu_{n-h}[\mbfs, \boldalpha]=\partial_j \tilde\mu_n[\mbfs,\boldalpha] \fstop
\end{align*}

By arbitrariness of $j$ this implies that $\tilde\zeta_n[\mbfs,\boldalpha]-\tilde\mu_n[\mbfs,\boldalpha]$ is constant as a function of $\mbfs$ (for fixed~$\boldalpha$), hence vanishing by choosing~$\mbfs=\zero$.
\end{proof}
\end{thm}

\emph{For the sake of completeness, let us point out that, in the case $\length{\boldalpha}=1$, proofs of the above proposition may also be inferred from~\cite[7.4]{KerTsi01} or~\cite[5.2--5.4]{Fen10}. Our proof is however elementary and does not rely on properties of the distribution.}


\subsection{Urns, beads and P\'olya Enumeration Theory}\label{ss:Polya}
An interpretation of the moments formula~\eqref{eq:l:MomDir} may be given in enumerative combinatorics, by means of \emph{P\'olya Enumeration Theorem}~\ref{t:Polya} below.
We refer the reader to~\cite{dBr64} for a didactic exposition on the subject; a minimal background is as follows.

\paragraph{Urns and beads} Consider the vector $\vec{\mbfn}$ as a set $[n]$ of \emph{labeled} beads, acted upon by a permutation group~$G$. An \emph{arrangement}~$f\in[k]^{[n]}$ is any way of placing the beads into $k$ different urns, labeled by weights~$s_j$ for $j\in[k]$. A \emph{pattern} is any equivalence class of arrangements up to the action of~$G$, i.e. with respect to the equivalence relation defined on~$[k]^{[n]}$ by setting $f_1\sim f_2$ if and only if~$f_1(g.x)=f_2(x)$ for~$g\in G$ and~$x\in[n]$. As an example, choosing~\mbox{$G=\mfC_n$} (the cyclic group of order $n$) we are considering the labeling of the beads modulo~$n$, while choosing~\mbox{$G=\mfS_n$} we are disregarding the labeling of the beads completely.
The \emph{pattern inventory} is then defined as the sum
\begin{align*}
\Pi^G_n[\mbfs]\eqdef\sum_{f^\bullet} \prod_{j\in[k]} s_{f(j)}
\end{align*}
where $f$ is any arrangement yielding the pattern $f^\bullet$ and the sum runs over all such patterns.

\begin{thm}[P\'olya Enumeration Theorem~{\cite[\S4]{Pol37}}]\label{t:Polya}
Standing the terminology and notation above, the pattern inventory is given by
\begin{align}\label{eq:Polya}
\Pi^G_n[\mbfs]=Z^G\quadre{p_1[\mbfs], p_2[\mbfs],\dotsc} \comma
\end{align}
where $p_i[\mbfs]\eqdef \length{\mbfs^{\compo i}}=\mbfs^{\compo i}\cdot\uno$ denotes the $i^\textrm{th}$ \emph{power sum symmetric polynomial} and $Z^G$ is the cycle index polynomial of~$G$ defined in Remark~\ref{r:CIP}.
\end{thm}

In the case of our interest $G=\mfS_n$, thus~\eqref{eq:Polya} reads~\mbox{$Z_n[\mbfs^{\compo 1}\cdot \uno, \dotsc, \mbfs^{\compo n}\cdot \uno]$}, that is, it coincides (up to normalization) with the moment~$\mu'_n[\mbfs,\uno]$ of the Dirichlet distribution~$\Dir[\uno]$ on~$\Delta^{k-1}$.
This provides a new proof of the aggregation property~\eqref{eq:Aggregation}, which simply translates into merging the content of two urns into one.

\begin{prop}[Aggregation property by means of P\'{o}lya Enumeration]\label{p:AggrPolya}
The Dirichlet distribution satisfies the aggregation property~\eqref{eq:AggrGen}.
\begin{proof}
Standing the metaphor of urns above, assume that~$s_i=s_{i+1}$ for some $i$ in~$[k]$, so that the corresponding urns are now indistinguishable. This amounts to replacing~$\mbfs$ with~$\mbft\eqdef \mbfs_{\hat \imath}$ in the pattern inventory count, which becomes
\begin{align*}
\tilde\mu'_n[\mbfs,\uno]=Z_n\quadre{\length{\mbfs^{\compo 1}}, \dotsc,\length{\mbfs^{\compo n}}}=Z_n[\mbft^{\compo 1}\cdot \uno_{\contra+i}, \dotsc, \mbft^{\compo n}\cdot \uno_{\contra+i}]=\tilde\mu'_n[\mbft, \uno_{\contra+i}]\comma
\end{align*}
where $\tilde\mu_n'$ denotes the normalized moment introduced in~\eqref{eq:l:MomDir}. That is, $\tilde\mu'_n[\mbfs,\uno]$ coincides with the normalized moment of the Dirichlet distribution with parameter $\uno_{\contra+i}$ on~$\Delta^{k-2}$.
Since~$\Dir$ is moment determinate (hence also determined by its normalized moments), this proves~\eqref{eq:Aggregation} for~$\Dir[\uno]$.

In light of the generalization of the aggregation property from~\eqref{eq:Aggregation} to~\eqref{eq:AggrGen}, the case when \mbox{$\boldalpha\in \N_1^k$} is similar.
Indeed, the normalized moment~$\tilde\mu'_n[\mbfs,\boldalpha]$ of the Dirichlet distribution coincides with the pattern inventory~$\Pi^{\mfS_n}_n$ when we allow for the urns' weights~$\set{s_1,\dotsc, s_k}$ to be actually the multi-set~$\bag{\mbfs_\boldalpha}$, with~$k$ types~$s_i$ of multiplicity~$\alpha_i$ each. In the metaphor, this amounts to display~$n$ (non-labeled) beads into $\length{\boldalpha}$ different urns, weighted according to the multi-set~$\bag{\mbfs_\boldalpha}$, and then merge all urns with the same weight into one.

Finally, recalling that, for integer-valued~$\boldalpha$, one has~$\#\bag{\mbfs_\boldalpha}=\length{\boldalpha}$, the case when~\mbox{$\boldalpha\in \Q^k_+$} amounts to only consider the occurrence rate ratios of elements in the multi-set, normalized with respect to~$\#\bag{\mbfs_\boldalpha}$. The general case of real-valued~$\boldalpha$ may be given the same interpretation by approximation.
\end{proof}
\end{prop}

Together with the metaphor of urns and beads, the above corollary provides a combinatorial interpretation for the right-hand side of~\eqref{eq:l:MomDir}. To see this, we regard the weights~$\mbfs$ as dummy variables, encoding for the type of each element in a multi-set of cardinality~$\length{\boldalpha}$, itself subdivided into~$k$ different multi-sets of cardinality~$\alpha_i$. By definition of multi-set coefficient one can rewrite~\eqref{eq:MomDir} as
\begin{align*}
\mu'_n[\mbfs,\boldalpha]=\frac{n!}{\Poch{\length{\boldalpha}}{n}} \sum_{\substack{\mbfm\in\N^k_0\\ \length{\mbfm}=n}} \mtbinom{\boldalpha}{\mbfm} \mbfs^\mbfm \comma
\end{align*}
which in turn makes sense for arbitrary non-negative~$\boldalpha$. In conformity with the established notation we set ~$\mtbinom{\boldalpha}{\mbfm}\eqdef \prod_i^k \mtbinom{\alpha_i}{m_i}$.

\paragraph{Operations on the set of urns} As suggested by any of the moments' formulae for $\Dir[\boldalpha]$, the case when~$\length{\boldalpha}=1$ is singled out as both computationally easiest (since~$\Poch{1}{n}=n!$ simplifies the expression of each moment, also cf. Prop.~\ref{p:LaurEGF} below) and most `meaningful',~$\boldalpha$ representing in that case a \emph{probability} distribution on the set~$[k]$ in the sense of Remark~\ref{r:Models}. For these reasons, this is usually the sole considered case (cf. e.g.~\cite{KerTsi01, LijPru10}). On the other hand though, the above metaphor of urns and beads suggests that the case when~$\boldalpha$ is any positive (integer) vector is equally interesting, since it allows for some natural operations on the multi-set of urns, corresponding to functionals of the distribution. Such operations are given in Table~\ref{tab:A}.

\begin{table}[!ht]
\caption{Some natural operations on multi-sets of urns with~$\boldalpha\in \R_+^k, \mbfs\in \R^k$, hence~$\Dir[\boldalpha]$ on~$\Delta^{k-1}$}
\begin{tabular}{ c c c }
operations & \multicolumn{1}{c}{on urns} & \multicolumn{1}{c}{on $\Dir[\boldalpha]$} \\
\midrule
type addition & $\bag{(\mbfs\oplus (s_{k+1}))_{\boldalpha\oplus (\alpha_{k+1})}}$ & $\Dir[\boldalpha\oplus (\alpha_{k+1})]$ \\
type deduction & $\bag{(\mbfs_{k-1})_{\boldalpha_{k-1}}}$ & $\Dir[\boldalpha_{k-1}]$ \\
urns merge & $\bag{(\mbfs_{\hat \imath})_{\boldalpha_{\contra+ i}}}$ & $\Dir[\boldalpha_{\contra+i}]$ \\
urn addition & $\bag{\mbfs_{\boldalpha+\mbfe^i}}$ & $\alpha_i \Dir[\boldalpha+\mbfe^i]$ \\
urn deduction & $\bag{\mbfs_{\boldalpha-\mbfe^i}}$ & $-(1-\length{\boldalpha})\Dir[\boldalpha-\mbfe^i]$
\end{tabular}
\label{tab:A}
\end{table}

We mean by the above table that, given a multi-set of weights $\bag{\mbfs_\boldalpha}$, the cycle index polynomial~$Z_n$ computed at the new multi-set~$\bag{\mbfs'_{\boldalpha'}}$ listed in the second column is the normalized moment of the (normalized) Dirichlet distribution listed in the third column. For the sake of notational simplicity, addition and deletion of types was given for the last type in a fixed order for~$\mbfs'$; this is however irrelevant as~$\mbfs'$ is naturally unordered when describing the multi-set~$\bag{\mbfs'_{\boldalpha'}}$; consistently, a typographically more general formulation may be inferred by quasi-exchangeability.
 The last two operations refer instead to the addition and deletion of urns of given weight~$s_i$, when urns of the same weight already exist in the configuration, i.e.~$\alpha_i>0$, and, concerning deletion only, the resulting configuration still has some urns of the same weight as that of the deleted one, i.e.~$\alpha_i>1$. 

The first two correspondences are straightforward; the third one was discussed in Proposition~\ref{p:AggrPolya}. The last two require however some more work, for they rely on a characterization of the \emph{dynamical symmetry algebra} of the function~${}_k\Phi_2$, which we address in the next section.

\subsection{Hypergeometric Lauricella functions of type D}\label{ss:Lauricella}
Most properties of the Dirichlet distribution may be inferred from its moment generating function~${}_k\Phi_2$, a confluent form of the \emph{fourth $k$-variate Lauricella hypergeometric function}~${}_kF_D$ introduced in~\cite{Lau1893} as a multivariate generalization of the more familiar Appell hypergeometric function~$F_1$ and Gauss hypergeometric function~${}_2F_1$.
We refer the reader to~\cite{Ext76} for a classical treatment of hypergeometric functions.

Many remarkable connections between (functionals of) the Dirichlet distribution/process and the function~${}_kF_D$ have been investigated in two series of mostly independent papers~\cite{Reg98,LijReg04},~\cite{KerTsi01,TsiVerYor00, TsiVerYor01} (see~\cite[\S{1\&3}]{LijReg04} for a time-line of the main developments). Both of them almost exclusively rely on integral representations of~${}_kF_D$. Basing instead on the series representation below, we enlarge the picture by providing some combinatorial and algebraic insights.

\paragraph{Some representations} Recall the multiple series and integral representations of $F_D$~\cite[\S2.1]{Ext76} 
\begin{align}
{}_kF_D[a;\mbfb;c;\mbfx]\eqdef& \sum_{\mbfm\in \N_0^k} \frac{\Poch{a}{\length{\mbfm}} \Poch{\mbfb}{\mbfm} \mbfx^\mbfm}{ \Poch{c}{\length{\mbfm}} \mbfm! } & \norm{\mbfx}_\infty<1\\
 \label{eq:LauIntRep}=& \frac{1}{\Beta[a,c-a]} \int_0^1 t^{a-1}(1-t)^{c-a-1}(\uno-t\mbfx)^{-\mbfb} \diff t & \Re c>\Re a>0
\end{align}
and its \emph{confluent} form (or \emph{second $k$-variate Humbert function}~\cite[\emph{ibid.}]{Ext76})
\begin{align}\label{eq:ConflLaur}
{}_k\Phi_2[\mbfb;c;\mbfx]\eqdef \lim_{\eps\rar0} {}_kF_D[1/\eps;\mbfb;c;\eps\mbfx]=\sum_{\mbfm\in \N_0^k} \frac{\Poch{\mbfb}{\mbfm} \mbfx^\mbfm}{ \Poch{c}{\length{\mbfm}} \mbfm! } && \mbfx\in \C^k \fstop
\end{align}

\begin{prop}\label{p:LaurEGF}
The function ${}_k\Phi_2[s\mbft;1;\boldalpha]$ is the exponential generating function of the cycle index polynomials of $\mfS_n$, in the sense that
\begin{align*}
{}_k\Phi_2[\boldalpha; 1; s\mbft]=& \EGF[Z_n[\seq{\mbft^{\compo 1}\cdot\boldalpha, \dotsc, \mbft^{\compo n}\cdot \boldalpha}]](s) & \mbft\in& \R^k, \boldalpha\in \R^k_+ \fstop
\end{align*}

\begin{proof}
Fixing $\length{\boldalpha}=1$ and recalling that ${}_k\Phi_2[\boldalpha;\length{\boldalpha};\mbfs]=\mcL[\Dir[\boldalpha]]$ by~\eqref{eq:LapDir}, we see that equation~\eqref{eq:l:MomDir} provides an exponential series representation for the Laplace transform of the Dirichlet distribution in terms of the cycle index polynomials of symmetric groups, viz.
\begin{align}\label{eq:LaplaceDir2}
\mcL[\Dir[\boldalpha]](\mbfs)=&\sum_{n=0}^\infty \frac{Z_n[\seq{\mbfs^{\compo 1}\cdot\boldalpha, \dotsc,\mbfs^{\compo n}\cdot\boldalpha}]}{n!}
\fstop
\end{align}

Choosing $\mbfs\eqdef s\mbft$ above and using~\eqref{eq:MultiBell} to extract the term $s^n$ from each summand the conclusion follows.
\end{proof}
\end{prop}

\paragraph{Asymptotic distributions in the space of parameters}
Profiting~\eqref{eq:LapDir}, it is possible to compute non-trivial confluent forms of the Humbert function ${}_k\Phi_2$ for appropriate choices of the parameters. We have for instance
\begin{align*}
\lim_{\beta\rar 0^+} {}_k\Phi_2[\beta\boldalpha;\length{\beta\boldalpha};\mbfs]=&\length{\boldalpha}^{-1} \boldalpha\cdot \exp^\compo[\mbfs] \comma\\
\lim_{\beta\rar +\infty} {}_k\Phi_2[\beta\boldalpha;\length{\beta\boldalpha};\mbfs]=& \exp[\length{\boldalpha}^{-1}\boldalpha\cdot\mbfs] \comma
\end{align*}
which results from the limit distributions of $\Dir[\boldalpha]$, obtained in the following proposition.

\begin{prop}\label{l:Asympt}
There exist the narrow limits
\begin{align*}
\lim_{\beta\rar 0^+} \Dir[\beta\boldalpha]=& \length{\boldalpha}^{-1} \sum_{i=1}^{k} \alpha_i\delta_{\mbfe^i}\comma &
\lim_{\beta\rar +\infty} \Dir[\beta\boldalpha]=&\delta_{\length{\boldalpha}^{-1}\boldalpha} \fstop
\end{align*} 

\begin{proof}
Since $\Dir[\boldalpha]$ is moment determinate, it suffices by Stone--Weierstra{\ss} Theorem and compactness of $\Delta^{k-1}$ to show the convergence of its moments. Thus,
\begin{align*}
\mu_n'[\mbfe^i,\beta\boldalpha]\eqdef& \frac{n!}{\Poch{\beta\length{\boldalpha}}{n}} Z_n[\beta \alpha_i \uno] =\frac{1}{\Poch{\beta\length{\boldalpha}}{n}} \sum_{r=1}^n \beta^r \sum_{\boldlambda \vdash_r n} \Mnom{\boldlambda} \alpha_i^{\length{\boldlambda}}\\
\underset{\beta\gg 1}{\sim} &\, \frac{1}{\beta^n \length{\boldalpha}^n}\beta^n \Mnom{n\mbfe^1} \alpha_i^{\length{n\mbfe^1}}= \length{\boldalpha}^{-n}\alpha_i^n \fstop
\end{align*}

The computation for $\beta\longrightarrow 0$ is analogous and therefore omitted.
\end{proof}
\end{prop}

\paragraph{The dynamical symmetry algebra of the second Humbert function}
\emph{Details and proofs for results in this paragraph are provided in the Appendix~\S\ref{app:A}}.

We are now ready to investigate the two operations on our urn model which were left as an open problem in~\S\ref{ss:Polya}. Since the Dirichlet distribution is characterized by its moment generating function, adding/deducting urns of a given type amounts to raise/lower indices of the vector~$\boldalpha$ in~${}_k\Phi_2[\boldalpha;\length{\boldalpha};\mbfs]$. In particular, understanding the mapping of~$\bag{\mbfs_\boldalpha}$ into~$\bag{\mbfs_{\boldalpha+\mbfe^i}}$ amounts to find a (differential) operator~$E_{i}=E_{+i}$, acting on the sole variables~$\mbfs$ (more precisely, not acting on~$\boldalpha$), independent of $\boldalpha$ and such that
\begin{align*}
(E_{i}\, {}_k\Phi_2)[\boldalpha;\length{\boldalpha};\mbfs]= C_{\boldalpha}\, {}_k\Phi_2[\boldalpha+\mbfe^i;\length{\boldalpha+\mbfe^i};\mbfs] \comma
\end{align*}
where $C_\boldalpha$ is some constant, possibly depending on~$\boldalpha$.

For practical reasons, it is convenient to collectively consider the linear span~$\mfl_0$ of the operators~$\set{E_{\pm 1}, \dotsc,  E_{\pm k}}$, endowed with the bracket induced by their composition. In analogy with quantum theory (see~\S\ref{app:A}), the smallest Lie algebra~$\mfl$ extending the bracket on~$\mfl_0$ is termed the \emph{dynamical symmetry algebra} of the function~${}_k\Phi_2[\boldalpha;\length{\boldalpha};\mbfs]$; we postpone an extensive treatment of such construction to~\S\ref{app:A}, where we show (Thm.~\ref{t:dsa}) that $\mfl\cong \mfs\mfl_{k+1}$, the Lie algebra of $(k+1)$-square matrices with vanishing trace, and that the operators above satisfy
\begin{equation}\label{eq:R/LOp}
\begin{aligned}
(E_{i}\, {}_k\Phi_2)[\boldalpha;\length{\boldalpha};\mbfs]=& \alpha_i\, {}_k\Phi_2[\boldalpha+\mbfe^i;\length{\boldalpha+\mbfe^i};\mbfs] \comma \\
 (-E_{-i}\, {}_k\Phi_2)[\boldalpha;\length{\boldalpha};\mbfs]=& (\length{\boldalpha}-1)\, {}_k\Phi_2[\boldalpha-\mbfe^i;\length{\boldalpha-\mbfe^i};\mbfs] \fstop
\end{aligned}
\end{equation}

Combining~\eqref{eq:l:MomDir} with the same argument in the proof of~Proposition~\ref{p:LaurEGF}, this yields e.g.
\begin{align*}
\alpha_i \diff_t^n\restr_{t=0}\, {}_k\Phi_2[\boldalpha+\mbfe^i;\length{\boldalpha+\mbfe^i}; t \mbfs]=\alpha_i \,\mu_n'[\mbfs,\boldalpha+\mbfe^i] \comma
\end{align*}
that is, the moment~$\mu_n^{\prime\, \alpha_i\!\Dir[\boldalpha]}[\mbfs\cdot (\emparg)]$ of the distribution $\alpha_i \Dir[\boldalpha+\mbfe^i]$. This proves the last two correspondences in Table~\ref{tab:A}.

\section{Gamma and Dirichlet--Ferguson random measures}\label{s:DirFer}
Following the same scheme as in the previous section for finite-dimen\-sion\-al distributions, we focus on the infinite-dimen\-sion\-al counterparts.

\paragraph{Notation} Everywhere in the following we denote by $(X,\tau,\mfB(\tau),\sigma)$ a Borel space of \emph{finite} \emph{diffusive} (i.e. atomless) measure \emph{fully supported} on a \emph{perfect} locally compact Polish space~$(X,\tau)$ and respectively by
\begin{itemize}
\item $\Gamma$ the space of configurations;
\item $\Meas$ the vector space of finite signed Radon measures;
\item $\Mp$ the space of non-negative Radon measures;
\item $\Mbp$ the space of non-negative finite (Radon) measures;
\item $\msP$ the space of probability (Radon) measures;
\item $\msK$ the space of discrete finite (Radon) measures;
\item $\overline\msK$ the space of discrete probability (Radon) measures
\end{itemize}
over $(X,\mfB(X))$. Each of the above is endowed with the vague topology (i.e. the weak*-topology) and  the relative (trace) $\sigma$-algebra, denoted by $\mfB$ (e.g. $\mfB(\Mp)$).
Notice that, since the measure~$\sigma$ is assumed to be diffusive and to have full topological support, \emph{perfectness} is in fact redundant; furthermore, it is well-known that any such measure is in fact Radon.
Finally, denote by~$\widehat X$ the (Polish) product space~$X\times \R$, endowed with the product topology, product \mbox{$\sigma$-algebra} and a product measure~$\sigma \otimes \lambda$ with~$\lambda$ satisfying~\eqref{eq:LevyMeas}, and set $\widehat \Gamma\eqdef \Gamma_{\widehat X}$.

\subsection{(Compound) Poisson Type processes}\label{ss:CPT}
The \emph{Poisson process} is defined as the unique (in law) measure-valued process on~$(X,\sigma)$ whose law~$\PP_\sigma$ satisfies the Fourier transform characterization
\begin{align*}
\mcF[\PP_\sigma](f)=&\exp\quadre{\av{e^{if}-1}_\sigma} & f\in \mcC_0 (X)\comma
\end{align*}
where $\mcC_0(X)$ denotes the space of functions on~$(X,\tau)$ vanishing at infinity. The existence of~$\PP_\sigma$ as a cylindrical measure follows by Theorem~\ref{t:BM}, while the existence of the process, as a point process, follows e.g. by \cite[\S2.1]{Kal83}.
When~$X$ is a Riemannian manifold, it is well-known (cf. e.g.~\cite[\S2]{AlbKonRoe98}) that~$\PP_\sigma$ --~constructed via Kolmogorov Extension~-- is in fact a Borel measure on the Polish space~$\Gamma$ (endowed with the vague topology and Borel $\sigma$-algebra) and a Borel measure on~$\Meas$ endowed with the $\sigma$-algebra~$\mfB(\Meas)=\mfB(\Meas, \tau(\Meas, \mcC_c(X)))$. We prove (Thm.~\ref{t:ConvLogPP} below) that under our general assumptions on~$X$, the law~$\PP_\sigma$ is in fact a Radon measure on the $\sigma$-algebra induced on~$\Meas$ by the topology of total variation.

Given a Poisson measure~$\PP_{\sigma\otimes \lambda}$ on~$\widehat X$, the diffusivity of~$\sigma$ implies that
\begin{align*}
\forall \gamma\eqdef \set{(x_i,s_i)}_{i\geq 0}\in \widehat\Gamma\quad \forall i,j\in \N,i\neq j \quad (x_i,s_i)\neq (x_j,s_j) \implies x_i\neq x_j \comma
\end{align*}
so that the map $\mathrm{H} \colon \gamma\mapsto \eta\eqdef \sum_{(x,s)\in \gamma} s\delta_x$ is a measurable isomorphism $\widehat\Gamma\rar \Mp$ on its image endowed with the appropriate trace $\sigma$-algebra. We define the \emph{CPT law with intensity measure~$\sigma$ and L\'{e}vy measure~$\lambda$} as the measure~$\mcR_{\sigma,\lambda}\eqdef \mathrm{H}_\pfwd \PP_{\sigma\otimes\lambda}$, satisfying the Fourier transform characterization
\begin{align}\label{eq:CharCPT}
\mcF[\mcR_{\sigma,\lambda}](f)=\exp\quadre{\int_{\widehat X} (e^{i s f(x)}-1)\diff\lambda(s)\diff\sigma(x)} \fstop
\end{align}

We will not discuss the properties of processes with law~$\mcR_{\sigma,\lambda}$. In this respect, we only notice that, in the case when~$\mcR_{\sigma,\lambda}$ may be regarded as supported on a nuclear dual space, the existence of a (unique in law) additive process with c\`adl\`ag paths whose law satisfies the above characterization follows by results in~\cite[\S3]{Ust84}.

\smallskip

\emph{Suitable definitions concerned with more general L\'{e}vy measures may be found in~\mbox{\cite{Lyt03, KonLytVer15, HagKonLytVer14}}.}

\paragraph{Cumulants and CPT laws} It is relevant to many applications wether the computation of moments can be made easier than that of cumulants or vice versa, so to retrieve, via duality, the ones from the others. Cumulants turn out to be preeminent over moments when the Fourier/Laplace transforms of the distributions in question are \emph{of exponential type}. In that case, the cumulant generating function is immediately seen to be \emph{linear}, rather than exponential (as the moment generating function), in the interesting quantities.

Courtesy of the \emph{L\'{e}vy--Khintchine formula}, a rich class of examples is provided by those (\emph{completely}) \emph{random measures} that are laws of L\'{e}vy processes, whose most representative instances --~together with Brownian motion~-- are (C)PT processes. For these processes, (analytically of exponential type by~\eqref{eq:CharCPT} as soon as their Laplace transform is finite), a straightforward application of~\eqref{eq:CumFourier} to~\eqref{eq:CharCPT} yields
\begin{align*}
\kappa_{\mbfm}^{\mcR_{\sigma,\lambda}}[\mbff]=
\mu_{\length{\mbfm}}^{\prime \, \lambda}[1] \,\mu_\mbfm^{\prime \, \sigma}[\mbff] \comma
\end{align*}
hence in particular, for~$\PP_\sigma=\mcR_{\sigma,\delta_1}$, the formula
\begin{align}\label{eq:PoiCumGen}
\kappa_{\mbfm}^{\PP_{\sigma}}[\mbff]=\mu_\mbfm^{\prime \, \sigma}[\mbff] \comma
\end{align}
where we take care to interpret~$\mbff=\seq{\scalar{\emparg}{f_1},\dotsc,\scalar{\emparg}{f_k}}$ in the left-hand side.
The moments of CPT laws, usually the interesting quantities in the applications, are then retrieved via duality from the cumulants computed above (see e.g.~Prop.~\ref{p:MomGP} and Rem.~\ref{r:MomPoi} below for two simple examples).

\begin{rem}[Convolution logarithm of Poisson laws]
Equation~\eqref{eq:PoiDelta} suggests a different way to look at the previous formula for cumulants, involving the convolution logarithm of~$\PP_\sigma$.
More precisely, let~$\PP_\sigma$ be a Poisson law with finite intensity~$\sigma$ and set~\mbox{$\delta\colon x\mapsto \delta_x$}. Fixing~$\mbff=\mbff_k$ a vector of continuous functions on~$X$ one has
\begin{equation}\label{eq:hint}
\begin{aligned}
\mu_\mbfm^{\prime\, \delta_\pfwd\sigma}[\mbff]=&\int_{\Mp} \diff\delta_\pfwd \sigma (\eta) \prod_{i=1}^k\scalar{\eta}{f_i}^{m_i}= \int_X \diff\sigma(x) \tonde{\prod_{i=1}^k\scalar{\emparg}{f_i}^{m_i}\circ\delta} [x]\\
=&\int_X \diff\sigma(x) \prod_{i=1}^k\scalar{\delta_x}{f_i}^{m_i}= \int_X \diff\sigma(x) \,\mbff^\mbfm[x]=\mu_\mbfm^{\prime \, \sigma}[\mbff] \comma
\end{aligned}
\end{equation}
where again we take care to interpret~$\mbff=\seq{\scalar{\emparg}{f_1},\dotsc,\scalar{\emparg}{f_k}}$ in the left-hand side. Together with~\eqref{eq:PoiCumGen} this shows that~$\kappa_\mbfm^{\PP_\sigma}[\mbff]=\mu^{\prime \, \delta_\pfwd \sigma}_\mbfm[\mbff]$.
A comparison of the latter equality with~\eqref{eq:PoiDelta} is reason for the claim that~$\ln^*\PP_\sigma=\delta_\pfwd \sigma$. A rigorous proof of the statement in full generality is however not immediate.
\end{rem}

\begin{thm}[Convolution logarithm of Poisson measures]\label{t:ConvLogPP}
Let~$(X,\tau, \mfB(\tau), \sigma)$ be a finite diffusive Borel measure space over a locally compact Polish space. Then
\begin{align*}
\PP_\sigma=\exp^*[\delta_\pfwd \sigma] \fstop
\end{align*}

\begin{proof}
\emph{For the sake of clarity, in this proof we state topologies and $\sigma$-algebras explicitly, withholding our usual conventions established in the beginning of the present section.}

Let $E\eqdef \mcC_c(X)$ denote the space of continuous compactly supported functions on~$X$ endowed with the supremum norm. By assumptions on $(X,\tau)$, the topological dual space~$E'$ may be identified with the Banach space~$\Meas$ endowed with the total variation norm. Denote further by~$E'_s$ the said space when endowed with its strong (i.e. norm) topology and the associated Borel $\sigma$-algebra, resp. by~$E'_{w^*}$ the same space endowed with the weak* topology~$\tau(E',E)$ and the associated Borel $\sigma$-algebra.
Recall as a standard fact in functional analysis that~$E'_{w^*}$ is completely regular Hausdorff.

Throughout the proof fix~$c\eqdef \abs{\sigma}$.

\smallskip

\emph{Step 1 ($\delta_\pfwd \sigma$ is $E'_s$-Radon).}
Recall as standard fact in measure theory that every finite Borel measure on a Polish space is Radon, hence~$\sigma$ is a Radon measure on~$X$.

The Dirac embedding~$\delta\colon X\rar \msP\subset E'_{w^*}$ is readily seen to be continuous (i.e. continuous with respect to the vague topology on~$\msP$), hence $(X,\sigma)$-$E'_{w^*}$-almost-continuous (in the sense of~\cite[411M]{Fre01}) by \cite[418I]{Fre01}. Every almost-continuous (Hausdorff-valued) image of a (locally) finite Radon measure is itself Radon~(\cite[418I]{Fre01}), thus
the measure~$\delta_\pfwd \sigma$ is $E'_{w^*}$-Radon, hence $\tau(E',E)$-smooth (in the sense of~\cite[\S{I.3.2}]{VakTarCho87}, also $\tau$-additive, $\tau$-regular) by~\cite[\S{I.3.2}, Prop.~3.1.c]{VakTarCho87}.

The measure~$\sigma$ is tight on~$X$ by Ulam's Theorem~(\cite[\S{I.3.2}, Thm.~3.1]{VakTarCho87}), hence, for every~$\eps>0$, there exists a compact set~$K_\eps\subset X$ such that $\sigma K_\eps>c-\eps$. By $E'_{w^*}$-continuity of~$\delta$, the set~$\delta[K_\eps]$ is itself compact, hence the measure~$\delta_\pfwd\sigma$ is itself $E'_{w^*}$-tight on~$E'_{w^*}$. Recall now that every $\tau$-smooth tight Baire measure on a completely regular Hausdorff space admits a unique Borel Radon extension~(\cite[\S{I.3.5}, Thm.~3.3]{VakTarCho87}); since the Baire $\sigma$-algebra of $E'_s$ coincides with the Borel $\sigma$-algebra of $E'_{w^*}$ (by~\cite[\S{I.2.2}, Prop.~2.7]{VakTarCho87}), the measure~$\delta_\pfwd\sigma$ uniquely extends to a Radon measure on~$E'_s$, denoted in the same way.

\emph{Step 2 (well-posedness of~$\exp^*$).}
In analogy with~\eqref{eq:ConvExp}, notice that~$\delta_\pfwd \sigma=c\delta_\pfwd \n\sigma$.

The convolution of two $\tau$-smooth (resp. Radon) probability measures on a topological vector space is itself a $\tau$-smooth (resp. Radon) probability measure by~\cite[\S{I.4.2}, Prop.~4.4]{VakTarCho87}, hence all of the convolution powers~$(\delta_\pfwd \n\sigma)^{*k}$ are Radon probability measures on~$E'_s$, in fact supported on~$E^{\prime\,+}=\Mbp(X)$ by non-negativity of~$\sigma$ (whence of~$\n\sigma$).
Denote by~$p_n[t]\eqdef \sum_{0\leq k\leq n} t^k/k!$ the truncation of the exponential power series up to~$n$ and by~$p_n^*$ its convolution analogue. The sequence of measures~$\varpi_n^\sigma\eqdef p_n[c]^{-1} p_n^*[\delta_\pfwd \sigma]$ satisfies
\begin{align*}
\norm{\varpi_n^\sigma-\varpi_m^\sigma}_\TV\leq& \abs{p_n[c]^{-1}-p_m[c]^{-1}} \norm{p_n^*[\delta_\pfwd \n\sigma]}_\TV+ p_n[c]^{-1}\norm{p_n^*[\delta_\pfwd \n\sigma]-p_m^*[\delta_\pfwd \n\sigma]}_\TV\\
\leq& \abs{p_n[c]-p_m[c]}+ \norm{p_n^*[\delta_\pfwd \n\sigma]-p_m^*[\delta_\pfwd \n\sigma]}_\TV\\
\leq& \abs{p_n[c]-p_m[c]}+ \abs{p_n[1]-p_m[1]}\comma
\end{align*}
hence it is a fundamental sequence in~$\Meas(E'_s)$, the space of finite Borel (not necessarily Radon) measures on~$E'_s$, endowed with the total variation norm. By completeness of the former  there exists the limit~$\nlim \varpi_n^\sigma\defeq\exp^*[\delta_\pfwd \sigma]$, itself a probability measure on~$E'_s$.

\smallskip

\emph{Step 3 (L\'{e}vy Continuity and identification). } Since the narrow topology on~$\Meas(E'_s)$ is coarser than the topology of total variation, the \emph{countable} set~$\Pi\eqdef \set{\varpi_n^\sigma}_{n\geq 0}$ is narrowly (sequentially) compact. As a consequence, by L\'{e}vy's Continuity Theorem in the form~\cite[Thm.~IV.3.1]{VakTarCho87}, as soon as the Fourier transforms~$\mcF[\varpi_n^\sigma](\emparg)$ converge pointwise on~$E$ to a functional $\varphi$, there exists a unique narrow (hence total variation) $E'_s$-Radon cluster point~$\varpi^\sigma$ of~$\Pi$, and such that~$\mcF[\varpi^\sigma]=\varphi$.
Now, usual properties of the Fourier transform of probability measures with respect to the convolution of measures (see~\cite[Prop.~IV.2.3]{VakTarCho87}) yield
\begin{align*}
\mcF[\varpi_n^\sigma](f)=&p_n[c]^{-1} \mcF[p_n^*[\delta_\pfwd \n\sigma]](f)=p_n[c]^{-1} p_n[\mcF[\delta_\pfwd \n\sigma]](f) && f\in E\comma
\intertext{where, reasoning as in~\eqref{eq:hint} above (and retaining once more the abuse of notation on linear functionals),}
\mcF[\delta_\pfwd \n\sigma](f)\eqdef &\int_E e^{i\scalar{\eta}{f}} \diff\delta_\pfwd \n\sigma(\eta)=\mcF[\n\sigma](f) && f\in E \fstop
\end{align*}

Since~$\nlim p_n[a]^{\pm 1}=e^{\pm a}$ for any constant~$a$, choosing respectively~$a=c$ and $a=\mcF[\delta_\pfwd \n\sigma](f)$ we have
\begin{align*}
\varphi\eqdef\nlim\mcF[\varpi_n^\sigma](f)=&\exp[\mcF[\sigma](f)-c]= \exp\quadre{\int_X (e^{i f(x)}-1) \diff\sigma(x)} && f\in E \fstop
\end{align*}

As a consequence of the above reasoning,~$\varpi^\sigma\eqdef\exp^*[\delta_\pfwd \sigma]$ satisfies
\begin{align*}
\mcF[\varpi^\sigma](f)=&\mcF[\PP_\sigma](f) && f\in E \comma
\end{align*}
hence, by~\cite[\S{IV.2.1},~Cor.~2]{VakTarCho87} we get~$\exp^*[\delta_\pfwd \sigma]=\mcP_\sigma$ as measures on~$E'_{w^*}$. We can thus regard~$\exp^*[\delta_\pfwd \sigma]$ as defined above as a Radon extension of~$\PP_\sigma$.
\end{proof}

A generalization of the statement to the case of~$\sigma$-finite intensity~$\sigma$ follows by approximation via a suitable generalization of L\'{e}vy's Continuity Theorem (e.g. de Acosta's Theorem~\cite[\S{IV.3.2}, Thm.~3.3]{VakTarCho87}). The proof may be greatly simplified if one assumes that~$E'_{w^*}$ is contained (both as a topological and measurable space) in~$\mcD'$ for some `nice' nuclear Gel'fand rigging~$\mcD\subset L^2_\sigma \subset \mcD'$, e.g. in the case of the standard Schwartz rigging~$\mcS(\R^d)\subset L^2(\R^d)\subset \mcS'(\R^d)$.
\end{thm}

\subsection{Gamma and Dirichlet--Ferguson random measures}\label{ss:PGDFrms}

We are now able to introduce the Gamma and Dirichlet--Ferguson processes, and the associated laws. While we are mainly interested in the measures, the processes provide useful heuristics and characterizations.

\paragraph{Notation} In addition to the notation already established in~\S\ref{ss:CPT}, we set the following.
For~$\mu$ in~$\Mb\setminus \set{0}$, we denote further by~$\ev_f\colon \mu\mapsto \scalar{\mu}{f}$ the \emph{evaluation map of $f$}, also setting $\ev_A\eqdef\ev_{\car_A}$ and $\abs{\eta}\eqdef \ev_X\eta$; by~$\n\mu$ the corresponding \emph{normalized} measure~$\mu/\abs{\mu}$ in $\msP$ and by~$N\colon \mu\mapsto \n\mu$ the \emph{normalization map}. For~$\eta\in \msK$ let further $\tau(\eta)\eqdef\set{x\in X\mid \eta\set{x}>0}$ be the \emph{pointwise support of $\eta$} and notice that $\tau(\eta)=\tau(\n\eta)$.
Finally, write~$X^{\times n}$ for the \mbox{$n$-fold} cartesian product of $X$ endowed with the product topology and $\sigma$-algebra, $\mbfx\eqdef(x_1,\dotsc,x_n)$ for a tuple in $X^{\times n}$, 
\begin{align*}
\widetilde X^{\times n}\eqdef\set{\mbfx\in X^{\times n}\mid x_i\neq x_j, i,j\in[n], i\neq j}
\end{align*}
and~$X^{\odot n}\eqdef \quotient{X^{\times n}}{\mfS_n}$, the group~$\mfS_n$ acting on~$X^{\times n}$ by permutation of co\"ordinates.
%
\paragraph{The Gamma measure} Momentarily, let~$X=\R^d$. The \emph{Gamma process} $\gamma_{\theta,k}$ with intensity~$\sigma$, \emph{shape} $\theta>0$ and \emph{scale}~$k>0$ is the unique $\Mp$-valued process on $\widehat X$ with L\'evy measure 
\begin{align*}
\diff\lambda_{\theta,k}(s)\eqdef \theta s^{-1}e^{-ks}\diff s \fstop
\end{align*}

The corresponding infinitely divisible law is the Gamma distribution $\Gam[\theta,k]$ and the corresponding intensity measure on $\widehat X$ is $\kai_{\sigma,\theta,k}\eqdef\sigma\otimes\lambda_{\theta,k}$.
As an equivalent definition~\cite[3.1]{HagKonLytVer14}, the \emph{Gamma measure} $\GP_{\sigma,\theta,k}\eqdef \law \gamma_{\theta,k}$ satisfies the Laplace transform characterization
\begin{align}\label{eq:LapGam}
\mcL[\GP_{\sigma,\theta,k}](f)=&\exp\quadre{-\theta\av{\quadre{\log\quadre{1-f}}}_\sigma} &  f\in \mcC_c(X)&\mid f<1\fstop
\end{align}

Existence (also in the general case when~$X$ is a Polish space) follows from that of CPT laws in the previous section (see also \cite{Kin67}). If $\gamma_k$ is the Gamma process with scale parameter~$k$ then~$\gamma_k=k\gamma_1$ in distribution, thus, following \cite{TsiVerYor01}, we restrict to the case $k=1$.
It is also apparent from \eqref{eq:LapGam} that $\GP_{\sigma,\theta,1}=\GP_{\theta\sigma, 1,1}$, thus, without loss of generality up to multiply~$\sigma$ by~$\theta$, we may further restrict ourselves to consider the random measure $\GP_\sigma\eqdef \GP_{\sigma, 1, 1}$, which is in turn characterized by the \emph{Mecke identity} (see~e.g.~\cite[3.2]{HagKonLytVer14})
\begin{equation}\label{eq:Mecke}
\int_{\msK}\diff \GP_\sigma(\eta)\int_X\diff\eta(x)F(\eta,x)=\int_{\msK}\diff\GP_\sigma(\eta)\int_{\widehat X} \diff \sigma(x) \diff s\, e^{-s} F(\eta+s\delta_x) \fstop
\end{equation}

Furthermore, it follows from Lukacs characterization of the Gamma distribution (see above) that the \emph{total charge} of the Gamma process is independent of the \emph{normalized Gamma} (or~\emph{Diri\-chlet--Ferguson}) \emph{process} 
\begin{align*}
\n\gamma_{\theta,k}\eqdef \frac{\gamma_{\theta,k}}{\abs{\gamma_{\theta,k}}}
\end{align*}
and satisfies (see Lemma~\ref{l:DistribTotMassG} below)
\begin{align*}
\abs{\gamma_{\theta,k}}\sim \Gam[\theta \abs{\sigma},1] \fstop
\end{align*}

The Gamma process is in fact characterized by this property (cf~\cite[2.1]{TsiVerYor01}), i.e. if for a L\'evy process the total charge and the normalized process are independent, then the processes is a Gamma process (compare with Lukacs' characterization of the Gamma distribution, see~\S\ref{ss:PoiGamDir} above).
Finally, as the Gamma process is a point process, the set of discrete Radon measures has full $\GP_\sigma$-measure.

\begin{lem}[Distribution of the total mass, properties of the support]\label{l:DistribTotMassG}
Let~$\sigma=\beta\n\sigma$.
Then, the total mass $\abs{\eta}$ is $\GP_\sigma(\diff\eta)$-a.e. $\Gam\quadre{\beta,1}$-distributed. Furthermore,~$\eta$ has full topological support.

\begin{proof}
\emph{Proofs of both statements are well-known, especially in the infinite-dimensional analysis Literature, (see e.g.~\cite{TsiVerYor01,Kin93,HagKonLytVer14}) and generally based on properties of the Gamma process. We give here a proof based on the Mecke identity~\eqref{eq:Mecke}, hence not relying on the existence of the process as a point process.}

For $L\subset [0,+\infty]$ let $\msK_L\eqdef \msK\cap\ev_X^{-1}(L)$, $\msK_u\eqdef \msK_{\set{u}}$ and set $\GP\eqdef \GP_{\beta\n\sigma}$ for notational purposes.
Let~$\set{\GP_u}_u$ be the (Rokhlin) disintegration of $\GP$ with respect to the total mass~$\abs{\emparg}$ and notice that 
\begin{align*}
\GP(\set{0})=\GP(\msK_\infty)=0
\end{align*}
(the latter by finiteness of $\sigma$).
For a measurable subset~$A\subset X$ set
\begin{align*}
C_A\eqdef \set{\eta\in \msK\mid \tau(\eta)\subset A}=\set{\eta\in\msK\mid \ev_A[\eta]=\abs{\eta}}\in\mfB(\msK)
\end{align*}
and notice that it is measurable since so is~$\ev_A$.
Choosing~$F(\eta,x)\eqdef \car_{\msK_{[0,t]}}\car_{C_A}$ in the Mecke identity~\eqref{eq:Mecke}
yields
\begin{align*}
\int_\msK\diff\GP(\eta) \abs{\eta} \car_{\msK_{[0,t]}}(\eta)\car_{C_A}(\eta)=&\\
=\int_\msK\diff\GP(\eta)&\int_X\diff \beta \n\sigma(x)\int_0^\infty\diff s \, e^{-s} \car_{\msK_{[0,t-s]}}(\eta)\car_{[0,t]}(s)\car_{C_A}(\eta)\\
\int_0^\infty\diff u \int_{\msK_u}\diff\GP_u(\eta) \abs{\eta} \car_{[0,t]}(\abs{\eta})\car_{C_A}(\eta)=&\\
=\beta\int_0^\infty\diff u &\int_{\msK_u}\diff\GP_u(\eta) \int_0^t\diff s \, e^{-s} \car_{C_A}(\eta)\car_{[0,t-s]}(\abs{\eta})\\
\int_0^\infty\diff u \, u \car_{[0,t]}(u) \int_{\msK_u}\diff\GP_u(\eta) \car_{C_A}(\eta)=&\\
=\beta\int_0^\infty\diff u &\int_0^t\diff s\, e^{-s} \car_{[0,t-s]}(u) \int_{\msK_u}\diff\GP_u(\eta) \car_{C_A}(\eta)\\
\int_0^t \diff u \, u g_A(u)=&\beta \int_0^t \diff s\,e^{-s}\int_0^{t-s}\diff u \, g_A(u) \comma
\end{align*}
where we set $g_A(u)\eqdef \av{\car_{C_A}}_{\GP_u}$.
Differentiating in $t$ twice yields 
\begin{align}\label{eq:ODEGamma}
t\dot g_A(t)=\tonde{\beta-1-t}g_A(t) \fstop
\end{align}

If $A$ is $\sigma$-co-negligible, monotonicity of measures implies 
\begin{align}\label{eq:ODEGammaIn}
1=\GP(\msK)=\lim_{t\rar\infty} \GP(\msK_{[0,t]} \cap C_A)\comma
\end{align}
and the unique solution $g_A=\Gam[\beta,1]$ to~\eqref{eq:ODEGamma} with initial condition~\eqref{eq:ODEGammaIn} is the required distribution of the total charge.
Assume now that $A$ is not of full $\sigma$-measure. Since $X$ is perfect (non-empty), there exists an \emph{infinite} family $\set{A_n}_{n\geq 0}$ of measurable subsets $A_n \subset X$, such that $A_0\eqdef A$, $\sigma(A_n)=\sigma(A_0)$ independently of $n$ and $A_i\cap A_j$ is $\sigma$-non-negligible for $i\neq j$. The corresponding sets $C_{A_n}$ form an infinite \emph{disjoint} family $\mcC$ in $\mfB(\msK)$. Since $\GP(C_{A_n})$ only depends on $\sigma(A_n)$ we have $\GP(C_{A_n})=c_0$ independently of $n$, hence the disjointness and cardinality of $\mcC$ together with the fact that~$\GP$ is a probability measure imply that $c_0=0$. This proves $\supp\eta=X$ for a set of full $\GP$-measure.
\end{proof}
\end{lem}

Now that we have the distribution of the total charge, the aforementioned independence allows to decompose~$\GP_\sigma$ as a product measure 
\begin{align}\label{eq:ProdMeas}
&\begin{aligned}
\GP_{\beta\n\sigma} =&\DF_\n\sigma\otimes \Gam[\beta,1](s)\diff s\\
N_\pfwd \GP_{\sigma}=& \DF_{\sigma}
\end{aligned}
& \sigma&=\beta\n\sigma\comma
\end{align}
where $\DF_\n\sigma$ is the \emph{simplicial part} and $\Gam[\beta,1](s)\diff s$ the \emph{conic part} of the measure (cf.~\cite{TsiVerYor00,TsiVerYor01}). As detailed in~\cite{TsiVerYor00}, the simplicial part is the law of the Dirichlet--Ferguson process below.

\paragraph{The Dirichlet--Ferguson measure}
Set
\begin{align*}
&\begin{aligned}
\ev_{\mbfX}^\beta\colon \msP &\longrightarrow \beta\Delta^{k-1}\subset \R^{k}\\
\nu &\longmapsto \beta\seq{ \nu X_1, \dotsc, \nu X_{k}}
\end{aligned} & \mbfX\eqdef \seq{X_1, \dotsc, X_{k}} \subset \mfB(X) \fstop 
\end{align*}

The \emph{Dirichlet--Ferguson measure $\DF_{\n\sigma}^\beta= \DF_{\beta\n\sigma}$ of parameter $\beta$ and reference measure $\sigma$}, introduced in~\cite[\S1,~Def.~1]{Fer73} (see also~\cite[\S2]{Set94} for an explicit construction), is the unique random measure on $\msP$ with intensity~$\sigma$ such that
\begin{align}\label{eq:PartX}
\tonde{\ev_{\mbfX}^1}_\pfwd \DF_{\n\sigma}^\beta\sim& \Dir\quadre{\ev_{\mbfX}^\beta\quadre{\n\sigma}} &
\begin{aligned}
\mbfX\eqdef& \seq{X_1, \dotsc, X_{k}} \subset \mfB(X) \\
X=&\sqcup_i X_i
\end{aligned}
\fstop
\end{align}

More explicitly (cf.~e.g.~\cite[6.1]{Stu11}), for every measurable function $u\colon \Delta^{k-1}\rar\R$
\begin{equation}\label{eq:d:DirFer}
\int_{\msP} \tonde{u\circ\ev_{\mbfX}^1} \quadre{\eta} \diff\DF^\beta_{\n\sigma}(\eta)=\int_{\Delta^k} u(\mbfy) \Dir\quadre{\ev_{\mbfX}^\beta\quadre{\n\sigma}}(\mbfy) \diff \mbfy \fstop
\end{equation}

\emph{Together with that of Dirichlet--Ferguson, the measure is known by a number of different names, including Dirichlet (e.g.~\cite{LijPru10}), Poisson--Dirichlet (e.g.~\cite{Kin06}) and Fleming--Viot with parent-independent mutation (e.g.~\cite{Fen10})}.

\smallskip

Existence was originally proved in~\cite{Fer73} by means of Kolmogorov Extension Theorem. Since~$X$ is Polish, it is in fact sufficient by the Portmanteau Theorem~\cite[\S8.2]{Bog07} to consider~$u$ continuous with $\norm{u}_{\infty}<1$ and~$\mbfX$ consisting of continuity sets for~$\sigma$ (cf.~e.g.~\cite[6.4]{Stu11}).
The measure, concentrated on $\n\msK$ by~\eqref{eq:ProdMeas} and Lemma~\ref{l:DistribTotMassG} (also cf.~\cite[property~\textbf{P2}]{Set94}), is the law of the normalized Gamma process~$\n\gamma_{1,1}$.

Observe the notational analogy $\length{\boldalpha}=\abs{\sigma}$; in light of Remark~\ref{r:Models}, the characterization~\eqref{eq:d:DirFer} is equivalent to the requirement that
\begin{align*}
(g_\pfwd)_\pfwd \DF_\sigma= \Dir[g_\pfwd \sigma]
\end{align*}
for arbitrary~$g\colon X\rar [k]$ constant on each $X_i$ (and such that every~$\sigma$-representative of $g$ is surjective). It is in fact straightforward to check that for any measurable space~$Z$ satisfying the same hypotheses as~$X$ (see above) and for any measurable map~$f\colon X\rar Z$, such that~$f_\pfwd \sigma$ is again diffusive, it holds
\begin{align*}
(f_\pfwd)_\pfwd\DF_\sigma=\DF_{f_\pfwd\sigma}
\end{align*}
(cmp. the same result for~$\PP_\sigma$,~\cite[\S2.3]{Kin93}).

\begin{rem}[Contractions, quasi-exchangeability and symmetry]\label{r:CQES}
Recall that a random measure~$\mcR$ on~$X$ is termed $\nu$-\emph{symmetric} or \emph{symmetric} with respect to a measure~$\nu$ on~$X$ if for every realization~$\rho$ of~$\mcR$ and every $\nu$-preserving map $f\colon X\rar X$ (i.e. such that~$f_\pfwd \nu=\nu$) it holds that~$f_\pfwd\rho=\rho$ in distribution. From~\cite[9.0]{Kal83} we know that $\mcR$ is $\nu$-symmetric iff the random vector $\rho^\compo\mbfX\eqdef\seq{\rho X_1,\dotsc, \rho X_k}$, with $\mbfX$ as in~\eqref{eq:PartX}, is distributed as a function of~$\nu^\compo\mbfX$ iff~$\rho^\compo\mbfX$ is exchangeable for every $\mbfX$ such that $\nu^\compo\mbfX=\nu X_1 \uno$. 
Choosing~$\mcR=\DF_\sigma$ and~$\nu=\sigma$ yields the symmetry of~$\DF_\sigma$ by (quasi-)exchangeability of~$\Dir[\sigma X_1 \uno]$.

Essentially because of this symmetry, \emph{contractions} analogous to the one in~\eqref{eq:AddContract} appear, elsewhere in the following, in different forms. The underlying common principle is that a contraction is but an operation on subsets of a set of indices (e.g.~$[k]$ in the description of the aggregation property~\eqref{eq:AggrGen}), thus naturally encoded by set partitions, and in particular by the set partition 
\begin{align*}
\vec\mbfL(\boldlambda)\eqdef \set{\set{1,\dotsc,\lambda_1}, \set{\lambda_1+1,\dotsc, \lambda_1+2\lambda_2}, \dotsc} \comma
\end{align*}
appearing e.g. in the statement of Rota--Wallstrom Theorem~\cite[6.1.1]{PecTaq11} (see also~\cite[\S6 \emph{passim}]{PecTaq11}). Wether this operation may be described via integer rather than set partitions depends on the possibility to freely permute the set of the indices, in the case of the Dirichlet distribution a consequence of quasi-exchangeability.
\end{rem}

Courtesy of the finite-dimensional projections~\eqref{eq:d:DirFer}, the moments of $\DF_\sigma$ may be inferred by~\eqref{eq:MomDir} and Theorem~\ref{l:lemma}. Indeed, if~$\mbfX$ is a measurable partition of~$X$, $s$ is a simple function attaining the value~$s_j$ on~$X_j$ and $\boldalpha\eqdef \ev_{\mbfX}[\sigma]$, then
\begin{align*}
\int_\n\msK \ev_s[\n\eta]^n \diff \DF_\sigma[\n\eta]=& \int_{\Delta^{k-1}} (\mbfs\cdot \mbfx)^n\Dir[\boldalpha](\mbfx)\diff \mbfx= \mu'_n[\mbfs,\boldalpha]
\fstop
\end{align*}

A simple density argument, Dominated Convergence Theorem and Theorem~\ref{l:lemma} yield for any bounded measurable $f\colon X\rar \R$
\begin{align}\label{eq:MomDF}
\mu^{\prime\,\DF^\beta_\sigma}_n[f]=& \frac{n!}{\Poch{\beta}{n}} Z_n\quadre{\av{(f\uno)^{\compo\vec\mbfn}}^\compo_\n\sigma} \fstop
\end{align}

As a consequence, fixing $\beta\eqdef \abs{\sigma}=1$ and arguing along the same lines as above for the Laplace transform yields, up to rearranging series,
\begin{align}\label{eq:LapDF}
\mcL[\DF^1_\sigma](\ev_f) =& \sum_{\mbfm} \prod_{i=1}^\infty \frac{\av{f^i}_\sigma^{m_i}}{i^{m_i}m_i!}\defeq \Phi [\sigma,1,f] \comma & f\in\mcC(X), \quad \norm{f}_\infty<1,
\end{align}
where $\mbfm$ is any non-negative finite multi-index in an arbitrary number of co\"ordinates.

\begin{prop}[Asymptotic behavior of the Dirichlet--Ferguson measure]\label{p:AsymptDF}
When $X$ is \emph{compact}, the asymptotic behavior in Proposition~\ref{l:Asympt} generalizes to the infinite-dimensional case as
\begin{align}\label{eq:AsymptDF}
\DF^0_\n\sigma\eqdef& \lim_{\beta\rar0} \DF^\beta_\n\sigma= \delta_\pfwd\sigma\comma & \DF^\infty_\n\sigma\eqdef&\lim_{\beta\rar\infty} \DF^\beta_\n\sigma=\delta_\n\sigma \comma
\end{align}
where $\delta\colon X\ni x\mapsto \delta_x\in \msP$ denotes the canonical embedding of $X$ into $\msP$.

\begin{proof} By Prokhorov Theorem 
the space $\msP$ is compact as well, hence there exist some weak cluster points~$\DF_\n\sigma^0$, resp.~$\DF_\n\sigma^\infty$, for the family $\{\DF_\n\sigma^\beta\}_\beta$ when $\beta$ tends to~$0$, resp.~$+\infty$. 

Uniqueness and representation are now straightforward from the finite-dimensional projection~\eqref{eq:d:DirFer} by approximation with cylinder functions, which are uniformly dense in the space of continuous functions by compactness of $\msP$ via Stone--Weierstra{\ss} Theorem.
\end{proof}
\end{prop}

Analogous asymptotic formulae for $\GP_\sigma$ readily follow from~\eqref{eq:ProdMeas}, viz.
\begin{align*}
\GP_0\eqdef&\lim_{\beta\rar 0} \GP_{\beta\n\sigma}=\delta_\pfwd\n\sigma \otimes \Gam[1,1] & \GP_\infty\eqdef&\lim_{\beta\rar \infty} \GP_{\beta\n\sigma}=\delta_\n\sigma \otimes \Gam[1,1] \fstop
\end{align*}

\begin{rem}[A Gibbsean interpretation]\label{r:MeccaStat} The asymptotic distributions~\eqref{eq:AsymptDF} may be interpreted --~at least formally~-- in the framework of statistical mechanics. In order to establish some lexicon, denote by~$G_\beta\eqdef Z_\beta^{-1}\exp[-\beta H]$ the \emph{Gibbs measure} of a physical system at inverse temperature~$\beta$ driven by a \emph{Hamiltonian}~$H$ with \emph{partition function}~$Z_\beta\eqdef \av{\exp[-\beta H]}$ and \emph{Helmholtz free energy}~$F_\beta\eqdef -\beta^{-1}\ln Z_\beta$. It was argued in~\cite[\S3.1]{vReStu09} that, heuristically,
\begin{align*}
\diff \DF^\beta_\n\sigma[\n\eta] = \frac{e^{-\beta \, S[\n\eta]}}{Z_\beta}\diff \DF^*_\n\sigma[\n\eta] \comma
\end{align*}
where $S$ is now an \emph{entropy functional} (rather than an energy functional), $Z_\beta$ a normalization constant and $\beta$ plays the role of the \emph{inverse temperature}. Here, $\DF^*_\n\sigma$ denotes a \emph{non-existing} (!) uniform distribution on $\msP$.
Borrowing again the terminology, one can say that for small~$\beta$ (i.e. large \emph{temperature}), the system \emph{thermalizes} towards the ``uniform'' distribution $\DF^0_\n\sigma$ induced by the (normalized) reference measure $\n\sigma$ on the base space, while for large $\beta$ it \emph{crystallizes} to $\n\sigma$, so that all randomness is lost.
\end{rem}

Next we turn to the moments of the Gamma measure. These are well-known, although never explicitly described in terms of Bell or cycle index polynomials. 

\begin{prop}[Moments of the Gamma measure]\label{p:MomGP}
The following identity holds
\begin{align}\label{eq:MomGP}
\mu_n^{\prime\, \GP_{\sigma}}[f]=n!\, Z_n\quadre{ \av{(f\uno)^{\compo\vec\mbfn}}^\compo_{\sigma}} \fstop
\end{align}

Furthermore,~$\GP_\sigma$ is analytically of exponential type.

\begin{proof}
The formula for the moments follows from~\eqref{eq:MomDF} using the product decomposition~\eqref{eq:ProdMeas} along with the moments~\eqref{eq:MomGam} of the Gamma distribution $\Gam[\beta,1]$.
Once established that $\GP_\sigma$ is \emph{analytically} of exponential type, the moments may also be computed directly (cf.~\cite{KonDaSStrUs98,KonLyt00}) by differentiating the Laplace transform~\eqref{eq:LapGam} along with~\eqref{eq:MomBell}, or trivially from the cumulants, in turn immediately computable as
\begin{align}\label{eq:CumGam}
\kappa^{\GP_{\sigma}}_n[f]\eqdef \diff_t^n\restr_{t=0} \ln\mcL[\GP_{\sigma}](tf)=\Gamma[n] \av{f^n}_{\sigma} \comma
\end{align}
whence~\eqref{eq:CumToMom} and~\eqref{eq:CI} together yield~\eqref{eq:MomGP}.

The fact that the Laplace transform~$\mcL[\GP_\sigma](tf)$ of $\GP_\sigma$ is analytic in $t$ was originally argued in~\cite{KonDaSStrUs98}. It may be readily deduced by combining~\eqref{eq:LapGam} with~\eqref{eq:EGFcycle} to get
\begin{align}\label{eq:LapGPMnom}
\mcL[\GP_\sigma](tf)=\exp \av{\EGF[\Mnom{\mbfe^n}](tf)}_\sigma \comma
\end{align}
where $f=f(x)$ and the expectation is taken with respect to $\diff\sigma(x)$. Comparing the above with~\eqref{eq:EGF3} allows to write the Taylor expansion of $\mcL[\GP_\sigma](tf)$ in $t$ as the exponential generating function in $t$ of Bell polynomials computed at the cumulants $\kappa_n^{\GP_\sigma}[f]$ given in~\eqref{eq:CumGam}.
\end{proof}
\end{prop}

\begin{rem}[Moments of the Poisson measure]\label{r:MomPoi}
The very same computation above also yields the moments of the Poisson measure
\begin{align}\label{eq:MomPoi}
\mu^{\prime\, \PP_\sigma}_n[f]=B_n\quadre{\av{(f\uno)^{\compo \vec\mbfn}}^\compo_{\sigma}} \fstop
\end{align}

Incidentally, notice that, in the case of $\PP_\sigma$, the appearance of Bell polynomials (instead of Touchard polynomials) is peculiar of the infinite-dimensional framework and accounts for the randomness; indeed, when~$\sigma$ degenerates to~$\delta_{x_0}$, the moment $\mu^{\prime\, \PP_\sigma}_n[f]$ reduces to $T_n[f(x_0)]$, for~$\mcL[\PP_\sigma]$ becomes the exponential generating function \eqref{eq:Touchard}, with $f(x_0)$ as parameter.
\end{rem}

Finally, combining~\eqref{eq:MomDF},~\eqref{eq:ProdMeas} and~\eqref{eq:MomToCum} one recovers the cumulants of $\DF_{\beta\n\sigma}$, viz.
\begin{align*}
\kappa^{\DF_{\beta\n\sigma}}_n[f]=B_n^*\quadre{\seq{\Poch{\beta}{1}^{-1} \mu^{\prime\, \GP_\sigma}_1[f],\dotsc,\Poch{\beta}{n}^{-1} \mu^{\prime\, \GP_\sigma}_n[f]}} \comma
\end{align*}
whence those of the Dirichlet distribution by choosing $f$ a simple function of values $\mbfs$ and applying~\eqref{eq:d:DirFer}.

\subsection{Urns \& beads and dynamical symmetry algebras in infinite dimensions}\label{ss:UBDSA}
As for the main properties of the Dirichlet distribution, also the content of~\S\S\ref{ss:Polya}~\&~\ref{ss:Lauricella} may be generalized to infinite dimensions, by a somewhat formal limiting procedure which we sketch in the next paragraphs.

\smallskip

\emph{Throughout this section we assume the space~$(X,\tau)$ to be compact.}

\paragraph{Urns \& beads} We can revisit our metaphor of urns and beads (see~\S\ref{ss:Polya}) for $\DF_\sigma$ as follows. We picture the urns as a family of measurable subsets of~$X$
\begin{align*}
\mcU\eqdef& \set{U_1^1,\dotsc, U_1^{\alpha_1}, \dotsc, U_j^1,\dotsc, U_j^{\alpha_j}, \dotsc, U_k^1,\dotsc, U_k^{\alpha_k} } \comma && \boldalpha\in \N_1^k\comma
\end{align*}
such that $\sigma U_i^j=1$, and arrange them in a measurable partition $\mbfX$ of the measure space~$(X,\sigma)$ with total mass~\mbox{$\abs{\sigma}=\length{\boldalpha}$}, in such a way that
\begin{align}\label{eq:InterprAlpha}
\sigma X_j=&\alpha_j \quad \textrm{ where } \quad X_j\eqdef \bigsqcup_{i\leq\alpha_j} U_j^i \comma j\in[k] \fstop
\end{align}

The weights can then be regarded as the simple function $s\colon X\rar \R$ attaining value $s_j$ on~$X_j$, and, dividing the total mass by $\length{\boldalpha}$, we obtain a probability space which we interpret as in~\S\ref{ss:Polya}.
The very same density argument used to compute the moments~\eqref{eq:MomDF} of~$\DF_\sigma$ allows to take as weights arbitrary (continuous) functions. Indeed, letting
\begin{equation}\label{eq:partition}
\begin{aligned}
\seq{\mbfX^h}_h=\seq{\mbfX^h_{k_h}}_h\eqdef& \seq{X_1^h,\dotsc, X_{k_h}^h}_{h\geq 0}\comma & \mbfX^h\subset& \mfB(\tau) \comma \\
\hlim \max_{j\in [k_h]} \diam X_j^h=&0 \comma & X=&\sqcup_i X_i 
\end{aligned}
\end{equation}
be a \emph{null-array} of measurable partitions of~$X$ (recall that $\diam X^h_j$ vanishes independently of the chosen metric on~$X$, cf.~\cite[\S2.1]{Kal83}), it suffices to set $\boldalpha_{k_h}^h\eqdef \sigma^\compo \mbfX^h_{k_h}$ and let $h$ tend to infinity. Letting~$\mbfX^h_{k_h}$ be such that~$\boldalpha^h_{k_h}=\sigma X^h_1 \uno_{k_h}$ (which is possible since~$\sigma$ is diffusive and~$X$ Polish), we can choose the said decomposition to be independent of $\sigma$.

\paragraph{Infinite-dimensional dynamical symmetry algebras} Within this framework, one can provide a characterization for some of the aforementioned algebraic constructions, as e.g. raising/lowering operators or the whole dynamical symmetry algebra of~${}_k\Phi_2$, in the limiting case when~$h\rar\infty$. However, care must be taken in that the limiting action of such operators, given in~\eqref{eq:R/LOp}, is degenerate, for~$\alpha_i^h$ is vanishing in~$h$. The resulting action makes thus sense only with respect to raising indices `on a measurable set of points', as we sketch in the following.

\smallskip

\emph{At present, Theorem~\ref{t:Borel0} below should be understood as a \emph{conjecture}. We provided only those arguments that allow us to formulate a further conjecture in~\S\ref{s:Conclusions}. We plan to address the details of these limiting constructions in future work.}  

\begin{thm}\label{t:Borel0}
The limiting action of raising operators in the dynamical symmetry algebra of~${}_k\Phi_2$ is meaningful in (Riemann) integral form given on measurable sets~$A$ by
\begin{align}\label{eq:Borel0}
(E_{A} \, \Phi_k)[\sigma,\abs{\sigma},f]\eqdef&\int_A\diff\sigma(x)\, \Phi [\sigma+\delta_x,\abs{\sigma}+1,f]  \fstop
\end{align}

\begin{proof}[Sketch of proof] Let~$\seq{\mbfX^h}_h$ be a sequence of partitions of~$X$ as in~\eqref{eq:partition} and in addition \emph{decreasing} (in the sense that~$\mbfX^{h_2}$ refines~$\mbfX^{h_1}$ for~$h_2\geq h_1$). Since~$\sigma$ is diffusive and~$X$ is Polish, one can additionally choose~$\mbfX^h$ consisting of continuity sets for~$\sigma$ (cf.~\eqref{eq:d:DirFer}), hence, in fact, consisting of closed sets (disjoint up to a $\sigma$-negligible set), compact by compactness of~$X$.
It follows by the finite intersection property that every decreasing sequence of sets~$\seq{X^h_{i_h}}_h$ such that~$X^h_{i_h}\in \mbfX^h$ admits a non-empty limit, which is a singleton because of the vanishing of diameters, hence itself a continuity set for~$\sigma$, since the measure is diffusive.
Vice versa, however chosen a sequence of partitions satisfying the above hypotheses, for every point~$x$ in~$X$ it is not difficult to construct a (possibly non-unique) sequence of closed sets~$X_{i_h}^h$ converging to~$x$ and such that~$X^h_{i_h}\in \mbfX^h$. Furthermore, letting $x$ be a point for which there exists more than one such sequence, we see that for every $h$ the point $x$ belongs to some intersection~$X_{i_{1}}^h\cap X_{i_{2}}^h\cap \dotsc$, hence, since every partition has disjoint interiors by construction,~$x\in \partial X_{i_1}^h\cap \partial X_{i_2}^h\cap \dotsc$. Since for every~$h$ and~$i\leq k_h$ the set~$X_{i}^h$ is a continuity set for~$\sigma$, the whole union~$\cup_{h\geq 0}\cup_{i\leq k_h} \partial X_{i}^h$ is $\sigma$-negligible, thus so is the set of points~$x$ considered above, so that, all in all, for $\sigma$-a.e.~$x$ there exists a unique sequence~$\seq{X^h_{i_h}}_h$ such that~$X^h_{i_h}\in \mbfX^h$ and~$\hlim X_{i_h}^h=\set{x}$.

Let now~$A$ be measurable and let $\seq{\mbfX^h}_h$ be a sequence of partitions as above and in addition refining~$A$. Fix~$f$ in $\mcC(X)$ and such that $\norm{f}_\infty<1$, set~$\boldalpha^h\eqdef \sigma^\compo \mbfX^h$ and notice that~\eqref{eq:R/LOp} yields by summation
\begin{align}\label{eq:Borel}
\tonde{\sum_{i\mid X^h_i\cap A\neq \emp} E_i \, {}_k\Phi_2}\quadre{\boldalpha^h;\length{\boldalpha^h};\mbfs^h}=\sum_{i\mid X^h_i\cap A\neq \emp} \alpha_i^h \, {}_k\Phi_2\quadre{\boldalpha^h+\mbfe^i; \length{\boldalpha^h}+1;\mbfs^h} \comma
\end{align}
where~$\seq{\mbfs^h}_h$ is the sequence of values of simple functions~$s^h$ approximating~$f$ on the partition~$\mbfX^h$. By letting $h\rar\infty$, on the sets~$X^h_{i_h}\rar \set{x}$ we have~$\boldalpha^h+\mbfe^i\rar \sigma+\delta_{x}$. By continuity of~$f$ and of~$\Phi[\sigma, \abs{\sigma},f]$ (i.e. narrow continuity in~$\sigma$, resp. with respect to uniform convergence in~$f$), the Riemann sum in the rhs of~\eqref{eq:Borel} converges to the rhs of~\eqref{eq:Borel0}. 
\end{proof}
\end{thm}

\begin{rem}[A Bayesian nonparametrics interpretation]
The formula~\eqref{eq:Borel0} has a natural interpretation in Bayesian nonparametrics (for the relevance of~$\DF_\sigma$ in Bayesian nonparametrics see e.g.~\cite{Fer73,LijPru10, LijReg04, Set94}). It is proved in~\cite[property \textbf{P3}, Thm.~4.3]{Set94} that the \emph{posterior distribution}~$\DF_\sigma^x$ of a random realization~$\rho$ of~$\DF_\sigma$ given an $X$-valued random variable~$x$ is given by~$\DF_\sigma^x=\DF_{\sigma+\delta_x}$. The formula~\eqref{eq:Borel0} is equivalent to this statement in the conjugate Fourier picture (i.e. at the level of Fourier/Laplace transform); indeed one can rewrite it as
\begin{align*}
E_A \mcL[\DF_\sigma]=\mbbE[ \mcL[\DF_\sigma^x] ; x\in A] \fstop
\end{align*}

In words, if~$x$ is an $A$-valued random variable ($A\subset X$), the Laplace transform of the posterior distribution of~$\rho$ given~$x$ is the Laplace transform of the prior distribution `augmented' by~$E_A$.
\end{rem}

As the above reasoning makes clear, limiting raising/lowering operators are parametrized by points~$x$ in~$X$ and morally act by charging~$x$ with a signed Dirac mass. Since we start with a diffusive probability~$\sigma$, lowering operators would reduce the reference measure to a (non-diffusive) \emph{signed} measure with vanishing mean, hence one would need to allow for such reference measures in the construction of~$\DF_\sigma$. As a consequence, one should in principle focus only on the Borel subalgebra of the dynamical symmetry algebra containing all raising operators.
In the limiting case, this may still be regarded as a subalgebra of an appropriate `limiting dynamical symmetry algebra'. We sketch a construction for the latter in the following theorem.

\begin{conj}\label{t:Borel1}
The limiting dynamical symmetry algebra of~$\Phi[\sigma;1;f]$ may be identified with the $*$-ring~$\msK_0(\msH_\sigma)$ of compact operators on~$\msH_\sigma\eqdef (L^2_\sigma)_\C$ with vanishing trace endowed with the Lie bracket induced by the composition of operators.

\begin{proof}[Motivation]
Standing the definitions and notations established in the proof of Theorem~\ref{t:Borel0}, set for the sake of simplicity~$\boldalpha\eqdef \boldalpha^h$.
Also let~$\mfl_k$ denote the dynamical symmetry algebra of~${}_k\Phi_2[\boldalpha;\length{\boldalpha};\mbfs]$ and notice that, up to conjugation via the Lie isomorphism~$\varphi$ in the proof of Proposition~\ref{p:chased}, the resulting embedding~$\mfl_{k_h} \hookrightarrow \mfl_{k_{h+1}}$ is a natural (though not diagonal) embedding in the sense of~\cite[2.2]{BarZhi99}. Since~$\mfl_{k_h}$ is simple by Theorem~\ref{t:dsa}, we \emph{conjecture} --~as a consequence of the classification in~\cite{BarZhi99}~-- that the Lie algebras direct limit of the system~$\seq{\mfl_{k_h}}_h$ with the said embeddings is a finitary simple Lie algebra. By the classification of such algebras~\cite[1.1]{BarStr02} and interpreting the vector~$\boldalpha=\boldalpha^h$ as in the proof of Theorem~\ref{t:Borel0}, the direct limit satisfies
\begin{align*}
\dirlim \mfl_{k_h}=\mfs\mfl_\infty(\C)\cong \mff\mfs\mfl(\msH_\sigma) \comma
\end{align*}
the finitary special linear algebra of finite rank operators with vanishing trace over the separable Hilbert space~$\msH_\sigma\eqdef (L^2_\sigma)_\C$.

Let now~$\tilde\mfl_k$ denote the vector space underlying to~$\mfl_k$, regarded as a subspace of~$\mfg\mfl_k(\C)$ and endowed with the spectral norm and the Hermitian transposition of matrices, making it a (finite-dimensional) $*$-ring. Notice that the Hermitian transposition is compatible with the Lie bracket up to sign, hence that the embedding $\mfl_{k_h} \hookrightarrow \mfl_{k_{h+1}}$ induces (the restriction of) a $C^*$-embedding~$\tilde\mfl_{k_h}\hookrightarrow \tilde\mfl_{k_{h+1}}$ compatible with the bracket and isometric (recall that any injective $C^*$-morphism is isometric). Essentially because of the above compatibility (cf.~\cite[\S1 and ref.s therein]{BarZhi99}) one has
\begin{align*}
\cl{\dirlim \tilde \mfl}=\cl{\mff\mfs\mfl(\msH_\sigma)}=\msK_0(\msH_\sigma)\comma
\end{align*}
where the closure is taken with respect to the operator norm topology in~$\msB(\msH_\sigma)$ and the first direct limit is here a direct limit of~normed $*$-rings. \phantom\qedhere
\end{proof}
\end{conj}

\section{Fock spaces and chaos representation}\label{ss:Fock}
\emph{For the chaos representation of PT processes we follow~\cite{Sur84} (a modern exposition may be found in~\cite[\S5]{PecTaq11}); for that of the Gamma process cf.~\cite[\S4]{KonDaSStrUs98} and \cite{KonLyt00}; we take the definition of Pascal type processes from~\cite{BarOueRia11}; for another (non-equivalent) definition of Pascal processes and some general results on CPT processes cf.~\cite[\S6 and passim]{Lyt03}.}

\emph{In order to reduce technicalities, we only deal with the case $X=\R^d$ with finite reference measure. The construction of a nuclear Gel'fand triple for $X$ a non-compact smooth orientable manifold with $\sigma$-finite reference measure --~the base step for any generalization of results in this section to the case of manifolds~-- is given in~\cite[\S2]{Lyt03}.
}

\paragraph{Notation} Let $X\eqdef \R^d$. Write~$\hoplus$, resp. $\oplus$, for the Hilbert, resp. algebraic direct sum of (Hilbert) spaces (analogously for tensor products) and set $\Phi^{\otimes 0}=\K$ for any $\K$-linear space~$\Phi$.

Given a finite diffusive reference measure~$\sigma$ on~$X$, denote by~$\msH_\sigma$ the space~$(L^2_\sigma)_{\C}$ of complex-valued square-integrable functions on~$X$, by~$\mcS(\R^d)$ the Schwartz space of real-valued rapidly decreasing functions, and set $\phi,\psi\in\msD\eqdef \mcS(\R^d)_{\C}$; also write~$L^2_{\odot n}\cong \msH_\sigma^{\hodot n}$ for the norm-closure in~$L^2_{\sigma^{\otimes n}}$ of the space~$\msD^{\odot n}$ of symmetric (i.e.~$\phi^{(n)}[\mbfx]=\phi^{(n)}[\mbfx_\pi]$) complex-valued Schwartz functions~$\phi^{(n)}$.
Term further
\begin{align*}
\msF^n(\msH_\sigma)\eqdef& \msH_\sigma^{\hodot n}\comma && \textrm{resp.} & \msF(\msH_\sigma)\eqdef& \overline{\bigoplus_{n\geq 0}}\, \msF^n(\msH_\sigma)\comma
\end{align*}
the \emph{$n$-par\-ticle space} (also $n^\textrm{th}$~\emph{Wiener--It\^o chaos}), resp. the (\emph{bosonic}) \emph{Fock space}, of~$\msH_\sigma$.

\subsection{Chaos representation}\label{ss:Chaos}
Denote now by~$\rho$ any realization of a completely random measure~$\mcR_\sigma$  (or~$\mcR_{\sigma,\lambda}$) on~$X$ and by
\begin{align*}
I_n^\rho \phi^{(n)}\eqdef \int_{X^n} \phi^{(n)}[\mbfx]\diff \rho^{\otimes n}[\mbfx]
\end{align*}
the multiple stochastic integral of~$\phi^{(n)}$ with respect to~$\rho$. We say that~$\mcR_\sigma$ has the \emph{chaos representation property} if the space of chaos decomposable exponential square-summable random variables
\begin{align*}
\msC_{\mcR_\sigma}\eqdef \cl_{L^2_{\mcR_\sigma}} \set{\sum_{n \geq 0} I_n^\rho \phi^{(n)} \mid \phi^{(n)}\in \msD^{\odot n} \comma \sum_{n\geq 0}n!\norm{\phi^{(n)}}_{\msF^n(\msH_\sigma)}^2 < \infty } \subset L^2_{\mcR_\sigma(\diff \rho)}
\end{align*}
coincides with the whole space $L^2_{\mcR_\sigma}$ (further details may be found in e.g.~\cite[\S4.2.1]{Mey93}).

\subsection{Fock spaces of some (C)PT processes}
\paragraph{Fock space of the Poisson point process}
Together with Gaussian processes, Poisson point processes are (virtually, see~\cite[\S4.2.1]{Mey93}) the sole processes possessing the CRP, which appears to have been rediscovered many times (cf. e.g.~\cite[\emph{ibid.}]{Mey93}).
In the case of the Poisson process we show, by reproving a result in~\cite{Sur84}, that multiple stochastic integrals with respect to a centered version of the process may be computed directly from the moments, which were in turn obtained from the cumulants in the previous section.

\begin{thm}[MSI's for PT processes]\label{t:Surgailis}
For any realization~$p$ of~$\PP_\sigma$ define the associated centered (also: \emph{compensated}) Poisson random measure~$q\eqdef p-\sigma$ (here the definition makes sense as it is, since we assumed $\abs{\sigma}<\infty$).
It~holds
\begin{align}\label{eq:MSI0}
I_n^q f^{(n)}=\sum_{k=0}^n \tbinom{n}{k} (-1)^{n-k} f^{(n)}_{\contra k} \comma
\end{align}
where
\begin{align*}
f^{(n)}_{\contra k} \eqdef &\sum_{\mbfx_k\subset \tau(p)\cap \widetilde{X}^n} \int_{X^{n-k}} f^{(n)}[\mbfx_k \oplus {}_k\mbfx] \diff \sigma^{\otimes (n-k)}({}_k \mbfx) & k=&0,\dotsc, n \comma
\end{align*}
having care to interpret the above for $k=n$, resp. $k=0$, as
\begin{align*}
f^{(n)}_{\contra n}=&f^{(n)}[\mbfx] \comma && \textrm{resp.} & f^{(n)}_{\contra 0}=&\int_{X^n} f^{(n)}[\mbfx] \diff\sigma^{\otimes n}(\mbfx) \fstop
\end{align*}

We caution that the well-posedness of this definition for $\sigma$-finite reference measures~$\sigma$ is non-trivial (cf.~\cite[\S4]{Sur84}) and is in fact meaningful only~$\PP_\sigma$-a.e.. Also, our notation differs from the one in~\cite{Sur84}:
The subscript~$\contra k$ denotes again a \emph{contraction}, here in the sense of~\cite[\S6.2]{PecTaq11} or~\cite[\S3.3.1]{Mey93} (where it is denoted by~$\smile^k$).
A proof of the statement is originally in~\cite[4.2]{Sur84}; we provide an independent proof, based on the raw/central moments duality.

\begin{proof}
Set $\mbff\eqdef \seq{f_1,\dotsc,f_n}$ with 
\begin{align}\label{eq:Sur84}
\sigma(\supp f_i\cap \supp f_j)=&0 & i,j\in [n], i\neq j 
\end{align}
and ${\mbff}^\odot\eqdef \tfrac{1}{n!}\sum_{\pi\in \mfS_n}(\mbff_\pi)^{\uno}$. Since the $\C$-linear span of such functions is dense in $L^2_{\odot n}$, we can choose without loss of generality $f^{(n)}=n!\mbff^\odot$ above, and substitute $\widetilde{X}^n$ with $X^n$ in the summation. This yields
\begin{align*}
f^{(n)}_{\contra k}=k!\int_{X^k} \mbff_k^{\odot}[\mbfy_k] \diff p^{\otimes k}[\mbfy_k] \times (n-k)!\int_{X^{n-k}} {{}_k\mbff}^{\odot}[{}_k\mbfy] \diff\sigma^{\otimes n-k}[{}_k\mbfy] \fstop
\end{align*}

Notice that the second integral does not depend on $p$. Withstanding the usual abuse of notation for linear functionals, taking the expectation with respect to $\PP_\sigma$ and recalling that~\mbox{$\av{f}_\sigma=\av{f}_{\PP_\sigma}$} by~\eqref{eq:MomPoi} yields
\begin{align*}
\tfrac{1}{n!}\tbinom{n}{k}\av{f^{(n)}_{\contra k}}_{\PP_\sigma}=&\mu^{\prime\, \PP_\sigma}_k[{\mbff}^\odot_k] \,\av{{}_k\mbff^\odot}_{\PP_\sigma} \fstop
\end{align*}

Again since $\av{I^q_i{f^{(i)}}}_{\PP_\sigma}= i! \av{\mbff_i^\odot}_{\PP_\sigma}$ for $i=k,n-k\in [n]$, applying~\eqref{eq:MomInvMultiA} with~\mbox{$\mbfm=\uno$} we get~\eqref{eq:MSI0} \emph{in $\PP_\sigma$-average}, thus~\eqref{eq:MSI0} for $\PP_\sigma$-a.e.~$p=q+\sigma$ by arbitrariness of~$\mbff$.
\end{proof}
\end{thm}

As readily checked, e.g. by approximation with simple functions, for~$f$ in~$\msF^n(\msH_\sigma)$
\begin{align*}
I^q_n f\in& L^2_{\PP_\sigma} \comma & \av{I^q_n f}_{\PP_\sigma}=&0 \comma\\
\av{(I^q_n f)^2}_{\PP_\sigma}=&n!\norm{f}_{L^2_n}^2 \comma & \av{I^q_n f \,\overline{\displaystyle{I^q_m} g}}_{\PP_\sigma}=&0 \qquad (n\neq m) \comma
\end{align*}
where the overline denotes here complex conjugation.
It follows by the last orthogonality property that the space~$\msF(\msH_\sigma)$ embeds isometrically in~$L^2_{\PP_\sigma}$. It is then a well-known fact (see e.g.~\cite[\S5.9]{PecTaq11} or~\cite[\S3.3]{Mey93}) that the said embedding is in fact surjective, hence that it is possible to realize the space $L^2_{\PP_\sigma}$ as unitarily isomorphic, via MSI, to the standard \emph{Fock space} over $L^2_\sigma$.

\paragraph{Extended Fock space of the Gamma process}
The \emph{extended Fock space} $\msF_{\ext}$ associated to the Gamma process was discovered in~\cite{KonDaSStrUs98} and thoroughly discussed in~\cite{KonLyt00}. We briefly review the construction in light of the established lexicon to provide some new combinatorial insights.

In analogy with~\eqref{eq:AddContract}, define the \emph{contraction}
\begin{equation}\label{eq:Contra}
\begin{aligned}
\phi^{(n)}_{\contra\boldlambda}[\mbfy]\eqdef \phi^{(n)}[(&\underbrace{y_1,\dotsc,y_{\lambda_1}}_{\lambda_1}, \underbrace{y_{\lambda_1+1}, y_{\lambda_1+1}, y_{\lambda_1+2}, y_{\lambda_1+2}, \dotsc, y_{\lambda_1+\lambda_2},y_{\lambda_1+\lambda_2}}_{2\lambda_2}, \dotsc,\\
& \underbrace{\underbrace{y_{\length{\boldlambda_{k-1}}+1},\dotsc,y_{\length{\boldlambda_{k-1}}+1}}_k, \underbrace{y_{\length{\boldlambda_{k-1}}+2},\dotsc, y_{\length{\boldlambda_{k-1}}+2}}_k,\dotsc, \underbrace{y_{\length{\boldlambda}},\dotsc, y_{\length{\boldlambda}}}_k }_{k\lambda_k})],
\end{aligned}
\end{equation}
where~$\boldlambda \vdash n$ and it is understood that all blocks of variables of length~$i$ are omitted whenever~$\lambda_i=0$.
Denote further (cf.~\cite[(5.19)]{Lyt03}) by $\msF_{\boldlambda}(\msH_\sigma)$ the strong closure in~$\msH_\sigma^{\hodot n}$ of the space of functions~$\phi^{(n)}_{\contra\boldlambda}$ varying~$\phi$ in~$\msD^{\odot n}$ and~$\boldlambda\vdash n$, by
\begin{align*}
\msF_{\ext}^n(\msH_\sigma)\eqdef \overline{\bigoplus_{\boldlambda\vdash n}} \Mnom{\boldlambda} \msF_\boldlambda(\msH_\sigma)
\end{align*}
the \emph{extended $n$-particle space} and by
\begin{align*}
\msF_{\ext}(\msH_\sigma)\eqdef \overline{\bigoplus_{n\geq 0}}\, n!\,\msF_{\ext}^{n}(\msH_\sigma)
\end{align*}
the \emph{extended Fock space} of~$\msH_\sigma$. By Proposition~\ref{p:MomGP} we can write
\begin{align}\label{eq:ScalGCoher}
\mcL[\GP_\sigma](st\phi\psi)=&\sum_{n=0}^\infty \frac{s^nt^n}{n!}\scalar{\phi^{\otimes n}}{\psi^{\otimes n}}_{n,\ext} & s,t\in\R, \, \phi,\psi\in \msD
\end{align}
whence we endow each space~$\msF_\boldlambda(\msH_\sigma)$ with the sesquilinear form defined on coherent states by
\begin{align}\label{eq:MomGamFock}
\scalar{\phi^{\otimes n}}{\psi^{\otimes n}}_{n,\ext}=\mu^{\prime\, \GP_\sigma}_n[\phi\psi] \fstop
\end{align}

By~\eqref{eq:MomGP} and polarization, the latter yields the scalar product of two arbitrary $n$-particle states
\begin{equation}\label{eq:ScalG}
\begin{aligned}
\scalar{\phi^{(n)}}{\psi^{(n)}}_{\msF_{\ext}^{(n)}(\msH_\sigma)}=&\sum_{\pi\in \mfS_n} \av{\overline{\phi^{(n)}_{\contra\boldlambda(\pi)}} \, \psi^{(n)}_{\contra\boldlambda(\pi)}}_{\sigma^{\otimes n}}\\
=& \sum_{\boldlambda\vdash n} (\Gamma^\compo[\vec\mbfn])^{\boldlambda} \av{\overline{\phi^{(n)}_{\contra\boldlambda}} \, \psi^{(n)}_{\contra\boldlambda}}_{\sigma^{\otimes n}} \comma
\end{aligned}
\end{equation}
(here $\overline{\phi}$ denotes the complex conjugate of $\phi$).

\begin{prop}[Two recursive identities]\label{p:Recursive}
The sesquilinear form~\eqref{eq:MomGamFock} satisfies
\begin{align*}
\scalar{\phi^{\otimes (n+1)}}{\psi^{\otimes (n+1)}}_{n,\ext}= \sum_{k=0}^n \tbinom{n}{k} \scalar{\phi^{\otimes k}}{\psi^{\otimes k}}_{k,\ext} \scalar{\phi^{n-k+1}}{\psi^{n-k+1}}_{\msH_\sigma} \fstop
\end{align*}

Additionally, for $\phi_i,\psi_i$ such that $\sigma(\supp \phi_i\cap \supp \psi_j)=0$ for $i,j=1,2$, $i\neq j$, one has
\begin{align*}
&\scalar{(\phi_1+\phi_2)^{\otimes n}}{(\psi_1+\psi_2)^{\otimes n}}_{n,\ext}=\\
&\qquad=\sum_{k=0}^n \tbinom{n}{k}\scalar{(\phi_1\psi_1)^{\otimes k}}{(\phi_1\psi_1)^{\otimes k}}_{k,\ext} \scalar{(\phi_2\psi_2)^{\otimes (n-k)}}{(\phi_2\psi_2)^{\otimes (n-k)}}_{n-k,\ext} \fstop
\end{align*}

\begin{proof}
Since the scalar product of coherent $n$-particle states is a moment (by~\eqref{eq:MomGamFock}), it may be expressed as a Bell polynomial computed at suitable cumulants (by~\eqref{eq:CumToMom}), hence satisfies the recursive identity~\eqref{eq:RecBell}.
Notice that this constitutes in fact a one-line proof for~\cite[Thm.~2.3]{BarOueRia11}, originally stated for the scalar product of the \emph{Pascal one-mode type interacting Fock space} (see below).

\smallskip

Concerning the second identity, the assumptions on $\phi_i,\psi_i$ yield
\begin{align*}
(\phi_1+\phi_2)(\psi_1+\psi_2)=&\phi_1\psi_1+\phi_2\psi_2 \comma\\
\av{(\phi_1+\phi_2)^i(\psi_1+\psi_2)^i}_\sigma=\av{(\phi_1\psi_1+\phi_2\psi_2)^i}_\sigma=&\av{(\phi_1\psi_1)^i}_\sigma+\av{(\phi_2\psi_2)^i}_\sigma && i\in [n] \fstop
\end{align*}

Thus, similarly as before, the scalar product satisfies the binomial type identity~\eqref{eq:BinomType}, viz.
\begin{align*}
&\scalar{(\phi_1+\phi_2)^{\otimes n}}{(\psi_1+\psi_2)^{\otimes n}}_{n,\ext}=\mu_n^{\prime \, \GP_\sigma}[(\phi_1+\phi_2)(\psi_1+\psi_2)]=\mu_n^{\prime \, \GP_\sigma}[\phi_1\psi_2+\phi_2\psi_2]\\
&=\sum_{k=0}^n \tbinom{n}{k}\scalar{(\phi_1\psi_1)^{\otimes k}}{(\phi_1\psi_1)^{\otimes k}}_{k,\ext} \scalar{(\phi_2\psi_2)^{\otimes (n-k)}}{(\phi_2\psi_2)^{\otimes (n-k)}}_{n-k,\ext} \fstop \qedhere
\end{align*}
\end{proof}
\end{prop}

\begin{rem}[A combinatorial interpretation of the \emph{extended} Fock space]\label{r:ExtFock}
The ter\-mi\-nol\-o\-gy of \emph{extended Fock space} is motivated by the fact, firstly noted in~\cite{KonLyt00}, that the canonical Fock space~$\msF(\msH_\sigma)$ associated to the Poisson point process embeds in~$\msF_\ext(\msH_\sigma)$. Indeed, setting~-- cf.~\eqref{eq:Sur84}~--
\begin{align*}
\msR^n\eqdef \Span_\C\set{\phi^{(n)}\in \msD^{\odot n}\mid \sigma(\supp \phi_i\cap \supp\phi_j)=0, \, i,j\in [n], i\neq j}
\end{align*}
yields $\cl_{\norm{\emparg}_{\msF^n(\msH_\sigma)}}\msR^n=\msF^n(\msH_\sigma)$. Moreover, by disjointness of supports one has
\begin{align*}
\scalar{\phi^{(n)}}{\psi^{(n)}}_{\msF^n(\msH_\sigma)}=&\scalar{\phi^{(n)}}{\psi^{(n)}}_{\msF_{\ext}^{n}(\msH_\sigma)} & \phi^{(n)},\psi^{(n)}\in& \msR^n\comma
\end{align*}
for the only non-vanishing term in the sum~\eqref{eq:ScalG} is the one given by $\boldlambda=\seq{{}_1n, 0, \dotsc, {}_n 0}$, for which $(\Gamma^\compo[\vec\mbfn])^\boldlambda=1$. Hence one can canonically identify $\oplus_{n\geq 0} \msR^n$ with a linear subspace of~$\msF_\ext(\msH_\sigma)$ and thus its closure~$\msF(\msH_\sigma)$ with a subspace of $\msF_\ext(\msH_\sigma)$.

Incidentally, let us notice that the above integer partition~$\boldlambda=n\mbfe_1$ uniquely corresponds to the set partition~$\mbfL=\set{\set{1},\dotsc, \set{n}}$, which is the infimum of the lattice of set partitions of~$[n]$.

As a consequence of~\eqref{eq:LapGPMnom}, each of the coefficients~$\Gamma[\ell]$ of the factor~$\Gamma^\compo[\vec\mbfn]$ in~\eqref{eq:ScalG} ought to be thought of as the multinomial number~$\Mnom{\mbfe^\ell}$, counting the number of cyclic permutations with length $\ell$. Indeed, computing the cycle index polynomial in~\eqref{eq:ScalGCoher} (cf.~\eqref{eq:MomGP}) amounts to compose all of the possible disjoint cyclic permutations together in such a way that the resulting permutation has order~$n$.
This was already partially understood in~\cite[\S1]{KonLyt00} and~\cite[Rem.~6.1]{Lyt03}, where the orbit of a cyclic permutation is termed a \emph{loop} and the summation in~\eqref{eq:ScalG} (see also~\cite[(4.6)]{KonDaSStrUs98}) is indexed over the set of all possible collections of non-intersecting loops, i.e. disjoint orbits of cyclic permutations acting over~$[n]$. 
\end{rem}

\paragraph{Other `non-standard' Fock spaces}
The contraction~\eqref{eq:Contra} may be regarded as a `folding' of~$X^{\times n}$, rather than of functions on~$X$. The main idea for such `foldings' is already in~\cite[Prop.~3]{RotWal97}, whereas the full extent of the construction is presented in~\cite[\S6]{Lyt03}, mostly by means of set partitions. We comparatively review this construction using both set and integer partitions.
Given a subset~$A\subset X^{\times n}$, let $A_\pi\eqdef \set{\mbfx_\pi\mid \mbfx\in A}$ and term~$A$ symmetric if~$A=A_\pi$ for~$\pi$ in~$\mfS_n$. Also set
\begin{align*}
X^{\times \mbfL}\eqdef \{\mbfx&\in X^{\times n}\mid \ell \in L_i \implies x_\ell =y_i, y_i\neq y_j, i,j\in [n], i \neq j\} & \mbfL\vdash& [n] \comma
\intertext{resp.}
X^{\times \boldlambda}\eqdef \{\mbfx&\in X^{\times n}\mid \mbfx=\mbfy\textrm{ as in~\eqref{eq:Contra}, } y_i\neq y_j, i,j\in[n], i\neq j\} & \boldlambda \vdash& n \fstop
\end{align*}

In analogy with~\cite[p.~1268]{RotWal97}, we say that such sets are \emph{$\boldlambda$-weakly triangular} (setting~$\boldlambda\eqdef \boldlambda(\mbfL)$ in the second case); in the notation of Remark~\ref{r:CQES} one has $X^{\times \boldlambda}=X^{\times\vec\mbfL(\boldlambda)}$.
It holds either that~$X^{\times \boldlambda}_\pi=X^{\times \boldlambda}_{\pi'}$ or $X^{\times \boldlambda}_{\pi}\cap X^{\times \boldlambda}_{\pi'}=\emp$ and for~$\mbfL$ and~$\boldlambda$ such that~$\boldlambda(\mbfL)=\boldlambda$ there exists a unique $\pi=\pi(\mbfL)$ in $\mfS_n(\boldlambda)\subset \mfS_n$ such that~$X^{\times\boldlambda}_\pi=X^{\times \mbfL}$. Thus
\begin{align*}
X^{\times n}=\bigsqcup_{\mbfL\vdash [n]} & X^{\times \mbfL}=\bigsqcup_{\boldlambda\vdash n} \,\,\bigsqcup_{\pi\in \mfS_n(\boldlambda)} X^{\times \boldlambda}_\pi \comma & \textrm{resp.} && X^{\odot n}=& \bigsqcup_{\boldlambda\vdash n} X^{\odot\boldlambda} \fstop
\intertext{where~$X^{\odot\boldlambda}\eqdef X^{\boldlambda}/\mfS_n$.
Given now~$\vartheta=\vartheta_n$ a symmetric measure on~$X^{\times n}$ (i.e. such that~$\vartheta A=\vartheta A_\pi$ for~$\pi\in\mfS_n$), let~$\vartheta^{{}^\odot}$ denote its pushforward on~$X^{\odot n}$ under the quotient map~$\pr_n$ by the action of~$\mfS_n$ and set~$\vartheta_{\mbfL}\eqdef \vartheta\restr_{X^{\times \mbfL}}$, resp.~$\vartheta^{{}^\odot}_\boldlambda\eqdef \vartheta^{{}^\odot}\restr_{X^{\odot \boldlambda}}$, whence, in analogy with the polynomials~$B_n$, resp.~$Z_n$,
}
\av{\phi^{\otimes n}}_{\vartheta}=&\sum_{\boldlambda\vdash n} \Bella{\boldlambda} \av{\phi^{\otimes n}}_{\vartheta_{\vec\mbfL(\boldlambda)}} \comma & \textrm{resp.} && \av{\phi^{\odot n}}_{\vartheta^{{}^\odot}}=& \sum_{\boldlambda \vdash n} \Mnom{\boldlambda} \av{\phi^{\odot n}}_{\vartheta^{{}^{\odot}}_\boldlambda} \fstop
\end{align*}
 
Following~\cite{Lyt03}, consider the measurable isomorphism
\functionnn{T_\boldlambda}{\widetilde X^{\times \length{\boldlambda}}}{X^{\times \boldlambda}}{\seq{y_1,\dotsc, y_{\length{\boldlambda}}}}{\mbfy \textrm{ as in~\eqref{eq:Contra}}}
and set~$T^{{}^{\odot}}_\boldlambda\eqdef \pr_n\circ T_\boldlambda$, and $\tilde\vartheta\eqdef (T_\boldlambda^{-1})_\pfwd \vartheta$, resp.~$\tilde\vartheta^{{}^\odot}\eqdef (T_\boldlambda^{{}^\odot-1})_\pfwd \vartheta^{{}^\odot}$. Directly from the above definitions one has that
\begin{align*}
\av{\displaystyle{\phi^{(n)}_{\contra\boldlambda}}}_{\tilde \vartheta_\boldlambda }=&\av{\phi^{(n)}}_{\vartheta_\boldlambda}\comma & \textrm{resp.} && \av{\displaystyle{\phi^{(n)}_{\contra\boldlambda}}}_{\tilde \vartheta^{{}^\odot}_\boldlambda }=&\av{\phi^{(n)}}_{\vartheta^{{}^\odot}_\boldlambda} \fstop
\end{align*}

When~$\vartheta=\rho^{\otimes n}$ with~$\rho$ a realization of some L\'{e}vy process with law~$\mcR$, one has~$\mcR(\diff\rho)$-a.e.
\begin{align}\label{eq:Lytvynov}
I_n^\rho\phi^{(n)}= \sum_{\boldlambda\vdash n}\Bella{\boldlambda} \av{ \phi^{(n)}_{\contra\boldlambda}}_{(\widetilde{\rho^{\otimes n}})_\boldlambda} \comma
\end{align}
where $I_n$ denotes MSI with respect to~$\mcR$.
As shown in~\cite[\S6]{Lyt03}, the study of Jacobi fields for general CPT laws and especially laws of Meixner type dwells on this and similar representations;
it is thus perhaps not surprising that the sesquilinear form~\eqref{eq:ScalGCoher}, defining the scalar product of extended $n$-particle spaces for~$\GP_\sigma$, appears as a constituent brick of non-standard Fock spaces associated to other processes.
As an example, recall the definition~\cite{BarOueRia11} of \emph{Pascal white noise measure}~$\Lambda_{r,\alpha}$ as the unique (analytically of exponential type) measure satisfying
\begin{align*}
\mcK[\Lambda_{r,\alpha}](\phi)=&\av{\ln\frac{\alpha}{1-\chi \exp[\phi]}}_{\nu_{r,\alpha}}&& \phi\in \msD\comma \norm{\phi}_\infty<-\ln \chi\comma
\end{align*}
where $\nu_{r,\alpha}$ is the \emph{negative binomial distribution} on~$\R$
\begin{align*}
\nu_{r,\alpha}\eqdef& \sum_{m\geq 0} \alpha^r \tbinom{-r}{m}(-\chi)^m \delta_m && r>0\comma\alpha\in (0,1) \comma \chi\eqdef 1-\alpha 
\end{align*}
(the sufficiency of the specification of the cumulant generating function follows from~\S\ref{ss:Radon}-\ref{ss:MomCum}).

It was proved in~\cite[Thm.~3.2]{BarOueRia11} (see also~\cite[p.~93]{Lyt03}) that
\begin{align*}
L^2_{\Lambda_{r,\alpha}}\cong \msF_{\mathsc{nb}}(\msH_{\nu_{r,\alpha}})\eqdef \overline{\bigoplus_{n\geq 0}}\, n!\tonde{\tfrac{\chi}{\alpha^2}}^n\msF_{\mathsc{nb}}^n(\msH_{\nu_{r,\alpha}}) \comma
\end{align*}
the \emph{Pascal one-mode type interacting Fock space} (here: $\mathsc{nb}$ stays for \emph{negative binomial},~cf.~\cite[2.2]{BarOueRia11}), where $\msF_{\mathsc{nb}}^n(\msH_\sigma)$ denotes the space $\msH_\sigma^{\hodot n}$ endowed with the scalar product induced by~$\scalar{\emparg}{\emparg}_{n,\ext}$.

\section{Conclusions and further developments}\label{s:Conclusions}
Our analysis strongly hints to cumulants as an interesting object in the framework of laws of random measures; despite some technicalities due to the setting of choice, the essentially combinatorial nature of the properties of cumulants, and in particular of dualities, grants suitable extensions of these properties to hold in virtually any infinite-dimensional setting. Whereas the said properties are useful mainly for computational purposes, their interplay with a possible combinatorial structure intrinsic of the measures in question may be used to deepen the study of the measures themselves. To this extent, the case of the Gamma measure and even more of the Dirichlet--Ferguson measure is exemplar, for the appearance of cycle index polynomials (consequence of the duality and of the moments of the Gamma distribution) provides new proofs for properties of the extended Fock space of the Gamma measure and unforeseen connections with P\'olya enumeration theory and the theory of dynamical symmetry algebras.

We undertook a detailed study of these connections to approach the quasi-invariance of the Dirichlet--Ferguson measure with respect to the action induced by self-transformations of the underlying space, which we motivate in the following.

\paragraph{A rigidity approach to the quasi-invariance of~$\DF_\sigma$} Starting with Kolmogorov Extension Theorem, finite-dimensional approximation is a common tool in the study of random measures and more generally of cylindrical measures on infinite-dimensional spaces; it is satisfactory to an almost surprising extent in the case when the approximating sequence of choice possesses some \emph{intrinsically linear structure}, making the direct limit procedure at its core particularly well-behaved.
Besides L\'{evy}'s construction of Brownian motion, two beautiful examples of the r\^ole played by linear spaces approximation are Feyel--\"{U}st\"{u}nel~\cite{FeyUst04} construction of optimal transport maps on the Wiener space as limits of solutions to the Monge--Kantorovich problem on finite-dimensional Hilbert spaces and Vershik's construction of the so-called infinite-dimensional Lebesgue measure~\cite{Ver07, VerGra09} as a limit of Haar measures on maximal toral subgroups of special linear groups.

In more general cases though, where no such linear structure is at hand, one usually strives for further rigidity, partly in order to allow for more general --~yet rigorous~-- limiting procedures. The said rigidity may be gained in several ways. Three meaningful examples are Kondratiev--Lytvynov--Vershik~\cite{KonLytVer15} \emph{chaos representation} of~$L^2_{\GP_\sigma}$ and realization of Laplace operators on it as second quantizations of operators on~$L^2_\sigma$; von Renesse--Sturm~\cite{vReStu09} construction of Wasserstein diffusion on one-dimensional spaces, basing on order properties of the line; Chodosh~\cite{Cho12} negative statement on Ricci curvature bounds for the entropic measure, via \emph{embedding} probability measures on the interval into the space of square-integrable functions over it.

Here, we propose a rigidity approach, via embedding, to the quasi-invariance of~$\DF_\sigma$ under the action via push-forward of the group of diffeomorphisms~$\Diffeo(X)$ when~$X$ is a smooth compact orientable Riemannian manifold, that is, the existence of the Radon-Nikod\'ym derivative satisfying $R_\psi [\n\eta]\diff \DF_\sigma(\n\eta)= \diff (\psi.)_\pfwd \DF_\sigma(\n\eta)$, where $\psi\in \Diffeo(X)$ and~$\psi.\n\eta\eqdef \psi_\pfwd\n\eta$ for $\n\eta$ in $\msP$. The connected problem of studying the quasi-invariance of~$\GP_\sigma$ under the same action has been negatively answered in~\cite[and ref.s therein]{KonLytVer15}, leading to the weaker notion of partial quasi-invariance, while the quasi-invariance of the entropic measure on the interval (also related to that of~$\DF_\sigma$) was given an affirmative answer in~\cite{vReStu09}.

We conjecture the existence of a group~$G$ of self-transformations of~$X$, such that 
\begin{align*}
\Isom(X) \lneq G<\Homeo(X)
\end{align*}
and satisfying the following. The group $G$, regarded as acting on~$L^2_\sigma$ by precomposition, whence on $\msK_0(\msH_\sigma)$ by conjugation, is isomorphic to a subgroup~$H$ of the Lie $*$-automorphisms of~$\msK_0(\msH_\sigma)$, itself acting by conjugation. Furthermore, the action of~$H$ is the limiting action --~in a suitable topology~-- of those of symmetric groups detailed in Proposition~\ref{p:chased}, corresponding, up to the aforementioned isomorphism~$\varphi$, to the symmetric group actions on standard simplexes described in Remark~\ref{r:Models}.

\section{Appendix: Dynamical symmetry algebras of Lauricella functions}\label{app:A}

\emph{The theory of dynamical symmetry algebras for Lauricella hypergeometric functions was developed in the series of papers~\cite{Mil72, Mil73, Mil73b}. Heuristics for such construction starting from \emph{quantum theory} and a motivation for the introduction of dummy variables are found in~\cite[\S1]{Mil73b}. Characterizations via systems of partial differential equations are found in~\cite[\S5]{Mil77}.
We refer the reader to~\cite{Bou81} and~\cite{Hum92} for the general theory of Lie algebras and Weyl groups respectively.}

\smallskip

We recall some of the results in~\cite{Mil72} and proceed to compute the dynamical symmetry algebra of the Humbert function~${}_k\Phi_2$ defined in~\eqref{eq:LauIntRep}.

\paragraph{The dynamical symmetry algebra of \texorpdfstring{${}_kF_D$}{kFD}}
Let~$f_{a,\mbfb,c}\colon \C^{2k+2}_{\mbfx,\mbfu, s,t}\longrightarrow \C$ be defined by
\begin{align}\label{eq:basis}
f_{a,\mbfb,c}\eqdef f_{a,\mbfb,c}(\mbfx, \mbfu, s,t)=& \Beta[a,c-a]\, {}_kF_D[a;\mbfb;c;\mbfx] s^a \mbfu^\mbfb t^c \comma & c\not\in \Z^-_0
\end{align}
where the normalization factor $\Beta[a,c-a]$ is for computational purposes (as it cancels the one in~\eqref{eq:LauIntRep}) and may in fact be omitted~(cf.~\cite[p.226]{Mil73}). On the space of holomorphic functions~$\mcO(\C^{2k+2}_{\mbfx,\mbfu, s,t})$ define the following differential operators
\begin{equation}\label{eq:DiffOpF}
\begin{aligned}
E_a\eqdef& s\tonde{\mbfx\cdot\nabla^\mbfx+s\partial_s} \comma\\
&E_{-a}\eqdef s^{-1}\tonde{(\mbfx\compo(\uno-\mbfx))\cdot \nabla^{\mbfx}+t\partial_t-s\partial_s-\mbfx\cdot (\mbfu\compo\nabla^\mbfu)}\comma\\
E_{b_i}\eqdef& u_i\tonde{x_i\partial_{x_i}+u_i\partial_{u_i}} \comma\\
&E_{-b_i}\eqdef u_i^{-1}\tonde{x_i(\uno-\mbfx)\cdot\nabla^\mbfx+t\partial_t -x_i s \partial_s-\mbfu\cdot\nabla^\mbfu} \comma\\
E_c\eqdef& t\tonde{(\uno-\mbfx)\cdot\nabla^{\mbfx} +t\partial_t-s\partial_s-\mbfu\cdot\nabla^{\mbfu}} \comma E_{-c}\eqdef t^{-1}\tonde{\mbfx\cdot\nabla^\mbfx+t\partial_t-1} \comma\\
E_{a,c}\eqdef& st\tonde{(\uno-\mbfx)\cdot\nabla^{\mbfx}-s\partial_{s}} \comma \\
& E_{-a,-c}\eqdef (st)^{-1}\tonde{(\mbfx\compo(\uno-\mbfx))\cdot \nabla^\mbfx-\mbfx\cdot(\mbfu\compo\nabla^\mbfu)+t\partial_t-1} \comma\\
E_{b_i,c}\eqdef& u_it\tonde{(x_i-1)\partial_{x_i}+u_i\partial_{u_i}}\comma\\
&E_{-b_i,-c}\eqdef (u_i t)^{-1}\tonde{(x_i-1)\mbfx\cdot\nabla^\mbfx+x_is\partial_s-t\partial_t+1} \comma\\
E_{a,b_i,c}\eqdef& su_it\partial_{x_i} \comma\\
E_{-a,-b_i,-c}&\eqdef (s u_i t)^{-1}\tonde{(\mbfx\compo(\uno-\mbfx))\cdot\nabla^\mbfx-t\partial_t+x_is\partial_s+(\mbfx-\mbfe^i)\cdot(\mbfu\compo\nabla^\mbfu)-x_i+1} \comma\\
E_{b_i,-b_j}\eqdef&u_iu_j^{-1}\tonde{(u_i-u_j)\partial_{x_i}+u_i\partial_{u_i}}\comma\\
J_a\eqdef& 2s\partial_s-t\partial_t \comma J_{b_i}\eqdef 2 u_i\partial_{u_i}-t\partial_t+\mbfu\cdot\nabla^\mbfu \comma\\
J_c\eqdef& 2t\partial_t-s\partial_s-\mbfu\cdot\nabla^\mbfu-1\comma
\end{aligned}
\end{equation}
where $i,j\in[k]$, $j\neq i$ and $\nabla^\mbfy\eqdef \seq{\partial_{y_1},\dotsc,\partial_{y_k}}$ for $\mbfy=\mbfu,\mbfx$. We choose slightly different operators with respect to to those in~\cite{Mil72} in order to obtain more standard commutation relations for the~$\mfs\mfl_2$-triples below. The operators do not depend on the parameters $a,\mbfb,c$: the subscripts indicate which indices they raise, resp. lower, viz.
\begin{equation}\label{eq:OpAct}
\begin{aligned}
E_a f_{a,\mbfb,c}=&(c-a-1)f_{a+1,\mbfb,c} \comma & E_{-a}f_{a,\mbfb,c}=&(a-1)f_{a-1,\mbfb,c} \comma\\
E_{b_i} f_{a,\mbfb,c}=& b_i f_{a,\mbfb+\mbfe^i,c} \comma & E_{-b_i} f_{a,\mbfb,c}=&(c-\length{\mbfb}) f_{a,\mbfb-\mbfe^i,c} \comma\\
E_c f_{a,\mbfb,c}=&(c-\length{\mbfb}) f_{a,\mbfb,c+1} \comma & E_{-c}f_{a,\mbfb,c}=&(c-a-1)f_{a,\mbfb,c-1} \comma\\
E_{a,c} f_{a,\mbfb,c}=& (\length{\mbfb}-c) f_{a+1,\mbfb,c+1} \comma & E_{-a,-c}f_{a,\mbfb,c}=& (a-1) f_{a-1,\mbfb,c-1} \comma\\
E_{b_i,c} f_{a,\mbfb,c}=& b_i f_{a,\mbfb+\mbfe^i,c+1} \comma & E_{-b_i,-c} f_{a,\mbfb,c}=& (a-c+1) f_{a,\mbfb-\mbfe^i,c-1} \comma\\
E_{a,b_i,c} f_{a,\mbfb,c}=& b_i f_{a+1,\mbfb+\mbfe^i,c+1} \comma & E_{-a,-b_i,-c} f_{a,\mbfb,c}=& (1-a)f_{a-1,\mbfb-\mbfe^i,c-1} \comma\\
E_{b_i,-b_j} f_{a,\mbfb,c}=& b_i f_{a,\mbfb+\mbfe^i-\mbfe^j,c} \comma & J_a f_{a,\mbfb,c}=& (2a-c) f_{a,\mbfb,c}\comma\\
J_{b_i} f_{a,\mbfb,c}=&\tonde{b_i+\length{\mbfb}-c}f_{a,\mbfb,c} \comma & J_c f_{a,\mbfb,c}=& \tonde{2c-a-\length{\mbfb}-1} f_{a,\mbfb,c} \comma
\end{aligned}
\end{equation}
where $i,j\in[k]$, $j\neq i$. Analogous computations are exemplified in Lemma~\ref{l:Calcoli} below.

It was stated in~\cite{Mil72} that the complex linear span of the operators~$J_{\pm,0}$ below, endowed with the bracket induced by their composition, is an \emph{irreducible} representation of a \emph{simple} Lie algebra~${}_k\mfg$ of dimension~$(k+3)^2-1$. By simplicity and dimension, ${}_k\mfg\cong\mfs\mfl_{k+3}(\C)$, the Lie algebra of $(k+3)$-square complex matrices with vanishing trace. Moreover, ${}_k\mfg$ has $\mfs\mfl_2$-triples
\begin{equation}\label{eq:triples}
\begin{aligned}
\{J_+,&J_-,J_0\}\eqdef\\
&\set{E_a,E_{-a},J_a}, \set{E_c, E_{-c}, J_c}, \set{E_{a,c},E_{-a,-c}, J_a+J_c}, \\ 
&\set{E_{b_i}, E_{-b_i}, J_{b_i}}, \set{E_{b_i,c},E_{-b_i,-c}, J_{b_i}+J_c} && i\in[k]\\
&\set{E_{a,b_i,c}, E_{-a,-b_i,-c}, J_a+J_{b_i}+J_c} && i\in[k]\\
&\set{E_{b_i-b_j}, E_{-b_i,b_j}, J_{b_i}-J_{b_j}} &&i,j\in [k], i\neq j \comma
\end{aligned}
\end{equation}
each satisfying
\begin{align*}
[J_+,J_-]=J_0\comma [J_0,J_{\pm}]=\pm 2J_\pm \comma
\end{align*}
thus isomorphic to $\mfs\mfl_2(\C)$ and corresponding to a root in the root system~${}_k\Psi$ of~${}_k\mfg$.
Detailing~\cite{Mil72}, one can check by tedious yet straightforward computations that the complex linear span of the vectors~$J_0$ is a \emph{maximal toral} sub-algebra~\mbox{${}_k\mfh<{}_k\mfg$}, that a set of positive roots induced by~${}_k\mfh^*$ and labeled by the above triples is given by
\begin{align*}
{}_k\Psi^+\eqdef& \set{\alpha_a, \alpha_c, \alpha_{b_i}, \alpha_{a,c}, \alpha_{b_i,c}, \alpha_{a,b_i,c}, \alpha_{b_i,-b_j}} & i,j&\in[k], j<i \comma
\end{align*}
and that a corresponding Dynkin diagram for ${}_k\mfg$ is as follows
\begin{align}\label{eq:Dynkin}
  \begin{tikzpicture}[scale=.4, baseline=(current bounding box.center)]
    \foreach \x in {0,...,5}
    \draw[xshift=3*\x cm,thick] (\x cm,0) circle (.2cm);
    \draw (0,-1) node[anchor=center] {$\alpha_a$};
    \draw (4,-1) node[anchor=center] {$\alpha_c$};
    \draw (8,-1) node[anchor=center] {$\alpha_{b_1}$};
    \draw (12,-1) node[anchor=center] {$\alpha_{b_2,-b_1}$};
    \draw (16,-1) node[anchor=center] {$\alpha_{b_3,-b_2}$};
    \draw (20,-1) node[anchor=center] {\quad$\alpha_{b_k,-b_{k-1}}$};
    %
    %
    \foreach \y in {0,...,3}
    \draw[xshift=3.875*\y cm,thick] (3.875*\y+0.2 cm,0) -- +(3.6 cm,0);
   \foreach \y in {4}
    \draw[xshift=3.875*\y cm,thick,dotted] (3.875*\y+0.2 cm,0) -- +(3.6 cm,0);
    \draw (25,0) node[anchor=east]{$\comma$};
  \end{tikzpicture}
\end{align}
where labels denote simple roots.

\paragraph{The dynamical symmetry algebra of~${}_k\Phi_2$}
Dynamical symmetry algebras of confluent forms of (hypergeometric) functions may be obtained by a formal \emph{contraction} procedure. In particular, as stated in~\cite[\S2]{Mil72}, dynamical symmetry algebras of confluent forms of~${}_kF_D$ are sub-algebras of
\begin{align*}
\mfg\mfl_{k+3}(\C)\cong \mfe_1(\C) \oplus \mfs\mfl_{k+3}(\C) 
\end{align*}
(where $\mfe_d$ denotes the (unique) abelian Lie algebra of dimension~$d$), resulting from dropping redundantly acting basis vectors and formally contracting along the basis vectors whose action depends on the parameter(s) with respect to which the confluent form of~${}_kF_D$ is taken.

For the case of~${}_k\Phi_2$ we proceed as follows. Firstly, we can immediately drop all basis vectors raising/lowering the parameter~$a$. Secondly, notice that we are in fact only interested in the case when~\mbox{$c=\length{\mbfb}$} (cf.~\eqref{eq:LapDir}), thus we can also drop the basis vectors raising/lowering~$c$ without affecting~$\mbfb$ and those only raising/lowering a \emph{single} element~$b_i$ in~$\mbfb$ without affecting~$c$. Consistently with this choice we see that, fixing~$c=\length{\mbfb}$, the action on~$\mcO(\C^{2k+2}_{\mbfx,\mbfu,s,t})$ of the dropped operators~$E_c$ and~$E_{-b_i}$ identically vanishes. Finally, we can drop $J_c$, for its action linearly depends on those of~$J_a$ and $J_{b_i}$. We are thus left with the basis
\begin{align*}
J_a\comma J_{b_i}\comma E_{b_i,-b_j}\comma E_{b_i,c}\comma E_{-b_i,-c} && i,j\in [k], i\neq j \comma
\end{align*}
whose complex linear span, of dimension~$k^2+2k+2$, we denote by~${}_k\mfg_0$. Setting
\begin{align*}
J_{a}'\eqdef a^{-1}J_a\comma E_{-b_i,-c}'\eqdef a^{-1} E_{-b_i,-c}
\end{align*}
and formally letting $a\rar \infty$, we obtain the dynamical symmetry algebra of ${}_k\Phi_2[\mbfb;\length{\mbfb}]$ as a new Lie algebra structure on ${}_k\mfg_0$, in general not isomorphic to a subalgebra of ${}_k\mfg$ (cf.~\cite[p.~1398]{Mil72}).

In order to identify~${}_k\mfg_0$ we first provide some heuristics. Fixing the parameter~$c$ has the effect of disconnecting the Dynkin diagram~\eqref{eq:Dynkin} (cf.~\cite[p.~228]{Mil73} for the case of Gauss hypergeometric function ${}_2F_1$). The contraction has then the effect of reducing the $\mfs\mfl_2$-triple corresponding to $a$ to the single basis vector~$J_a'$ and of making the span of the remaining vectors closed with respect to the new Lie bracket. Since~$J_a'$ acts as the identity on~$\mcO(\C^{2k+2}_{\mbfx,\mbfu,s,t})$, it belongs to the center~$\mfz({}_k\mfg_0)$ of~${}_k\mfg_0$. This yields
\begin{align*}
\begin{tikzpicture}[scale=.4, baseline=(current bounding box.center)]
    \draw[fill](0,0) circle (.05cm);
    \foreach \x in {1,...,5}
    \draw[xshift=3*\x cm,thick] (\x cm,0) circle (.2cm);
    \draw (0,-1) node[anchor=center] {$J_a'$};
    \draw (4,-1) node[anchor=center] {$\alpha_c$}; \draw (4,0) node[cross] {};
    \foreach \y in {2,...,3}
    \draw[xshift=3.875*\y cm,thick] (3.875*\y+0.2 cm,0) -- +(3.6 cm,0);
   \foreach \y in {4}
    \draw[xshift=3.875*\y cm,thick,dotted] (3.875*\y+0.2 cm,0) -- +(3.6 cm,0);
    \draw (25,0) node[anchor=east]{$\comma$};
  \end{tikzpicture}
\end{align*}
that is~${}_k\mfg_0\cong \mfe_1(\C)\oplus \mfs\mfl_{k+1}(\C)=\mfg\mfl_{k+1}(\C)$.

A rigorous proof of the statement requires starting from scratch.
Let~$f_{\mbfb,c}\colon \C^{2k+1}_{\mbfx,\mbfu, t}\longrightarrow \C$ be defined by
\begin{align}\label{eq:basis2}
f_{\mbfb,c}\eqdef f_{\mbfb,c}(\mbfx, \mbfu, t)=& {}_k\Phi_2[\mbfb;c;\mbfx] \mbfu^\mbfb t^c \comma & c\not\in \Z^-_0 \semicolon
\end{align}
on the space of holomorphic functions~$\mcO(\C^{2k+1}_{\mbfx,\mbfu, t})$ define, with the same notation of~\eqref{eq:DiffOpF}, the following differential operators
\begin{equation}\label{eq:DiffOpPhi}
\begin{aligned}
\tilde E_{b_i,-b_j}\eqdef& E_{b_i,-b_j} \comma \tilde J_{b_i,c}\eqdef t\partial_t+u_i\partial_{u_i}-1\comma\\
\tilde E_{b_i,c}\eqdef& u_it(x_i\partial_{x_i}+u_i\partial_{u_i}-(\mbfx\cdot \nabla^\mbfx)\partial_{x_i})\comma \\
&\tilde E_{-b_i,-c}\eqdef (u_i t)^{-1}(x_i-\mbfx\cdot\nabla^{\mbfx}-t\partial_t+1) && i,j\in[k]\comma i\neq j \fstop
\end{aligned}
\end{equation}

In analogy with~\eqref{eq:OpAct}, one has the following.

\begin{lem}[Raising/lowering actions]\label{l:Calcoli}
The following identities hold
\begin{align*}
\tilde E_{b_i,-b_j} f_{\mbfb,c}=& b_i f_{\mbfb+\mbfe^i-\mbfe^j,c} \comma\\
\tilde E_{b_i,c} f_{\mbfb,c}=& b_i f_{\mbfb+\mbfe^i,c+1} \comma \\
\tilde E_{-b_i,-c} f_{\mbfm,c}=&(1-c)f_{\mbfb-\mbfe^i,c-1} \fstop
\end{align*}

\begin{proof} A proof reduces to the following computations.
\begin{align*}
\tilde E_{b_i,-b_j} f_{\mbfb,c}=& \mbfu^\mbfb t^c \quadre{\sum_{\mbfm\geq \zero} \frac{\Poch{\mbfb}{\mbfm} (m_i+b_i) \mbfx^\mbfm}{\Poch{c}{\length{\mbfm}} \mbfm!} - \sum_{\mbfm\geq \zero} \frac{\Poch{\mbfb}{\mbfm} m_i \mbfx^{\mbfm-\mbfe^i+\mbfe^j}}{\Poch{c}{\length{\mbfm}} \mbfm!} }\\
=& \mbfu^\mbfb t^c \quadre{\sum_{\mbfm\geq \zero} \frac{\Poch{\mbfb}{\mbfm} (m_i+b_i) \mbfx^\mbfm}{\Poch{c}{\length{\mbfm}} \mbfm!} - \sum_{\mbfm\geq \zero} \frac{\Poch{\mbfb}{\mbfm+\mbfe^i-\mbfe^j} (m_i+1)\mbfx^{\mbfm}}{\Poch{c}{\length{\mbfm}} (\mbfm+\mbfe^i-\mbfe^j)!} }\\
=& \mbfu^\mbfb t^c \frac{b_i}{b_j-1} \times\\
&\times\quadre{\sum_{\mbfm\geq \zero} \frac{\Poch{\mbfb+\mbfe^i-\mbfe^j}{\mbfm-\mbfe^i+\mbfe^j} (m_i+b_i) \mbfx^\mbfm}{\Poch{c}{\length{\mbfm}} \mbfm!} - \sum_{\mbfm\geq \zero} \frac{\Poch{\mbfb+\mbfe^i-\mbfe^j}{\mbfm} \mbfx^{\mbfm}}{\Poch{c}{\length{\mbfm}} (\mbfm-\mbfe^j)!} }\\
=& \mbfu^\mbfb t^c \frac{b_i}{b_j-1}\times\\
&\times \quadre{\sum_{\mbfm\geq \zero} \frac{\Poch{\mbfb+\mbfe^i-\mbfe^j}{\mbfm} (m_j+b_j-1) \mbfx^\mbfm}{\Poch{c}{\length{\mbfm}} \mbfm!} - \sum_{\mbfm\geq \zero} \frac{\Poch{\mbfb+\mbfe^i-\mbfe^j}{\mbfm} m_j \mbfx^{\mbfm}}{\Poch{c}{\length{\mbfm}} \mbfm!} }\\
=&b_i f_{\mbfb+\mbfe^i-\mbfe^j,c} \comma
\\
%
\tilde E_{b_i,c} f_{\mbfb,c}=& \mbfu^\mbfb t^c\quadre{\sum_{\mbfm\geq \zero} \frac{\Poch{\mbfb}{\mbfm} (m_i+b_i) \mbfx^\mbfm}{\Poch{c}{\length{\mbfm}} \mbfm!} - \sum_{\mbfm\geq \zero} \frac{\Poch{\mbfb}{\mbfm} m_i(\length{\mbfm}-1) \mbfx^{\mbfm-\mbfe^i}}{\Poch{c}{\length{\mbfm}} \mbfm!} }\\
=&\mbfu^\mbfb t^c\quadre{\sum_{\mbfm\geq \zero} \frac{\Poch{\mbfb}{\mbfm+\mbfe^i} \mbfx^\mbfm}{\Poch{c}{\length{\mbfm}} \mbfm!} - \sum_{\mbfm\geq \zero} \frac{\Poch{\mbfb}{\mbfm} (\length{\mbfm}-1) \mbfx^{\mbfm-\mbfe^i}}{\Poch{c}{\length{\mbfm}} (\mbfm-\mbfe^i)!} }\\
=&\mbfu^\mbfb t^c\quadre{\sum_{\mbfm\geq \zero} \frac{\Poch{\mbfb}{\mbfm+\mbfe^i} \mbfx^\mbfm}{\Poch{c}{\length{\mbfm}} \mbfm!} - \sum_{\mbfm\geq \zero} \frac{\Poch{\mbfb}{\mbfm+\mbfe^i} \length{\mbfm} \mbfx^{\mbfm}}{\Poch{c}{\length{\mbfm}+1} \mbfm!} }\\
=&\mbfu^\mbfb t^c\frac{b_i}{c}\quadre{\sum_{\mbfm\geq \zero} \frac{\Poch{\mbfb+\mbfe^i}{\mbfm} \mbfx^\mbfm}{\Poch{c+1}{\length{\mbfm}-1} \mbfm!} - \sum_{\mbfm\geq \zero} \frac{\Poch{\mbfb+\mbfe^i}{\mbfm} \length{\mbfm} \mbfx^{\mbfm}}{\Poch{c+1}{\length{\mbfm}} \mbfm!} }\\
=&\mbfu^\mbfb t^c\frac{b_i}{c}\quadre{\sum_{\mbfm\geq \zero} \frac{\Poch{\mbfb+\mbfe^i}{\mbfm} \mbfx^\mbfm (c+\length{\mbfm})}{\Poch{c+1}{\length{\mbfm}} \mbfm!} - \sum_{\mbfm\geq \zero} \frac{\Poch{\mbfb+\mbfe^i}{\mbfm} \length{\mbfm} \mbfx^{\mbfm}}{\Poch{c+1}{\length{\mbfm}} \mbfm!} }\\
=&b_i f_{\mbfb+\mbfe^i,c+1} \comma
\\
%
\tilde E_{-b_i,-c} f_{\mbfm,c}=& \mbfu^\mbfb t^c\quadre{\sum_{\mbfm\geq \zero} \frac{\Poch{\mbfb}{\mbfm} \mbfx^{\mbfm+\mbfe^i}}{\Poch{c}{\length{\mbfm}} \mbfm!} -\sum_{\mbfm\geq \zero} \frac{\Poch{\mbfb}{\mbfm} \mbfx^{\mbfm}}{\Poch{c}{\length{\mbfm}} \mbfm!}(\length{\mbfm}+c-1) }\\
=& \mbfu^\mbfb t^c \quadre{\sum_{\mbfm\geq \zero} \frac{\Poch{\mbfb}{\mbfm-\mbfe^i} m_i \mbfx^\mbfm}{\Poch{c}{\length{\mbfm}-1} \mbfm!} -\sum_{\mbfm\geq \zero} \frac{\Poch{\mbfb}{\mbfm} \mbfx^{\mbfm}}{\Poch{c}{\length{\mbfm}} \mbfm!} (\length{\mbfm}+c-1) }\\
=& \mbfu^\mbfb t^c \frac{c-1}{b_i-1}\times\\
&\times\quadre{\sum_{\mbfm\geq \zero} \frac{\Poch{\mbfb-\mbfe^i}{\mbfm} m_i \mbfx^\mbfm}{\Poch{c-1}{\length{\mbfm}} \mbfm!} -\sum_{\mbfm\geq \zero} \frac{\Poch{\mbfb-\mbfe^i}{\mbfm} (b_i+m_i-1) \mbfx^{\mbfm}}{\Poch{c-1}{\length{\mbfm}} (c+\length{\mbfm}-1) \mbfm!} (\length{\mbfm}+c-1) }\\
=&(1-c)f_{\mbfb-\mbfe^i,c-1} \fstop \qedhere
\end{align*}
\end{proof}
\end{lem}

Denote now by~${}_k\tilde\mfg$ the complex linear span of the vectors~\eqref{eq:DiffOpPhi} endowed with the Lie bracket induced by their composition. Having set
\begin{align*}
\set{\tilde J_+,\tilde J_-,\tilde J_0}\eqdef \set{\tilde E_{b_i,-b_j}, \tilde E_{-b_i,b_j}, \tilde J_{b_i,c}-\tilde J_{b_j,c}}, \set{\tilde E_{b_i,c}, \tilde E_{-b_i,-c}, \tilde J_{b_i,c}} \comma
\end{align*}
one can straightforwardly check that~${}_k\tilde\mfg$ is the $k(k+2)$-dimensional dynamical symmetry algebra of the function ${}_k\Phi_2[\mbfb;\length{\mbfb};\mbfx]$, with  commutation relations
\begin{align}\label{eq:Commu}
[\tilde J_+,\tilde J_-]=\tilde J_0 \comma [\tilde J_\pm, \tilde J_0]=\pm 2\tilde J_\pm \comma [\tilde E_{\pm b_i,\mp b_j}, \tilde E_{\mp b_i, \mp c}]=\tilde E_{\mp b_j,\mp c} \fstop
\end{align}

\begin{thm}\label{t:dsa}
The dynamical symmetry algebra of~${}_k\Phi_2[\mbfb;\length{\mbfb};\mbfs]$ satisfies~${}_k\tilde\mfg\cong \mfs\mfl_{k+1}(\C)$.
\begin{proof}
For~$k=1$ (with ${}_1\Phi_2\eqdef {}_2F_1$) the statement is readily seen to hold by inspection, hence by induction on~$k$. Argue by contradiction that there exists~$h>1$ such that~${}_h\tilde \mfg$ is simple yet~${}_{h+1}\tilde \mfg$ is not simple. Thus, there exists~$(0)\lneq\mfi \lhd {}_{h+1}\tilde \mfg$ a nontrivial ideal of~${}_{h+1}\tilde \mfg$. It holds,
\begin{align*}
{}_{h+1}\tilde \mfg={}_{h}\tilde \mfg\oplus \Span_\C \set{\tilde E_{\pm b_{h+1}, \pm c}, \tilde E_{\pm b_{h+1}, \mp b_j}, \tilde J_{b_{h+1},c}\mid j\in[h] }
\end{align*}
as vector spaces, hence~$\mfi$ contains at least one of the basis vectors~$\tilde E_{\pm b_{h+1}, \pm c}, \tilde E_{\pm b_{h+1}, \mp b_j}, \tilde J_{b_{h+1},c}$. By the commutation relations~\eqref{eq:Commu}, and since $\mfi$ is an ideal,
\begin{align*}
(0)\subsetneq [{}_{h+1}\tilde \mfg,[{}_{h+1}\tilde \mfg,\mfi]]\cap {}_{h}\tilde \mfg \subset \mfi\cap {}_h\tilde\mfg\comma
\end{align*}
yet~$\mfi\cap {}_h\tilde\mfg=(0)$ by simplicity of~${}_h\tilde \mfg$, a contradiction. 
Thus~${}_k\tilde\mfg\cong \mfs\mfl_{k+1}(\C)$ by simplicity and dimension.
\end{proof}
\end{thm}

Finally, we are now able to characterize the action of the symmetric group in the definition of quasi-exchangeability for~$\Dir[\boldalpha]$ as a subgroup of the Weyl group of the dynamical symmetry algebra of the second Humbert function (for other interpretations see Rem.~\ref{r:Models}).

\begin{prop}[Another model for~$\mfS_k$]\label{p:chased}
The symmetric group~$\mfS_k$ acting on the standard simplex~$\Delta^{k-1}$ in the definition of quasi-exchangeability for~$\Dir[\boldalpha]$ acts as a unique subgroup of the Weyl group of~${}_k\tilde\mfg$ (described in the proof).

\begin{proof}
Let $\mfh<\mfs\mfl_{k+1}(\C)$ be the sub-algebra of diagonal matrices with vanishing trace and denote by~$W\cong\mfS_{k+1}$ the Weyl group of the root system~$\Psi$ of $\mfs\mfl_{k+1}(\C)$ induced by~$\mfh$.
The action of~$W$ on~$\Psi\subset \mfh^*$ may be regarded as dual to the action of~$\mfS_{k+1}$ on~$\mfh$ via conjugation by permutation matrices in~$\mfP_{k+1}\cong\mfS_{k+1}<GL_{k+1}(\C)$.

For the dynamical symmetry algebra of~${}_k\Phi_2[\mbfb;\length{\mbfb};\mbfx]$, a Lie isomorphism as in Theorem~\ref{t:dsa}, say $\varphi$, may be given in such a way that
\begin{align*}
\varphi\colon \tilde{J}_{b_i,c}\mapsto \diag (1,0,\dotsc,{}_{i-1}0,{}_i -1, {}_{i+1}0, \dotsc,{}_{k+1}0)\in \mfh\fstop
\end{align*}

The action of~$\mfS_k$ on Dirichlet distributions on~$\Delta^{k-1}$, given by $\pi.\Dir[\boldalpha]\eqdef\Dir[\boldalpha_\pi]$, corresponds to permuting~$\tilde{J}_{b_i,c}\mapsto \tilde{J}_{b_{\pi(i)},c}$ for $i\in [k]$. Thus, $\mfS_k$ acts as the subgroup $\mfP_k< \mfP_{k+1}$ fixing the first row and column of matrices in~$\mfh$.
\end{proof}
\end{prop}


\thispagestyle{empty}
{\small
\bibliographystyle{plain}
\markboth{References}{References}
\bibliography{bibliography}
}

\end{document}